\documentclass[12pt]{article}

\usepackage{indentfirst}
\usepackage{amsfonts}
\usepackage{amsmath}
\usepackage{amssymb}
\usepackage{amsbsy}
\usepackage[margin=1in]{geometry}

\newtheorem{dfn}{Definition}[section]
\newtheorem{theorem}[dfn]{Theorem}
\newtheorem{lemma}[dfn]{Lemma}
\newtheorem{prop}[dfn]{Proposition}

\newtheorem{corollary}[dfn]{Corollary}
\newtheorem{conj}[dfn]{Conjecture}
\newtheorem{question}[dfn]{Question}
\newenvironment{pf}{\noindent{\bf Proof.}}
{\enspace\vrule height5pt depth0pt width5pt} 

\def\F {{\mathcal F}}
\def\Se {{\mathcal S}}

\def\X {{\mathcal X}}
\def\W {{\mathcal W}}
\def\C {{\mathcal C}}
\def\L {{\mathcal L}}
\def\G {{\mathcal G}}
\def\D {{\mathcal D}}

\def\Y {{\mathcal Y}}
\def\Q {{\mathcal Q}}

\def\HH {{\mathcal H}}

\def\ex {{\rm ex}}

\def\flap {{\rm flap}}
\def\dup {{\rm dup}}
\def\hom {{\rm hom}}
\def\ind {{\rm ind}}

\begin{document}

\title{Homomorphism counts in robustly sparse graphs}
\author{Chun-Hung Liu\thanks{chliu@math.tamu.edu. Partially supported by NSF under Grant No.~DMS-1954054.} \\
\small Department of Mathematics, \\
\small Texas A\&M University,\\
\small College Station, TX 77843-3368, USA \\
}

\maketitle

\begin{abstract}
For a fixed graph $H$ and for arbitrarily large host graphs $G$, the number of homomorphisms from $H$ to $G$ and the number of subgraphs isomorphic to $H$ contained in $G$ have been extensively studied in extremal graph theory and graph limits theory when the host graphs are allowed to be dense.
This paper addresses the case when the host graphs are robustly sparse and proves a general theorem that solves a number of open questions proposed since 1990s and strengthens a number of results in the literature.

We prove that for any graph $H$ and any set ${\mathcal H}$ of homomorphisms from $H$ to members of a hereditary class ${\mathcal G}$ of graphs, if ${\mathcal H}$ satisfies a natural and mild condition, and contracting disjoint subgraphs of radius $O(\lvert V(H) \rvert)$ in members of ${\mathcal G}$ cannot create a graph with large edge-density, then an obvious lower bound for the size of ${\mathcal H}$ gives a good estimation for the size of ${\mathcal H}$.
This result determines the maximum number of $H$-homomorphisms, the maximum number of $H$-subgraphs, and the maximum number $H$-induced subgraphs in graphs in any hereditary class with bounded expansion up to a constant factor; it also determines the exact value of the asymptotic logarithmic density for $H$-homomorphisms, $H$-subgraphs and $H$-induced subgraphs in graphs in any hereditary nowhere dense class.
Hereditary classes with bounded expansion include (topological) minor-closed families and many classes of graphs with certain geometric properties; nowhere dense classes are the most general sparse classes in sparsity theory.
Our machinery also allows us to determine the maximum number of $H$-subgraphs in the class of all $d$-degenerate graphs with any fixed $d$. 
\end{abstract}

\section{Introduction} \label{sec:intro}

A {\it homomorphism} from a graph\footnote{All graphs are finite and simple in this paper.} $H$ to a graph $G$ is a function $\phi: V(H) \rightarrow V(G)$ such that $\phi(u)\phi(v) \in E(G)$ whenever $uv \in E(H)$.
We denote by $\hom(H,G)$ the number of homomorphisms from $H$ to $G$.
Counting homomorphisms between two graphs appears in many areas, such as extremal graph theory, statistical physics, and computer science.
For example, Lov\'{a}sz \cite{l_iso} showed that if $\hom(H,G_1)=\hom(H,G_2)$ for every graph $H$, then $G_1$ and $G_2$ are isomorphic.
One way to define the limit of a sequence $(G_n)_{n \geq 1}$ of graphs is to consider the homomorphism density $\hom(H,G_n)/\lvert V(G_n) \rvert^{\lvert V(H) \rvert}$ for all graphs $H$.
This notion leads to a rich theory for dense host graphs $G_n$ (see \cite{bclsv_survey,l_book}).
But this definition for homomorphism density does not work well for sparse graphs because $\hom(H,G)$ is often $o(\lvert V(G) \rvert^{\lvert V(H) \rvert})$ for sparse graphs $G$.
A class of graphs of bounded maximum degree is an example of classes of sparse graphs.
If $G$ has bounded maximum degree and $H$ is connected, it is easy to show $\hom(H,G)=O(\lvert V(G) \rvert)$.
So one way to define convergence of a sequence $(G_n)_{n \geq 1}$ of graphs with uniformly bounded maximum degree is to consider the convergence of $\hom(H,G_n)/\lvert V(G_n) \rvert$ for all connected graphs $H$.
This notion of convergence coincides with a notion introduced implicitly by Aldous \cite{a} and explicitly by Benjamini and Schramm \cite{bs} (see \cite{bclsv_survey}).
Even though limit theory was developed for various sparse graph classes, such as for planar graphs \cite{gn}, graphs of bounded average degree \cite{l_limit_bd_av}, graphs of bounded tree-depth \cite{no_limit}, and nowhere dense graphs \cite{no_limit_nd}, it remained unknown what the correct exponent $k$ should be in order to make $\hom(H,G_n)/\lvert V(G_n) \rvert^{k}$ meaningful when $G_n$ is sparse.
This is one motivation of this paper.

A central direction in extremal graph theory is solving Tur\'{a}n-type questions which ask for the maximum number of copies of a fixed graph in graphs in a fixed graph class.
Tur\'{a}n \cite{t_turan} determined the maximum number of edges in $K_t$-free graphs.
Erd\H{o}s and Stone \cite{es} generalized it by determining the maximum number of edges in $n$-vertex $L$-free graphs, for any fixed graph $L$, up to an $o(n^2)$ error. 
In fact, the number of edges equals the number of subgraphs isomorphic to $K_2$.
Zykov \cite{z} determined the maximum number of subgraphs isomorphic to $K_s$ in $K_t$-free graphs, for any fixed integers $s$ and $t$.
In general, the maximum number of subgraphs isomorphic to a fixed graph $H$ in $L$-free graphs for another fixed graph $L$ has been extensively studied.
For example, the cases when $H$ and $L$ are cycles, trees, complete graphs, or complete bipartite graphs were considered in \cite{as,ggmz,g,hhknr}.
It is far from a complete list of known results of this type.
We refer readers to \cite{fm} for a survey.

Tur\'{a}n-type questions are closely related to counting homomorphisms.
Counting the number of edges is equivalent to counting the number of $K_2$-subgraphs, and the number of $K_s$-subgraphs in $G$ equals $\hom(K_s,G)/s!$.
More generally, the number of $H$-subgraphs in $G$ equals the number of injective homomorphisms from $H$ to $G$ divided by the size of the automorphism group of $H$.
In addition, it is known that the difference between $\hom(H,G)$ and the number of injective homomorphisms from $H$ to $G$ is at most $O(\lvert V(G) \rvert^{\lvert V(H) \rvert-1})$. 
So the homomorphism density approximates the density of injective homomorphisms with $o(1)$ additive errors in dense graphs. 
This allows a powerful machinery for the use of homomorphism inequalities to solve Tur\'{a}n-type problems in dense graphs \cite{r}.
On the other hand, it was unclear how a similar machinery can be applied when the host graphs $G$ are sparse because $\hom(H,G)=o(\lvert V(G) \rvert^{\vert V(H) \rvert})$, and it was unknown what the correct exponent $k$ to normalize homomorphism counts to properly define homomorphism densities $\hom(H,G)/\lvert V(G) \rvert^k$ is.
This leads to another motivation of this paper.

More precisely, we are interested in for every fixed graph $H$, determining the value $k$ such that $\max_{G \in \G, \lvert V(G) \rvert=n}\hom(H,G)=\Theta(n^k)$ and the analogous problem for the number of $H$-subgraphs in $G$, where $\G$ is a class of sparse graphs.
For some graphs $H$ and some classes $\G$ of sparse graphs, it is not hard to determine this value $k$; in those cases, the coefficient of $n^k$ has been studied, such as when $H$ is a complete graph and $G$ is either of bounded maximum degree \cite{c,cr}, of bounded Euler genus \cite{dfjsw}, $K_t$-minor free \cite{fw_clique_minor}, or $K_t$-topological minor free \cite{fw_clique_subdiv,lo}. 
The subgraph counts when $G$ is planar and $H$ is either a cycle \cite{gpstz_c5,hs}, a complete bipartite graph \cite{ac}, or a 4-vertex path \cite{gpstz_p4} were also studied.

However, determining the exponent $k$ for any fixed graph $H$ and fixed class of sparse graphs seems challenging in general.
For example, the following old question of Eppstein \cite{e} remained open (until this paper), where only some special cases for trees \cite{hw} or for classes of graphs embeddable in surfaces \cite{e,hjw,w_3conn} were known; in fact, it was even unknown whether the value $k$ in the following question can be chosen to be independent with $n$. 

\begin{question}[\cite{e}] \label{que_minor_poly}
Are there proper minor-closed families $\F$ and graphs $H$ such that the maximum number of sugraphs isomorphic to $H$ in an $n$-vertex graph in $\F$ is $\Theta(n^k)$ for some non-integer $k$?
\end{question}

The main result of this paper (Theorem \ref{main}) provides a negative answer to Question \ref{que_minor_poly}.
Theorem \ref{main} is actually much general in the following aspects.
	\begin{itemize}
		\item Theorem \ref{main} determines the precise value of this integer $k$ for any hereditary class with bounded expansion (Theorem \ref{bdd_expan_intro}) and determines $k$ up to an additive $o(1)$ error term for any hereditary nowhere dense class (Theorem \ref{nowhere_dense_intro}).
			Nowhere dense classes are the most general sparse graph classes in sparsity theory; classes with bounded expansion are common generalizations of (topological) minor-closed families and many extensively studied graph classes with some geometric properties.
			Further motivation and formal definition of those classes are included in Section \ref{subsec:main_results}.
			(See Figure \ref{fig_sparse} for a relationship between classes of sparse graphs, and see Section \ref{sec:concrete} for concrete examples of applications of Theorem \ref{bdd_expan_intro} that solves other open questions in the literature.)	
			\begin{figure} \label{fig_sparse}
				\begin{picture}(100,180) (0,30)
					\thicklines

					\put(0,200){Bounded}
					\put(0,185){tree-depth}
					\put(-3,215){\line(1,0){59}}
					\put(56,215){\line(0,-1){38}}
					\put(56,177){\line(-1,0){59}}
					\put(-3,177){\line(0,1){38}}

					\put(56,196){\vector(1,0){43}}

					\put(0,112){Trees}
					\put(-3,127){\line(1,0){33}}
					\put(30,127){\line(0,-1){22}}
					\put(30,105){\line(-1,0){33}}
					\put(-3,105){\line(0,1){22}}

					\put(30,116){\vector(1,0){37}}
					\put(30,127){\vector(3,2){73}}

					\put(102,200){Bounded}
					\put(102,185){tree-width}
					\put(99,215){\line(1,0){59}}
					\put(158,215){\line(0,-1){38}}
					\put(158,177){\line(-1,0){59}}
					\put(99,177){\line(0,1){38}}

					\put(158,196){\vector(1,0){39}}

					\put(70,112){Planar}
					\put(67,127){\line(1,0){40}}
					\put(107,127){\line(0,-1){22}}
					\put(107,105){\line(-1,0){40}}
					\put(67,105){\line(0,1){22}}
					
					\put(107,116){\vector(1,0){35}}

					\put(200,200){Bounded}
					\put(200,185){layered tree-width}
					\put(197,215){\line(1,0){99}}
					\put(296,215){\line(0,-1){38}}
					\put(296,177){\line(-1,0){99}}
					\put(197,177){\line(0,1){38}}

					\put(296,196){\vector(1,0){43}}

					\put(145,120){Bounded}
					\put(145,105){Euler genus}
					\put(142,135){\line(1,0){65}}
					\put(207,135){\line(0,-1){38}}
					\put(207,97){\line(-1,0){65}}
					\put(142,97){\line(0,1){38}}

					\put(207,116){\vector(1,0){35}}
					\put(190,135){\vector(1,1){41}}

					\put(245,120){Minor-}
					\put(245,105){closed}
					\put(242,135){\line(1,0){42}}
					\put(284,135){\line(0,-1){38}}
					\put(284,97){\line(-1,0){42}}
					\put(242,97){\line(0,1){38}}

					\put(284,116){\vector(1,0){33}}
					\put(284,135){\vector(2,1){55}}

					\put(50,50){Bounded}
					\put(50,35){maximum degree}
					\put(47,65){\line(1,0){93}}
					\put(140,65){\line(0,-1){38}}
					\put(140,27){\line(-1,0){93}}
					\put(47,27){\line(0,1){38}}

					\put(140,46){\vector(1,0){50}}

					\put(193,42){Immersion-closed}
					\put(190,57){\line(1,0){94}}
					\put(284,57){\line(0,-1){22}}
					\put(284,35){\line(-1,0){94}}
					\put(190,35){\line(0,1){22}}

					\put(284,57){\vector(1,1){39}}

					\put(342,200){Admitting}
					\put(342,185){strongly sublinear}
					\put(342,170){balanced separators}
					\put(339,215){\line(1,0){109}}
					\put(448,215){\line(0,-1){53}}
					\put(448,162){\line(-1,0){109}}
					\put(339,162){\line(0,1){53}}

					\put(420,162){\vector(1,-1){26}}

					\put(320,120){Topological}
					\put(320,105){minor-closed}
					\put(317,135){\line(1,0){70}}
					\put(387,135){\line(0,-1){38}}
					\put(387,97){\line(-1,0){70}}
					\put(317,97){\line(0,1){38}}

					\put(387,116){\vector(1,0){30}}

					\put(420,120){Bounded}
					\put(420,105){expansion}
					\put(417,135){\line(1,0){56}}
					\put(473,135){\line(0,-1){38}}
					\put(473,97){\line(-1,0){56}}
					\put(417,97){\line(0,1){38}}

					\put(450,97){\vector(0,-1){32}}

					\put(425,50){Nowhere}
					\put(425,35){dense}
					\put(422,65){\line(1,0){50}}
					\put(472,65){\line(0,-1){38}}
					\put(472,27){\line(-1,0){50}}
					\put(422,27){\line(0,1){38}}
				\end{picture}
				\caption{Relationship between extensively studied classes of sparse graphs. ``A $\rightarrow$ B'' means that any class satisfying property A satisfies property B.}
			\end{figure}
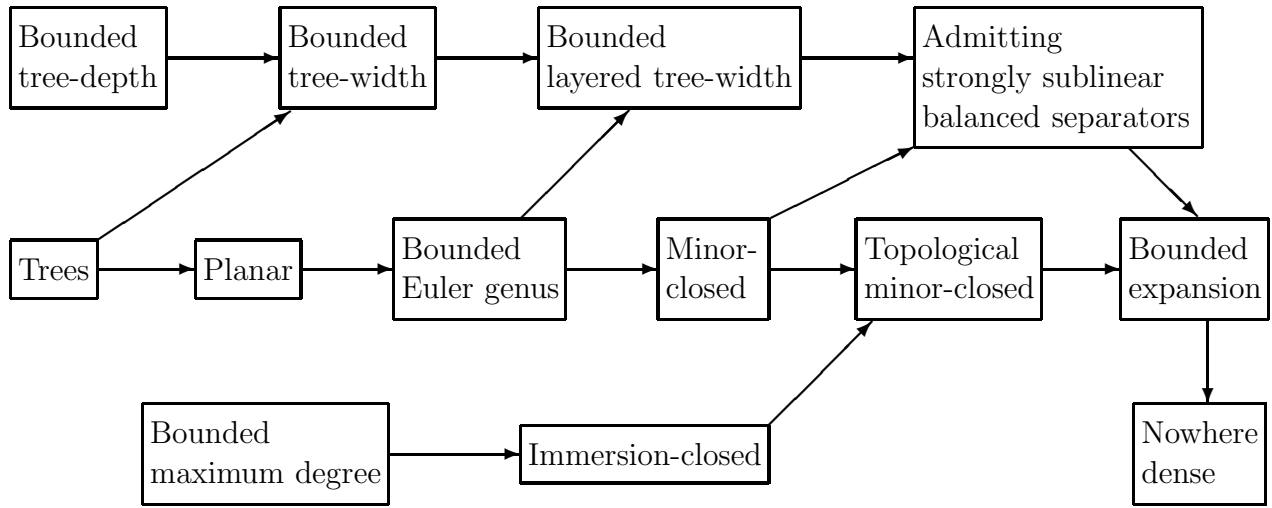
		\item Moreover, this value $k$ equals an ``obvious lower bound'' obtained by duplication and can be decided in finite time (if an oracle for membership testing is given) and can be explicitly stated (for some classes that are explicitly described, see Section \ref{sec:concrete}).
			See Section \ref{sec:obvious_lower_bdd} for the precise definition of this obvious lower bound.\footnote{We remark that we state Theorems \ref{bdd_expan_intro} and \ref{nowhere_dense_intro} in a general setting that will be described in the following bullets.
			The general setting requires a number of definitions stated in Sections \ref{subsec:label}-\ref{subsec:lower_bdd_general}.
			Readers who are interesting in counting subgraphs only can skip Sections \ref{subsec:label}-\ref{subsec:lower_bdd_general} and jump to Section \ref{subsec:main_results} after reading Proposition \ref{lower_bound_easy}, and can treat all labelled graphs as graphs and treat $\ex(H,f_H,\HH,\G,n)$ and $\dup_\HH(H,f_H,\G)$ stated in Theorems \ref{bdd_expan_intro} and \ref{nowhere_dense_intro} as $\ex(H,\G,n)$ and $\max_\C\lvert \C \rvert$ over all collections $\C$ stated in Proposition \ref{lower_bound_easy}, respectively.
			All applications of Theorems \ref{bdd_expan_intro} and \ref{nowhere_dense_intro} explicitly stated in this paper only rely on this special case for counting subgraphs in unlabelled graphs.} 
		\item Theorem \ref{main} counts the size of any set $\HH$ of homomorphisms from any fixed graph $H$ to graphs in any fixed class of sparse graphs mentioned above, as long as $\HH$ satisfies some natural and mild conditions for consistency with respect to isomorphisms. 
			(See Section \ref{subsec:consistent} for the precise definition of consistent sets of homomorphisms.)
			This allows us to count the number of $H$-homomorphisms (without other restrictions), the number of $H$-subgraphs, the number of induced $H$-subgraphs, the number of $H$-homomorphisms in which every vertex in the image is mapped by a bounded number of vertices, and more.
		\item Theorem \ref{main} works for graphs whose some cliques are labelled.
			(See Section \ref{subsec:label} for the formal definition of the labelling.)
			This allows us to count objects beyond graphs, such as directed graphs, edge-colored graphs, hypergraphs, and relational structures.
		\item For any fixed graph $H$, Theorem \ref{main} can count the number of homomorphisms from $H$ to graphs in a class $\F$, as long as any shallow minor in graphs in $\F$ with depth $O(\lvert V(H) \rvert)$ has low edge-density, in contrast to nowhere dense classes or bounded expansion classes that require low edge-density for shallow minors with all depths.
			Our proof does not use machineries that are frequently applied in the study of sparsity theory such as low tree-depth coloring, center coloring or weak coloring.
			Our machinery applies to $d$-degenerate graphs for any fixed integer $d$ (Theorem \ref{bdd_degen_intro}), which is a somewhere dense class of graphs, and solves a conjecture of Huynh and Wood in \cite{hw}.
	\end{itemize}

\subsection{An obvious lower bound} \label{sec:obvious_lower_bdd}

We describe an obvious lower bound for the number of homomorphisms from a fixed graph to graphs in a graph class in this subsection. 

A {\it separation} of a graph $G$ is an ordered pair $(A,B)$ of subsets of $V(G)$ such that $A \cup B=V(G)$ and there exists no edge of $G$ between $A-B$ and $B-A$.
The {\it order} of $(A,B)$ is $\lvert A \cap B \rvert$.
A collection $\C$ of separations of $G$ is {\it independent} if 
	\begin{itemize}
		\item for every member $(A,B)$ of $\C$, $A-B \neq \emptyset$, and 
		\item for any distinct members $(A,B)$ and $(C,D)$ of $\C$, we have $A \subseteq D$ and $C \subseteq B$. 
	\end{itemize}
We remark that we do not require $B \neq \emptyset$ for members $(A,B)$ in an independent collection, so $\{(V(G),\emptyset)\}$ is an independent collection.
Note that for any two distinct members $(A,B)$ and $(C,D)$ of an independent collection, $A-B \subseteq D-C$, so $A-B$ and $C-D$ are disjoint, and there exists no edge between $A-B$ and $C-D$.
Hence every independent collection of separations of $H$ has size at most the independence number of $H$.
(The {\it independence number} of $H$, denoted by  $\alpha(H)$, is the maximum size of a set of pairwise non-adjacent vertices in $H$.)

Let $H$ be a graph.
Let $Z \subseteq V(H)$.
Let $k$ be a positive integer.
We define $H \wedge_k Z$ to be the graph obtained from a union of $k$ disjoint copies of $H$ by for each $z \in Z$, identifying the $k$ copies of $z$ into a vertex.
So $\lvert V(H \wedge_k Z) \rvert = k(\lvert V(H) \rvert-\lvert Z \rvert)+\lvert Z \rvert$.
Let $\C$ be an independent collection of separations of a graph $H$.
We define $H \wedge_k \C$ to be $H \wedge_k (\bigcap_{(A,B) \in \C}B)$. 
It is easy to see that if $G=H \wedge_k \C$, then there are at least $k^{\lvert \C \rvert} \geq (\frac{\lvert V(G) \rvert}{\lvert V(H) \rvert})^{\lvert \C \rvert} \geq \frac{1}{\lvert V(H) \rvert^{\alpha(H)}}\lvert V(G) \rvert^{\lvert \C \rvert}$ induced subgraphs of $G$ isomorphic to $H$.
Hence we obtain the following obvious lower bound. 

\begin{prop} \label{lower_bound_easy}
For every graph $H$, there exists a real number $c=c(H)$ such that if $\G$ is a class of graphs and $\C$ is an independent collection of separations of $H$ such that $H \wedge_\ell \C \in \G$ for infinitely many positive integers $\ell$, then $\ex(H,\G,n)  \geq cn^{\lvert \C \rvert}$ for infinitely many integers $n$, where $\ex(H,\G,n)$ is the maximum number of subgraphs isomorphic to $H$ contained in an $n$-vertex graph in $\G$.
\end{prop}

This obvious lower bound actually works for a more general setting (Proposition \ref{lower_bound}).
We need a number of definitions to state our general setting and Proposition \ref{lower_bound}.

\subsubsection{Labelled graphs} \label{subsec:label}

Let $G$ be a graph.
A {\it march} in $G$ is a sequence over $V(G)$ with distinct entries.
A {\it labelling} of $G$ is a function such that every element in its domain is a march in $G$.
A labelling $f$ of $G$ is {\it legal} if for every march in the domain of $f$, its entries form a clique.

A {\it labelled graph} is a pair $(G,f_G)$, where $G$ is a graph and $f_G$ is a labelling of $G$.
A labelled graph is {\it legal} if its labelling is legal.

A {\it quasi-order} is a pair $(X,\preceq)$, where $X$ is a set and $\preceq$ is a reflexive and transitive binary relation on $X$.
For a quasi-order $Q=(X,\preceq)$, a {\it $Q$-labelling} of a graph $G$ is a labelling of $G$ whose image is a subset of $X$, and a {\it $Q$-labelled graph} is a labelled graph whose labelling is a $Q$-labelling.

Each of following objects can be encoded as legal $Q$-labelled graphs for an antichain $Q$: 
	\begin{itemize}
		\item Directed graphs. (Map each march corresponding to a directed edge to an element in $Q$ representing the direction of this edge.)
		\item Edge-colored graphs (such as in the setting for $k$-common graphs \cite{cy,jst,knnvw} or for rainbow Tur\'{a}n problems \cite{gmmp,kmsv}).
		\item $(m,n)$-colored mixed graphs (such as in the setting in \cite{nr}).
		\item Hypergraphs. (Construct a graph by making each hyperedge a clique and label this clique). 
		\item Relational structures. (Consider their Gaifman graphs, see \cite{no_book}.) 
	\end{itemize}

For any function $f$ and a sequence $s=(s_1,s_2,...,s_n)$ whose every entry is in the domain of $f$, we define $f(s)$ to be the sequence $(f(s_1),f(s_2),...,f(s_n))$.
Similarly, if $S$ is a subset of the domain of $f$, we define $f(S)$ to be the set $\{f(x):x \in S\}$.

Let $Q$ be a quasi-order.
An {\it isomorphism} from a $Q$-labelled graph $(H,f_H)$ to a $Q$-labelled graph $(G,f_G)$ is an isomorphism\footnote{An {\it isomorphism} from an (unlabelled) graph $H$ to an (unlabelled) graph $G$ is a function $\iota: V(H) \rightarrow V(G)$ such that $\iota$ is a bijection, and for every $u,v \in V(H)$, $uv \in E(H)$ if and only if $\iota(u)\iota(v) \in E(G)$.} $\phi$ from $H$ to $G$ such that for every march $m$ in the domain of $f_H$, $\phi(m)$ is in the domain of $f_G$ with $f_H(m)=f_G(\phi(m))$, and for every march $m'$ in the domain of $f_G$, $\phi^{-1}(m')$ is in the domain of $f_H$.
Two $Q$-labelled graphs are {\it isomorphic} if there is an isomorphism between them.

Let $Q$ be a quasi-order.
A $Q$-labelled graph $(H,f_H)$ is an {\it induced $Q$-labelled subgraph} (or simply {\it induced subgraph}) of a $Q$-labelled graph $(G,f_G)$ if $H$ is an induced subgraph of $G$, the domain of $f_H$ is a subset of the domain of $f_G$, and for each march $m$ in the domain of $f_G$ whose all entries are in $V(H)$, $m$ is in the domain of $f_H$ and $f_H(m)=f_G(m)$.
We say that a class $\G$ of $Q$-labelled graphs is {\it hereditary} if every induced $Q$-labelled subgraph of a member of $\G$ is a member of $\G$.
A $Q$-labelled graph $(H',f_{H'})$ is a {\it $Q$-labelled subgraph} (or simply {\it subgraph}) of a $Q$-labelled graph $(G,f_G)$ if $H'$ is a subgraph of $G$, the domain of $f_{H'}$ is a subset of the domain of $f_G$, and for every march $m$ in the domain of $f_{H'}$, $f_{H'}(m)=f_G(m)$.
And we say that $(H',f_{H'})$ is a {\it spanning subgraph} of $(G,f_G)$ if $(H',f_{H'})$ is a subgraph of $(G,f_G)$ and $V(H')=V(G)$.

Let $Q=(X,\preceq)$ be a quasi-order.
A {\it homomorphism} from a legal $Q$-labelled graph $(H,f_H)$ to a $Q$-labelled graph $(G,f_G)$ is a homomorphism $\phi$ from $H$ to $G$ such that for every march $m$ in the domain of $f_H$, $\phi(m)$ is in the domain of $f_G$ and $f_H(m) \preceq f_G(\phi(m))$.
Note that since $(H,f_H)$ is legal, every march $m$ in the domain of $f_H$ form a clique, so entries of $\phi(m)$ are distinct and hence $\phi(m)$ is a march.

For simplicity, we say that $\HH$ is a {\it set of homomorphisms from $(H,f_H)$ to a class $\G$} if every member of $\HH$ is a homomorphism from $(H,f_H)$ to a member of $\G$.

A {\it separation} of a labelled graph $(G,f)$ is a separation of $G$.
Note that if $f$ is legal, then for any separation $(A,B)$ of $G$, there exists no march in the domain of $f$ containing entries in both $A-B$ and $B-A$.

\subsubsection{Consistent sets of homomorphisms} \label{subsec:consistent}

Let $Q$ be a quasi-order.
Let $(H,f_H)$ be a legal $Q$-labelled graph, and let $\G$ be a hereditary class of legal $Q$-labelled graphs.
Let $\HH$ be a set of homomorphisms from $(H,f_H)$ to members in $\G$.
We say that $\HH$ is {\it consistent} if the following conditions hold.
	\begin{itemize}
		\item[(CON1)] For every homomorphism $\phi \in \HH$ from $(H,f_H)$ to a member $(G,f)$ in $\G$, 
			\begin{itemize}
				\item if $(G',f')$ is an induced subgraph of $(G,f)$ with $\phi(V(H)) \subseteq V(G')$, then the function obtained from $\phi$ by restricting the codomain to be $V(G')$ is a homomorphism from $(H,f_H)$ to $(G',f')$ and belongs to $\HH$;
				\item if $(G,f)$ is an induced subgraph of another member $(G',f')$ in $\G$, then the function obtained from $\phi$ by changing the codomain to be $V(G')$ is a homomorphism from $(H,f_H)$ to $(G',f')$ and belongs to $\HH$;
			\end{itemize}
		\item[(CON2)] For any members $(G_1,f_1)$ and $(G_2,f_2)$ of $\G$ with an isomorphism $\iota$ from $(G_1,f_1)$ to $(G_2,f_2)$, if $\phi$ is a homomorphism in $\HH$ from $(H,f_H)$ to $(G_1,f_1)$, then $\iota \circ \phi$ is a homomorphism in $\HH$ from $(H,f_H)$ to $(G_2,f_2)$.  
		\item[(CON3)] For any members $(G_1,f_1)$ and $(G_2,f_2)$ of $\G$ such that there exists an isomorphism $\iota$ from $(G_2,f_2)$ to a spanning subgraph $(G_1',f_1')$ of $(G_1,f_1)$, if $\phi$ is a homomorphism in $\HH$ from $(H,f_H)$ to $(G_1,f_1)$ and $\phi$ is also an homomorphism\footnote{Note that we only assume that $\G$ is hereditary, so $(G_1',f_1')$ is not necessarily in $\G$ and hence $\phi$ is not necessarily in $\HH$.} from $(H,f_H)$ to $(G_1',f_1')$, then $\iota^{-1} \circ \phi$ is a homomorphism from $(H,f_H)$ to $(G_2,f_2)$, and $\iota^{-1} \circ \phi \in \HH$. 
	\end{itemize}
Note that (CON1) and (CON2) are simply conditions for ``consistency'' with respect to the induced subgraph relation and isomorphisms, respectively. 
(CON3) is a condition stating that redundant edges can be dropped.
Even though (CON3) might seem artificial at the first glance, it is actually a combination of analogies of (CON1) and (CON2) if $\G$ is assumed to be closed under taking subgraphs (this assumption can always be made if one only considers homomorphisms with no restriction on non-adjacent pairs of vertices); but here we only assume that $\G$ is hereditary (that is, closed under taking induced subgraphs) which it is a weak assumption; this weaker assumption is more natural when we consider homomorphisms that might have restrictions on non-adjacent pairs, such as counting the number of induced subgraphs.
In fact, (CON2) is a special case of (CON3), but we state them separately since we only need (CON2) in many lemmas.

The following sets are examples of consistent sets of homomorphisms:
	\begin{itemize}
		\item The set of all homomorphisms from $(H,f_H)$ to members in $\G$. 
		\item The set of homomorphisms from $(H,f_H)$ to members in $\G$ corresponding to subgraph embeddings (or induced subgraph embeddings, respectively).
		\item For any fixed integer $t$, the set consisting of all homomorphisms from $(H,f_H)$ to members of $\G$ such that for every member $\phi$ of this set and $v \in V(G)$ with $(G,f_G) \in \G$, the size of the set $\{x \in V(H): \phi(x)=v\}$ is at most $t$.
	\end{itemize}

Let $Q$ be a quasi-order.
Let $(H,f_H)$ be a legal $Q$-labelled graph.
Let $\HH$ be a set of consistent homomorphisms from $(H,f_H)$ to a hereditary class $\G$ of legal $Q$-labelled graphs.
For any member $(G,f_G) \in \G$, we define $\ex(H,f_H,\HH,G,f_G)$ to be the number of homomorphisms in $\HH$ from $(H,f_H)$ to $(G,f_G)$.
For every positive integer $n$, we define 
$$\ex(H,f_H,\HH,\G,n) := \max_{(G,f_G) \in \G, \lvert V(G) \rvert=n}\ex(H,f_H,\HH,G,f_G).$$

\subsubsection{Duplication}

For a labelled graph $(G,f_G)$, an {\it independent collection of separations} of $(G,f_G)$ is an independent collection of separations of $G$.

Let $(H,f_H)$ be a labelled graph.
Let $\L$ be an independent collection of separations of $(H,f_H)$.
For a positive integer $w$, we define $(H,f_H) \wedge_w \L$ to be the labelled graph $(G,f_G)$ such that there exist $w$ disjoint isomorphic copies $(G_1,f_{G_1}), (G_2,f_{G_2}),...,(G_w,f_{G_w})$ of $(H,f_H)$ such that
	\begin{itemize}
		\item $G$ is obtained from $G_1,G_2,...,G_w$ by for each vertex $v \in \bigcap_{(A,B) \in \L}B$, identifying the copies of $v$ in $G_1,G_2,...,G_w$ into a vertex $v^G$, 
		\item if $m$ is a march in the domain of $f_{G_i}$ for some $i \in [w]$, then $m^G$ is in the domain of $f_G$ and $f_G(m^G)=f_{G_i}(m)$, where $m^G$ is the march in $G$ obtained from $m$ by replacing each entry $v$ of $m$ with $v \in \bigcap_{(A,B) \in \L}B$ by $v^G$, and
		\item for every march $m$ in the domain of $f_G$, there exist $i \in [w]$ and a march $m_i$ in the domain of $f_{G_i}$ such that $m=m_i^G$ and $f_G(m)=f_{G_i}(m_i)$. 
	\end{itemize}
For a class of labelled graphs $\G$, we say that $\L$ is {\it $\G$-duplicable} if there exist infinitely many positive integers $k$ such that $(H,f_H) \wedge_k \L \in \G$.

For a graph $G$ and a subset $S$ of $V(G)$, we define $N_G[S]$ to be the set of all vertices of $G$ that are either in $S$ or are adjacent in $G$ to some vertex in $S$.

Let $Q$ be a quasi-order.
Let $(H,f_H)$ be a legal $Q$-labelled graph, and let $\L$ be an independent collection of separations of $(H,f_H)$.
Let $\HH$ be a set of homomorphisms from $(H,f_H)$ to a class $\G$ of $Q$-labelled graphs.
We say that $\L$ is {\it $(\G,\HH)$-duplicable} if there exist a $Q$-labelled graph $(G,f_G)$ and a homomorphism $\phi \in \HH$ from $(H,f_H)$ to $(G,f_G)$ such that
	\begin{itemize}
		\item $\phi$ is onto,
		\item $\phi(A-B) \cap \phi(B) = \emptyset$ for any $(A,B) \in \L$, and
		\item $\{(N_G[\phi(A-B)], V(G)-\phi(A-B)): (A,B) \in \L\}$ is a $\G$-duplicable independent collection of separations of $(G,f_G)$ with size $\lvert \L \rvert$. 
	\end{itemize}
Define 
$$\dup_\HH(H,f_H,\G) := \max_\C \lvert \C \rvert,$$ where the maximum is over all $(\G,\HH)$-duplicable independent collections $\C$ of separations of $(H,f_H)$.

\subsubsection{Obvious lower bound in the general form} \label{subsec:lower_bdd_general}

Now we are ready to state the obvious lower bound in the general form.
Its proof is almost identical to the proof of Proposition \ref{lower_bound_easy}, so we omit its proof. 

\begin{prop} \label{lower_bound}
For every graph $H$, there exists a positive real number $c = c(H)$ such that the following holds.
Let $f_H$ be a legal $Q$-labelling for some quasi-order $Q$.
Let $\G$ be a hereditary class of legal $Q$-labelled graphs. 
Let $\HH$ be a consistent set of homomorphisms from $(H,f_H)$ to members of $\G$.
Then there are infinitely many positive integers $n$ such that 
$$\ex(H,f_H,\HH,\G,n) \geq cn^{\dup_\HH(H,f_H,\G)}.$$
\end{prop}

The main results of this paper show that the lower bound in Proposition \ref{lower_bound} is also an upper bound, up to a constant factor or an $n^{o(1)}$ multiplicative error, when $\G$ has bounded expansion or is nowhere dense, respectively.

An immediately following question is whether $\dup_\HH(H,f_H,\G)$ is decidable.
It is not hard to give explicit descriptions for $\dup_\HH(H,f_H,\G)$ when $\G$ has nice structures (see examples in Section \ref{sec:concrete}) or when $H$ has only few vertices or has nice structures (such as for the cases studied in Tur\'{a}n-type questions, like when $H$ is a complete graph, a complete bipartite graph, a forest or a cycle).
In general, deciding $\dup_\HH(H,f_H,\G)$ only requires to check whether each independent collection $\C$ of separations of $(H,f_H)$ is $(\G,\HH)$-duplicable or not; and there are only finitely many such collections.
To check whether a collection $\C$ is $(\G,\HH)$-duplicable or not, since $\G$ is hereditary and $\HH$ is consistent, it suffices to check for each $Q$-labelled graph $(H',f_{H'})$ on at most $\lvert V(H) \rvert$ vertices satisfying that there exists an onto homomorphism in $\HH$ from $(H,f_H)$ to $(H',f_{H'})$, whether the corresponding independent collection of $(H',f_{H'})$ is $\G$-duplicable or not.
Note that there are only finitely many such graphs $H'$, and if $Q$ is finite (such as for the case we encode directed graphs or hypergraphs), then there are only finitely many such $Q$-labelled graphs $(H',f_{H'})$.
Moreover, since $\G$ is hereditary, for every labelled graph $(H,f_H)$, there exists an integer $k$ such that any independent collection $\L$ of separations of $(H,f_H)$ is $\G$-duplicable if and only if $(H,f_H) \wedge_k \L \in \G$.
Hence $\dup_\HH(H,f_H,\G)$ can be decided in finite time if $Q$ is finite, provided an oracle to test the membership of $\G$ is given.

\subsection{Main results} \label{subsec:main_results}

One of the main results of this paper shows that the obvious lower bounded mentioned in Proposition \ref{lower_bound} actually gives an upper bound for graphs with a robust sparsity condition.
We first describe such a condition.

Planar graphs and graphs with bounded degeneracy are sparse in the sense that the number of their edges are at most a linear function of the number of their vertices.
But these two kinds of graphs behave very differently.
Planarity is a robust sparsity in the sense that contracting any connected subgraph preserves the planarity, so contracting disjoint connected subgraphs in planar graphs cannot create a dense graph.
On the other hand, bounded degeneracy is not a robust sparsity, since subdividing every edge of a large complete graph results in a 2-degenerate graph, and contracting disjoint stars from this 2-degenerate graph results in a very dense graph.

Such robust sparsity conditions are the hearts of the extensively studied sparsity theory (for example, see \cite{no_book}) and can be described in terms of conditions of the edge-density of shallow minors.

The {\it radius} of a connected graph $G$ is the minimum $k$ such that there exists a vertex $v$ of $G$ such that every vertex of $G$ has distance from $v$ at most $k$. 
Let $\ell \in \mathbb{Z} \cup \{\infty\}$.
We say that a graph $G$ contains another graph $H$ as an {\it $\ell$-shallow minor} (or equivalently, {\it $H$ is an $\ell$-shallow minor of $G$}) if $H$ is isomorphic to a graph that can be obtained from a subgraph $G'$ of $G$ by contracting disjoint connected subgraphs of $G'$ of radius at most $\ell$.
In other words, every branch set of an $\ell$-shallow minor is a connected subgraph of radius at most $\ell$.
Note that $0$-shallow minors are exactly subgraphs, up to isomorphism.
We say that a graph $H$ is a {\it minor} of $G$ if $G$ contains $H$ as an $\infty$-shallow minor.
Note that this definition for minors is equivalent to the usual definition for minors in the literature.

Theorem \ref{main} is the main theorem of this paper.
It determines $\ex(H,f_H,\HH,\G,n)$ up to a constant factor when the edge-density of every $O(\lvert V(H) \rvert)$-shallow minor in a graph in $\G$ is small.
In order to serve as a source to derive almost all other results in this paper and to solve open questions in the literature, it is stated in a general but technical form.
We postpone its formal statement until Section \ref{sec:homo_density}.
We include some simple applications of Theorem \ref{main} in this section.
More concrete applications of Theorem \ref{main} are included in Section \ref{sec:concrete}.

A class $\F$ of graphs has {\it bounded expansion} if there exists a function $f: {\mathbb N} \cup \{0\} \rightarrow {\mathbb R}$ such that for every nonnegative integer $d$, every $d$-shallow minor of a graph in $\F$ has average degree at most $f(d)$. 

Classes of bounded expansion are very general.
They not only include any class of bounded maximum degree but also include a number of classes of graphs with some natural geometric properties.
(See Figure \ref{fig_sparse} for a relationship between bounded expansion classes and other classes.)
Well-known examples of classes of bounded expansion include any proper minor-closed family \cite{m} (such as the classes of graphs of bounded Euler genus and the class of knotless embeddable graphs), any proper topological minor-closed family \cite{d} (such as the classes of bounded crossing number), any class of graphs with bounded queue number or bounded stack number \cite{now}, and any class of graphs admitting strongly sublinear balanced separators \cite{dn} (such as any class of string graphs with a forbidden bipartite subgraph \cite{fp} and any class of intersection graphs of sets in ${\mathbb R}^d$ with certain geometric properties \cite{dmn,mttv}). 
Moreover, for every $p>0$, there exists a class with bounded expansion such that an Erd\H{o}s-R\'{e}nyi random graph ${\mathbb G}(n,\frac{p}{n})$ belongs to this class asymptotically almost surely \cite{now}.

Our main result (Theorem \ref{main}) determines $\ex(H,f_H,\HH,\G,n)$ up to a constant factor for any class of bounded expansion.
We say that a class of labelled graphs has {\it bounded expansion} if the class of their underlying graphs has bounded expansion.

\begin{theorem} \label{bdd_expan_intro}
For any quasi-order $Q$ with finite ground set, hereditary class $\G$ of legal $Q$-labelled graphs with bounded expansion, legal $Q$-labelled graph $(H,f_H)$, consistent set $\HH$ of homomorphisms from $(H,f_H)$ to $\G$, we have
$$\ex(H,f_H,\HH,\G,n) =\Theta(n^{\dup_\HH(H,f_H,\G)}).$$
\end{theorem}

Theorem \ref{bdd_expan_intro} solves a number of open questions in the literature.
Most of them will be stated in Section \ref{sec:concrete}.
We give two examples in this section.

For any class $\G$ of graphs and graph $H$, we define $\ex(H,\G,n)$ to be the maximum number of subgraphs isomorphic to $H$ contained in $G$, where the maximum is over all $n$-vertex graphs $G$ in $\G$.
A class of graphs is {\it minor-closed} if every minor of a member of this class belongs to this class; a minor-closed family is {\it proper} if it does not contain all graphs.
Every proper minor-closed family is hereditary and has bounded expansion \cite{m}.
So Theorem \ref{bdd_expan_intro} immediately gives a negative answer to Question \ref{que_minor_poly}.

A variant of Question \ref{que_minor_poly} is the following conjecture.

\begin{conj}[{\cite[Conjecture 2.6]{gpstz_turan}}] \label{conj_planar_turan}
For any finite set of graphs $\F$ and for any graph $H$, if $\G$ is the set of all planar graphs with no subgraph isomorphic to any member in $\F$, and\footnote{Note that the statement in \cite{gpstz_turan} does not include the condition $H \in \G$. But it is required. If $H \not \in \G$, then no graph in $\G$ can contain $H$ as a subgraph, so $\ex(H,\G,n)=0$.} $H \in \G$, then $\ex(H,\G,n)=\Theta(n^k)$ for some integer $k$. 
\end{conj}

Note that the class $\G$ in Conjecture \ref{conj_planar_turan} is not minor-closed, so it is incomparable with Question \ref{que_minor_poly}.
But this class $\G$ is still hereditary and has bounded expansion, so Theorem \ref{bdd_expan_intro} immediately gives a positive answer to Conjecture \ref{conj_planar_turan}. 
In fact, we do not require the finiteness of $\F$, and $\G$ can be replaced by an arbitrary hereditary class of graphs with bounded expansion.
Moreover, instead of forbidding members of $\F$ as subgraphs, we can obtain stronger results by forbidding members of $\F$ as induced subgraphs.
Note that for every class $\F$, there exists a class $\F'$ such that forbidding all members of $\F$ as subgraphs is equivalent to forbidding all members of $\F'$ as induced subgraphs; but there exists a class $\F''$ such that forbidding all graphs in $\F''$ as induced subgraphs cannot be described by forbidding graphs in any class as subgraphs.
The following corollary of Theorem \ref{bdd_expan_intro} solves both Question \ref{que_minor_poly} and Conjecture \ref{conj_planar_turan} simultaneously.

\begin{corollary}
For any proper minor-closed family $\C$ and any (not necessarily nonempty or finite) set $\F$ of graphs, if $\G$ is the set of all graphs in $\C$ with no induced subgraph isomorphic to any member in $\F$, and $H \in \G$, then $\ex(H,\G,n)=\Theta(n^k)$ for some integer $k$. 
\end{corollary}

The most general classes in sparsity theory are the nowhere dense classes.
A class $\F$ of graphs is {\it somewhere dense} if there exists an integer $t$ such that every complete graph is a $t$-shallow minor of a graph in $\F$.
A class of graphs is {\it nowhere dense} if it is not somewhere dense.
Even though the definitions for these two kinds of classes look artificial, there are a number of equivalent definitions showing that these two classes are the right notions for capturing the dichotomy about sparse and dense graphs \cite{no_nowhere_dense,no_book}.

Theorem \ref{bdd_expan_intro} cannot be generalized to nowhere dense classes even when $H=K_2$: the class of all graphs whose maximum degree are at most their girth is nowhere dense and hereditary, and $n$-vertex graphs in this class have $O(n^{1+o(1)})$ edges, but there exists no constant $c$ such that every graph in this class has at most $cn$ edges \cite{no_book}. 
However, our main result (Theorem \ref{main}) applies to nowhere dense classes as well, by giving a slightly weaker estimate comparing to Theorem \ref{bdd_expan_intro}.

For any quasi-order $Q$, legal $Q$-labelled graph $(H,f_H)$, class $\G$ of $Q$-labelled graphs, and consistent set $\HH$ of homomorphisms from $(H,f_H)$ to $\G$, the {\it asymptotic logarithmic density for $(H,f_H,\G,\HH)$} is defined to be 
$$\limsup_{n \to \infty}\frac{\log(\ex(H,f_H,\HH,\G,n))}{\log n}.$$
Note that if the asymptotic logarithmic density for $(H,f_H,\G,\HH)$ is $k$, then $\ex(H,f_H,\HH,\G,n) \allowbreak \leq c n^{k+f(n)}$ for some constant $c$ and function $f(n)=o(1)$.
So this notion is slightly weaker than the estimate in Theorem \ref{bdd_expan_intro}.

Our main result (Theorem \ref{main}) implies the following analog of Theorem \ref{bdd_expan_intro} in terms of nowhere dense classes and asymptotic logarithmic density.

\begin{theorem} \label{nowhere_dense_intro}
Let $Q$ be a quasi-order with finite ground set.
Let $(H,f_H)$ be a legal $Q$-labelled graph.
Let $\G$ be a hereditary nowhere dense class of legal $Q$-labelled graphs.
Let $\HH$ be a consistent set of homomorphisms from $(H,f_H)$ to $\G$.
Then the asymptotic logarithmic density for $(H,f_H,\G,\HH)$ equals $\dup_\HH(H,f_H,\G)$.
\end{theorem}

Theorem \ref{nowhere_dense_intro} strengthens a result of Ne\v{s}et\v{r}il and Ossona de Mendez \cite{no_count} who proved that the asymptotic logarithmic density for the number of $H$-induced subgraphs in graphs in a hereditary nowhere dense class is an integer between $0$ and $\alpha(H)$.
Note that the set of induced subgraph embeddings is a consistent set of homomorphisms, and any independent collection of separations of $H$ has size at most $\alpha(H)$.

The machinery developed for the proof of Theorem \ref{main} works for some somewhere dense classes as well.

For any nonnegative integer $d$, a graph is {\it $d$-degenerate} if every its subgraph has minimum degree at most $d$.
The class of $d$-degenerate graphs is somewhere dense when $d \geq 2$, since this class contains all graphs that can be obtained by subdividing every edge of a complete graph. 
The following is an equivalent statement of a conjecture of Huynh and Wood \cite{hw}.
(For a nonnegative integer $d$ and graph $H$, $\flap_d(H)$ is defined to be the maximum size of an independent collection of separations of $H$ of order at most $d$.)

\begin{conj}[{\cite[Conjecture 16]{hw}}] \label{deg_conj}
Let $d$ be a nonnegative integer.
If $\D_d$ is the class of $d$-degenerate graphs and $H$ is a $d$-degenerate graph, then $\ex(H,\D_d,n)=\Theta(n^{\flap_d(H)})$.
\end{conj}

We disprove Conjecture \ref{deg_conj} and prove the correct exponent for $n$.
For a nonnegative integer $d$ and a graph $H$, we define $\alpha_d(H)$ to be the maximum size of a set of pairwise non-adjacent vertices of degree at most $d$ in $H$.

\begin{theorem} \label{bdd_degen_intro}
Let $d$ be a nonnegative integer.
If $\D_d$ is the class of $d$-degenerate graphs and $H$ is a $d$-degenerate graph, then $\ex(H,\D_d,n)=\Theta(n^{\alpha_d(H)})$.
Moreover, there are infinitely many $d$-degenerate graphs $H$ with $\alpha_d(H) < \flap_d(H)$ when $d \geq 2$.
\end{theorem}

\subsection{Organization of this paper}

In Section \ref{sec:counting}, we will prove a key lemma about counting the number of homomorphisms: if $G$ is a graph with a ``nice'' ordering, then one can bound the size of any given set of homomorphisms from a graph $H$ to $G$ and construct an independent collection of separations of $H$ matching the size for this bound.
Graphs with bounded degeneracy have such a nice ordering.
In Section \ref{sec:bdd_degen}, we prove Theorem \ref{bdd_degen_intro} by showing that the independent collection obtained in the previous lemma is duplicable.

The next goal is to prove our main theorem (Theorem \ref{main}).
The strategy is to prove that graphs mentioned in Theorem \ref{main} has a nice ordering as stated in the lemma in Section \ref{sec:counting} and prove that the corresponding independent collection is duplicable.
In Section \ref{sec:struct_sm}, we prove structure theorems for graphs whose shallow minors have low edge-density.
These structure theorems will be used in Section \ref{sec:homo_density} to prove Theorem \ref{main}.
We will also show how to deduce Theorems \ref{bdd_expan_intro} and \ref{nowhere_dense_intro} from Theorem \ref{main} in Section \ref{sec:homo_density}.

In Section \ref{sec:concrete}, we will use Theorem \ref{main} (or more precisely, Corollary \ref{cor_bdd_expan} which is a more informative version of Theorem \ref{bdd_expan_intro}) to deduce results about $\ex(H,\G,n)$, including solutions of some open questions in the literature.
Results in Section \ref{sec:concrete} can be read without any knowledge in earlier sections except the statement of Corollary \ref{cor_bdd_expan}. 

\subsection{Notations}

Now we define some notations that will be used in this paper.
For any vertex $x$ of a graph $G$ and any (possibly negative) real number $\ell$, we define $N_G^{\leq \ell}[x]$ to be the set of all the vertices in $G$ whose distance to $x$ is at most $\ell$; in particular, $N_G^{\leq 0}[x]=\{x\}$ and $N_G^{\leq -1}[x]=\emptyset$.
We also denote $N_G^{\leq 1}[x]$ by $N_G[x]$.
For every subset $S$ of $V(G)$, we define $N_G(S)=\{v \in V(G)-S: uv \in E(G)$ for some $u \in S\}$; note that $N_G[S]=N_G(S) \cup S$.
When $S$ consists of only one vertex $x$, we write $N_G(S)$ as $N_G(x)$.
For a subset $S$ of $V(G)$, we define $G[S]$ to be the subgraph of $G$ induced by $S$.
For a function $f$ whose domain is a set of marches over a set $S$, if $T$ is a subset of $S$, then we define $f|_T$ to be the function obtained from $f$ by restricting the domain to be $\{m: m$ is a march in the domain of $f$, and all entries of $m$ are in $T\}$. 
For a positive integer $k$, $[k]$ denotes the set $\{1,2,...,k\}$.

\section{A counting lemma} \label{sec:counting}

Let $G$ be a graph.
An {\it ordering} of $G$ is a bijection from $V(G)$ to $[\lvert V(G) \rvert]$.
Let $\sigma$ be an ordering of $G$.
Given $i \in [\lvert V(G) \rvert]$, we define $G_{\sigma, \geq i}$ to be the subgraph of $G$ induced by $\{v \in V(G): \sigma(v) \geq i\}$.
Let $p$ be a nonnegative integer.
Let $\Se=(S_i: i \in [\lvert V(G) \rvert])$ be a sequence such that for each $i \in [\lvert V(G) \rvert]$, $S_i \subseteq \{v \in V(G): \sigma(v) \geq i+1\}$.
For each $i \in [\lvert V(G) \rvert]$, the {\it $p$-basin (with respect to $\sigma$ and $\Se$) at $i$} is the set $N_{G_{\sigma, \geq i}-S_i}^{\leq p}[\sigma^{-1}(i)]$.
That is, the $p$-basin at $i$ consists of the vertices in $G_{\sigma,\geq i}$ that can be reached from $\sigma^{-1}(i)$ by a path in $G_{\sigma,\geq i}$ of length at most $p$ disjoint from $S_i$.

\begin{lemma} \label{container}
For any integers $k \geq 0, h \geq 1, N \geq 0$, there exists a positive integer $c=c(k,h,N)$ such that the following holds.
Let $H$ be a graph on $h$ vertices.
Let $G$ be a graph, and let $\sigma$ be an ordering of $G$.
Let $\Se=(S_i: i \in [\lvert V(G) \rvert])$ be a sequence such that $S_i \subseteq \{v \in V(G): \sigma(v) \geq i+1\}$ and $\lvert S_i \rvert \leq k$ for every $i \in [\lvert V(G) \rvert]$.
For every $i \in [\lvert V(G) \rvert]$, let $B_i$ be the $(h-1)$-basin with respect to $\sigma$ and $\Se$ at $i$.
Let $b$ be a positive real number.
If $\lvert B_i \rvert \leq b$ for every $i \in [\lvert V(G) \rvert]$, then for every set $\HH$ of homomorphisms from $H$ to $G$, there exists an integer $t$ such that the following hold.
	\begin{enumerate}
		\item $\lvert \HH \rvert \leq cb^h\lvert V(G) \rvert^t$, 
		\item For every $\phi \in \HH$, there exist an independent collection $\L_\phi$ of separations of $H$ with $\lvert \L_\phi \rvert \leq t$ and an injection $\iota_\phi: \L_\phi \rightarrow [\lvert V(G) \rvert-N]$ such that for every $(X,Y) \in \L_\phi$, if $i=\iota_\phi((X,Y))$, then 
			\begin{enumerate}
				\item $\sigma^{-1}(i) \in \phi(X-Y) \subseteq B_i$, 
				\item $\phi(X \cap Y) \subseteq S_i$, 
				\item for every component $C$ of $H[X-Y]$, $\sigma^{-1}(i) \in \phi(V(C))$, and 
				\item every vertex in $X \cap Y$ is adjacent in $H$ to some vertex in $X-Y$.
			\end{enumerate}
		\item There exists $\phi \in \HH$ with $\lvert \L_\phi \rvert=t$. 
	\end{enumerate}
\end{lemma}

\begin{pf}
Let $k \geq 0, h \geq 1, N \geq 0$ be integers.
Let $c_1=h(k+h)^{2h-1}+N$. 
Define $c=c_1^h$. 

Let $H,G,\sigma,\Se=(S_i: i \in [\lvert V(G) \rvert])$, $b$, $B_1,B_2,...,B_{\lvert V(G) \rvert}$ be as stated in the lemma.
For each $i \in [\lvert V(G) \rvert]$, let $v_i$ be the vertex of $G$ with $\sigma(v_i)=i$.

Let $\HH$ be a set of homomorphisms from $H$ to $G$.
For every $\phi \in \HH$, we define the following.
	\begin{itemize}
		\item Let $H_{\phi,0}=H$, $Z_{\phi,0}=\emptyset$ and $I_{\phi,0}=\emptyset$.
		\item For each $i \in [\lvert V(G) \rvert]$, 
			\begin{itemize}
				\item if $v_i \in \phi(V(H_{\phi,i-1}))$, then 
					\begin{itemize}
						\item let $Z_{\phi,i}=Z_{\phi,i-1} \cup \{v_i\}$,
						\item let $\Q_{\phi,i}$ be the collection of the components of $H_{\phi,i-1}-\{u \in V(H): \phi(u) \in S_i\}$ containing a vertex $u \in V(H)$ with $\phi(u)=v_i$, 
						\item let $H_{\phi,i}=H_{\phi,i-1}-\bigcup_{Q \in \Q_{\phi,i}}V(Q)$, 
						\item define $I_{\phi,i}$ as follows:
							\begin{itemize}
								\item if $i \leq \lvert V(G) \rvert-N$ and $\bigcup_{j \in [i-1], v_j \in Z_{\phi,i-1}}S_j \cap \bigcup_{Q \in \Q_{\phi,i}}\phi(V(Q)) = \emptyset$, then let $I_{\phi,i}=I_{\phi,i-1} \cup \{v_i\}$,
								\item otherwise, let $I_{\phi,i}=I_{\phi,i-1}$;
							\end{itemize}
					\end{itemize}
				\item otherwise, let $Z_{\phi,i}=Z_{\phi,i-1}$, $H_{\phi,i}=H_{\phi,i-1}$ and $I_{\phi,i}=I_{\phi,i-1}$. 
			\end{itemize}
		\item Let $Z_\phi=Z_{\phi,\lvert V(G) \rvert}$ and $I_\phi=I_{\phi,\lvert V(G) \rvert}$. 
	\end{itemize}

\noindent{\bf Claim 1:} For every $\phi \in \HH$, $\phi(V(H)) \subseteq \bigcup_{v \in Z_\phi}B_{\sigma(v)}$.

\noindent{\bf Proof of Claim 1:}
Let $\phi \in \HH$.
Let $w$ be a vertex of $H$.
Since $H_{\lvert V(G) \rvert}$ is empty, there exists $i_w \in [\lvert V(G) \rvert]$ such that $w \in V(H_{\phi,i_w-1})-V(H_{\phi,i_w})$.
So there exists $Q \in \Q_{\phi,i_w}$ such that $Q$ is a component of $H_{\phi,i_w-1}-\{u \in V(H): \phi(u) \in S_{i_w}\}$ containing $w$ and a vertex $u \in V(H)$ with $\phi(u)=v_{i_w} \in Z_{\phi,i_w}$.
Hence there exists a path $P$ in $G_{\sigma, \geq i_w}-S_{i_w}$ from $v_{i_w}$ to $\phi(w)$ with length at most $\lvert V(H) \rvert-1=h-1$.
So $\phi(w) \in B_{i_w}=B_{\sigma(v_{i_w})} \subseteq \bigcup_{v \in Z_\phi}B_{\sigma(v)}$.
This proves the claim.
$\Box$

\medskip

For every $\phi \in \HH$, 
	\begin{itemize}
		\item let $D_\phi$ be the directed graph with $V(D_\phi)=V(G)$ and $E(D_\phi)=\{(v_i,u): i \in [\lvert V(G) \rvert], u \in S_i\} \cup \{(u,v_i): i \in [\lvert V(G) \rvert], v_i \in Z_\phi, u \in \bigcup_{Q \in \Q_{\phi,i}}\phi(V(Q))\}$,
		\item let $R_\phi=\{v \in V(D_\phi):$ there exists a directed path in $D_\phi$ from a vertex in $I_\phi$ to $v$ of length at most $2h-2\} \cup \{v \in V(D_\phi): \sigma(v) \geq \lvert V(G) \rvert-N+1\}$, and 
		\item let $W_\phi=\bigcup_{v \in R_\phi}B_{\sigma(v)}$. 
	\end{itemize}

\noindent{\bf Claim 2:} For every $\phi \in \HH$ and for every $v \in Z_\phi-I_\phi$ with $\sigma(v) \leq \lvert V(G) \rvert-N$, there exist $u \in Z_\phi$ with $\sigma(u)<\sigma(v)$ and a directed path in $D_\phi$ from $u$ to $v$ with length at most two.

\noindent{\bf Proof of Claim 2:}
Let $\phi \in \HH$.
Let $v \in Z_\phi-I_\phi$ with $\sigma(v) \leq \lvert V(G) \rvert-N$.
Since $v \not \in I_\phi$, by the definition of $I_\phi$, $\bigcup_{j \in [\sigma(v)-1], v_j \in Z_{\phi,\sigma(v)-1}}S_j \cap \bigcup_{Q \in \Q_{\phi,\sigma(v)}}\phi(V(Q)) \neq \emptyset$.
So there exist $u \in Z_{\phi,\sigma(v)-1} \subseteq Z_\phi$ with $\sigma(u)<\sigma(v)$ and $w \in S_{\sigma(u)} \cap \bigcup_{Q \in \Q_{\phi,\sigma(v)}}\phi(V(Q))$.
Hence $uwv$ is a directed path in $D_\phi$ from $u$ to $v$ of length at most two.
$\Box$

\medskip

\noindent{\bf Claim 3:} For every $\phi \in \HH$, $\phi(V(H)) \subseteq W_\phi$.

\noindent{\bf Proof of Claim 3:}
Let $\phi \in \HH$.
By Claim 2, for every $v \in Z_\phi$ with $\sigma(v) \leq \lvert V(G) \rvert-N$, there exists a directed path in $D_\phi$ from a vertex in $I_\phi$ to $v$ with length at most $2(\lvert Z_\phi \rvert-1) \leq 2(\lvert V(H) \rvert-1)=2h-2$.
Hence $Z_\phi \subseteq R_\phi$.
By Claim 1, $\phi(V(H)) \subseteq \bigcup_{v \in Z_\phi}B_{\sigma(v)} \subseteq \bigcup_{v \in R_\phi}B_{\sigma(v)}=W_\phi$. 
$\Box$

\medskip

\noindent{\bf Claim 4:} For every $\phi \in \HH$, every vertex of $D_\phi$ has out-degree at most $k+h$.

\noindent{\bf Proof of Claim 4:}
Let $\phi \in \HH$.
Let $x \in V(D_\phi)$.
Since $\lvert Z_\phi \rvert \leq \lvert V(H) \rvert$, there exists at most $h$ integers $i \in [\lvert V(G) \rvert]$ with $v_i \in Z_\phi$ such that $x \in \bigcup_{Q \in \Q_{\phi,i}}\phi(V(Q))$.
So the out-degree of $x$ in $D_\phi$ is at most $\lvert S_{\sigma(x)} \rvert+h \leq k+h$.
$\Box$

\medskip

Define $t=\max_{\phi \in \HH} \lvert I_\phi \rvert$.
Since $\lvert I_\phi \rvert \leq \lvert Z_\phi \rvert \leq \lvert V(H) \rvert$ for every $\phi \in \HH$, $t \leq h$.

\medskip

\noindent{\bf Claim 5:} There exist ${\lvert V(G) \rvert \choose t}$ subsets of $V(G)$ with size at most $c_1b$ such that for any $\phi \in \HH$, $\phi(V(H))$ is contained in at least one of those sets.

\noindent{\bf Proof of Claim 5:}
By the definition of $t$, $\lvert I_\phi \rvert \leq t$ for every $\phi \in \HH$.
So there exist ${\lvert V(G) \rvert \choose t}$ subsets of $V(G)$ with size $t$ such that for any $\phi \in \HH$, $I_\phi$ is contained in at least one of them.
By Claim 4, for every $\phi \in \HH$, every vertex of $D_\phi$ has out-degree at most $k+h$, so for every $v \in V(D_\phi)$, there exist at most $\sum_{j=0}^{2h-2}(k+h)^j \leq (k+h)^{2h-1}$ directed paths in $D_\phi$ of length at most $2h-2$ starting from $v$.
So there exist ${\lvert V(G) \rvert \choose t}$ subsets of $V(G)$ with size at most $t(k+h)^{2h-1} \leq h(k+h)^{2h-1}$ such that for any $\phi \in \HH$, $R_\phi-\{v \in V(D_\phi): \sigma(v) \geq \lvert V(G) \rvert-N+1\}$ is contained in at least one of those sets.
Hence there exist ${\lvert V(G) \rvert \choose t}$ subsets of $V(G)$ with size at most $h(k+h)^{2h-1}+N =c_1$ such that for any $\phi \in \HH$, $R_\phi$ is contained in at least one of those sets.
Since $\lvert B_i \rvert \leq b$ for every $i \in [\lvert V(G) \rvert]$, there exist ${\lvert V(G) \rvert \choose t}$ subsets of $V(G)$ with size at most $c_1 \cdot b$ such that for any $\phi \in \HH$, $W_\phi$ is contained in at least one of those sets.
Then this claim follows from Claim 3.
$\Box$

\medskip

Note that for each subset $U$ of $V(G)$, there are at most $\lvert U \rvert^{\lvert V(H) \rvert}$ homomorphisms from $H$ to $G$ whose image is contained in $U$.
So by Claim 5, $\lvert \HH \rvert \leq {\lvert V(G) \rvert \choose t} \cdot (c_1b)^h \leq cb^h\lvert V(G) \rvert^t$.
This proves Statement 1 of this lemma.

Now we prove Statements 2 and 3.
Note that for every $\phi \in \HH$ and $v \in Z_\phi$, $v$ is the unique vertex in $Z_{\phi,\sigma(v)}-Z_{\phi,\sigma(v)-1}$, and $\Q_{\phi,\sigma(v)}$ is defined; we denote $\Q_{\phi,\sigma(v)}$ by $\Q_{\phi,v}$ for simplicity. 

\medskip

\noindent{\bf Claim 6:} For every $\phi \in \HH$ and for any $v \in I_\phi$ and $Q \in \Q_{\phi,v}$, $N_H(V(Q)) \subseteq \{u \in V(H): \phi(u) \in S_{\sigma(v)}\}$.  

\noindent{\bf Proof of Claim 6:}
Let $\phi \in \HH$.
Let $v \in I_\phi$.
Let $Q \in \Q_{\phi,v}$.
Suppose to the contrary that there exists a vertex $w \in N_H(V(Q)) - \{u \in V(H): \phi(u) \in S_{\sigma(v)}\}$. 
Since $V(Q)$ is the vertex-set of a component of $H_{\phi,\sigma(v)-1}-\{u \in V(H): \phi(u) \in S_{\sigma(v)}\}$, $w \in V(H)-V(H_{\phi,\sigma(v)-1})$.
So there exists $i_w \in [\sigma(v)-1]$ such that $w \in V(H_{\phi,i_w-1})-V(H_{\phi,i_w})$.
Hence there exists a component $Q' \in \Q_{\phi,v_{i_w}}$ of $H_{\phi,i_w-1}-\{u \in V(H): \phi(u) \in S_{i_w}\}$ containing $w$ and a vertex $w' \in V(H)$ with $\phi(w')=v_{i_w} \in Z_{\phi,i_w} \subseteq Z_{\phi,\sigma(v)-1}$.
Since $w \in N_H(V(Q))$, $w \in V(Q')$ is adjacent in $H$ to some vertex $w''$ in $Q \subseteq H_{\phi,\sigma(v)-1} = H_{\phi,\sigma(v)-1}-V(Q') \subseteq H_{\phi,i_w-1}-V(Q')$, so $w'' \in \{u \in V(H): \phi(u) \in S_{i_w}\}$.
Hence $\phi(w'') \in S_{i_w} \cap \phi(V(Q)) \subseteq \bigcup_{j \in [\sigma(v)-1], v_j \in Z_{\phi,\sigma(v)-1}}S_j \cap \bigcup_{Q'' \in \Q_{\phi,v}}V(Q'')$. 
Therefore, $v \not \in I_\phi$, a contradiction.
$\Box$

\medskip

For every $\phi \in \HH$ and every $v \in I_\phi$, let $X_{\phi,v} = N_H[\bigcup_{Q \in \Q_{\phi,v}}V(Q)]$, and let $Y_{\phi,v}=V(H)-\bigcup_{Q \in \Q_{\phi,v}}V(Q)$, so $(X_{\phi,v},Y_{\phi,v})$ is a separation of $H$ with $X_{\phi,v} \cap Y_{\phi,v} \subseteq \{u \in V(H): \phi(u) \in S_{\sigma(v)}\}$ by Claim 6. 
Note that $X_{\phi,v} \cap Y_{\phi,v} = N_H(\bigcup_{Q \in \Q_{\phi,v}}V(Q))$ and $X_{\phi,v}-Y_{\phi,v} = \bigcup_{Q \in \Q_{\phi,v}}V(Q)$, so every vertex in $X_{\phi,v} \cap Y_{\phi,v}$ is adjacent in $H$ to some vertex in $X_{\phi,v}-Y_{\phi,v}$.

For every $\phi \in \HH$, let $\L_\phi = \{(X_{\phi,v},Y_{\phi,v}): v \in I_\phi\}$, and for every $(X_{\phi,v},Y_{\phi,v}) \in \L_\phi$, let $\iota_\phi((X_{\phi,v},Y_{\phi,v})) = \sigma(v)$.
Clearly, $\iota_\phi$ is an injection from $\L_\phi$ to $[\lvert V(G) \rvert]$.
Note that $v \in \phi(V(Q))$ for every $Q \in \Q_{\phi,v}$.
Since $X_{\phi,v}-Y_{\phi,v} = \bigcup_{Q \in \Q_{\phi,v}}V(Q) = V(H_{\phi,\sigma(v)-1})-V(H_{\phi,\sigma(v)})$, $\L_\phi$ is an independent collection of separations of $H$.
Moreover, every component of $H[X_{\phi,v}-Y_{\phi,v}]$ is a member of $\Q_{\phi,v}$, so $v$ is contained in $\phi(V(C))$ for every component $C$ of $H[X_{\phi,v}-Y_{\phi,v}]$.

For every $\phi \in \HH$ and $v \in I_\phi$, since $\lvert V(H) \rvert=h$ and $v \in \phi(\bigcup_{Q \in \Q_{\phi,v}}V(Q)) \subseteq V(G_{\sigma, \geq \sigma(v)})-S_{\sigma(v)}$, we know $\phi(\bigcup_{Q \in \Q_{\phi,v}}V(Q)) \subseteq B_{\sigma(v)}$.
So for any $\phi \in \HH$ and $(X,Y) \in \L_\phi$, there exists $v \in I_\phi$ with $\sigma(v)=\iota_\phi((X,Y))$ such that $v \in \phi(X-Y) \subseteq B_{\sigma(v)}$, $\phi(X \cap Y) \subseteq S_{\sigma(v)}$, and for every component $C$ of $H[X-Y]$, $v \in \phi(V(C))$. 
Note that $v \in I_\phi$ implies that $\sigma(v) \leq \lvert V(G) \rvert-N$.
So the image of $\iota_\phi$ is contained in $[\lvert V(G) \rvert-N]$.

Since $\lvert \L_\phi \rvert = \lvert I_\phi \rvert$, Statements 2 and 3 of this lemma follow from the definition of $t$.
\end{pf}

\section{Bounded degeneracy} \label{sec:bdd_degen}

We will prove Theorem \ref{bdd_degen_intro} in this section.
A {\it subgraph embedding} from a graph $H$ to a graph $G$ is an injective homomorphism from $H$ to $G$.
Note that the ratio between the number of subgraph embeddings from $H$ to $G$ and the number of $H$-subgraphs of $G$ are upper bounded and lower bounded by constants only depending on $H$.

Recall that for any nonnegative integer $d$ and a graph $G$, $\alpha_d(G)$ is the maximum size of a set of pairwise non-adjacent vertices of degree at most $d$ in $G$.

\begin{lemma} \label{deg_upper}
For any integers $d \geq 0$ and $h \geq 1$, there exists an integer $c=c(d,h)$ such that the following holds.
If $H$ is a graph on $h$ vertices and $G$ is a $d$-degenerate graph, then there are at most $c\lvert V(G) \rvert^{\alpha_d(H)}$ subgraph embeddings from $H$ to $G$.
\end{lemma}

\begin{pf}
Let $d \geq 0$ and $h \geq 1$ be integers.
Define $c=c_{\ref{container}}(d,h,0)$, where $c_{\ref{container}}$ is the number $c$ mentioned in Lemma \ref{container}.

Let $H$ be a graph on $h$ vertices.
Let $G$ be a $d$-degenerate graph.
Since $G$ is $d$-degenerate, there exists an ordering $\sigma$ of $G$ such that for every $v \in V(G)$, $v$ has at most $d$ neighbors in $G_{\sigma,\geq \sigma(v)}$. 
For each $i \in [\lvert V(G) \rvert]$, let $S_i = N_{G_{\sigma,\geq i}}(\sigma^{-1}(i))$, so $S_i \subseteq \{v \in V(G): \sigma(v) \geq i+1\}$ and $\lvert S_i \rvert \leq d$.
Let $\Se=(S_i: i \in [\lvert V(G) \rvert])$.
For each $i \in [\lvert V(G) \rvert]$, let $v_i$ be the vertex of $G$ with $\sigma(v_i)=i$, and let $B_i$ be the $(h-1)$-basin with respect to $\sigma$ and $\Se$ at $i$.
Note that each $B_i$ equals $\{v_i\}$, so $\lvert B_i \rvert=1$. 
By Lemma \ref{container}, there exists an integer $t$ such that there are at most $c\lvert V(G) \rvert^t$ subgraph embeddings from $H$ to $G$, and there exist a subgraph embedding $\phi$ from $H$ to $G$ and an independent collection $\L_\phi$ of separations of $H$ with $\lvert \L_\phi \rvert=t$ such that for every $(X,Y) \in \L_\phi$, there exists $i_X \in [\lvert V(G) \rvert]$ such that $v_{i_X} \in \phi(X-Y) \subseteq B_{i_X}=\{v_{i_X}\}$ and $\phi(X \cap Y) \subseteq S_{i_X}$.
For every $(X,Y) \in \L_\phi$, $\phi(X-Y)=\{v_{i_X}\}$, so there exists $u_X \in X-Y$ with $\phi(u_X)=v_{i_X}$.
Since $\phi$ is a subgraph embedding, $\phi$ is injective.
So for every $(X,Y) \in \L_\phi$, $X-Y=\{u_X\}$; since $\phi(X \cap Y) \subseteq S_{i_X}$, the degree of $u_X$ in $H$ is at most $\lvert S_{i_X} \rvert \leq d$.
Therefore, $\{u_X: (X,Y) \in \L_\phi\}$ is a set of pairwise non-adjacent vertices of degree at most $d$ in $H$.
Hence $t=\lvert \L_\phi \rvert \leq \alpha_d(H)$
\end{pf}

\begin{theorem} \label{bdd_degen}
Let $d \geq 0$ be an integer.
Let $\D_d$ be the class of all $d$-degenerate graphs.
Let $H$ be a graph.
If $H$ is $d$-degenerate, then there exist real numbers $c_1,c_2$ such that $c_1n^{\alpha_d(H)} \leq \ex(H,\D_d,n) \leq c_2n^{\alpha_d(H)}$; if $H$ is not $d$-degenerate, then $\ex(H,\D_d,n)=0$.
\end{theorem}

\begin{pf}
If $H$ is not $d$-degenerate, then $H$ cannot be a subgraph of any graph in $\D_d$, so $\ex(H,\D_d,n)=0$.
Hence we may assume that $H$ is $d$-degenerate.

By Lemma \ref{deg_upper}, there exists a real number $c_2$ such that $\ex(H,\D_d,n) \leq c_2 n^{\alpha_d(H)}$.
Let $S$ be a set of pairwise non-adjacent vertices of degree at most $d$ in $H$. 
Then $\L=\{(N_H[v], V(G)-\{v\}): v \in S\}$ is an independent collection of separations of $H$ of order at most $d$.
For any positive integer $k$, let $H_k = H \wedge_k \L$.

We shall show that for every positive integer $k$, $H_k$ is $d$-degenerate.
Let $k$ be a positive integer and let $R$ be a subgraph of $H_k$.
If $R$ contains a copy of some vertex in $S$, then $R$ contains a vertex of degree at most $d$.
If $R$ does not contain any copy of a vertex in $S$, then $R$ is a subgraph of $H-S$ and hence contains a vertex of degree at most $d$.
So $H_k$ is $d$-degenerate.

Hence, $H \wedge_k \L \in \G$ for all positive integers $k$.
Note that every graph is a legal $Q$-labelled graph for some quasi-order of size 1.
So $\L$ is $\G$-duplicable by definition.
Let $\HH$ be the set of all subgraph embeddings from $H$ to graphs in $\G$.
Then $\L$ is $(\G,\HH)$-duplicable by definition.
By Proposition \ref{lower_bound}, there exists a real number $c_1$ such that $\ex(H,\D_d,n) \geq c_1 n^{\alpha_d(H)}$ since $\lvert \L \rvert = \alpha_d(H)$.
Therefore, $c_1 n^{\alpha_d(H)} \leq \ex(H,\D_d,n) \leq c_2 n^{\alpha_d(H)}$.
\end{pf}

\bigskip

Theorem \ref{bdd_degen_intro} is an immediate corollary of Theorem \ref{bdd_degen} and the following proposition.

\begin{prop}
Let $d$ be a positive integer with $d \geq 2$.
Then there are infinitely many $d$-degenerate graphs $H$ with $\alpha_d(H) < \flap_d(H)$.
\end{prop}

\begin{pf}
Define $H$ to be a graph whose vertex-set is a union of three disjoint set $X,Y,Z$ such that $H[Z]=K_{d+1}$, $Y$ is a stable set in $H$ with size $d$, every vertex in $Y$ is adjacent to exactly $d$ vertices in $Z$ such that every vertex in $Z$ is adjacent to at least one vertices in $Y$ (it is possible since $d \geq 2$), $X$ is a nonempty stable set in $H$, and the neighborhood of any vertex in $X$ is $Y$.

We first show that $H$ is $d$-degenerate.
Let $R$ be a subgraph of $H$.
If $R$ contains at least one vertex in $X$, then $R$ contains a vertex of degree at most $\lvert Y \rvert=d$.
If $V(R) \cap X=\emptyset$ and $V(R) \cap Y \neq \emptyset$, then the degree in $R$ of any vertex in $V(R) \cap Y$ is at most $d$.
If $V(R) \cap X = V(R) \cap Y = \emptyset$, then $R$ is a subgraph of $H[Z]=K_{d+1}$, so $R$ contains a vertex of degree at most $d$.
Hence $H$ is $d$-degenerate.

Since $X$ is a stable set in $H$, and every vertex in $X$ has degree at most $d$, we know $\alpha_d(H) \geq \lvert X \rvert$.
Since every vertex of $H$ with degree at most $d$ is contained in $X$, $\alpha_d(H)=\lvert X \rvert$.
But $\{(Y \cup Z, X \cup Y), (\{x\} \cup Y, Y \cup Z): x \in X\}$ is an independent collection of separations of $H$ of order at most $\lvert Y \rvert=d$.
So $\flap_d(H) \geq \lvert X \rvert+1>\alpha_d(H)$.
Note that $\lvert X \rvert$ can arbitrary, so there are infinitely many such graphs $H$.
\end{pf}

\section{Structure for graphs with forbidden shallow minors} \label{sec:struct_sm}

In this section we prove structural theorems for graphs whose shallow minors have low edge-density.

Let $s,t$ be positive integers.
Let $G$ be a graph.
For $p,q \in {\mathbb N} \cup \{\infty\}$, a {\it $(p,q)$-model for a $K_{s,t}$-minor in $G$} is an ordered pair of two collections $(\X,\Y)$ such that
	\begin{itemize}
		\item $\lvert \X \rvert=s$, $\lvert \Y \rvert=t$, 
		\item $\X \cup \Y$ is a collection of pairwise disjoint connected subgraphs of $G$, 
		\item for each $X \in \X$ and $Y \in \Y$, there exists an edge of $G$ between $V(X)$ and $V(Y)$,
		\item every member of $\X$ contains at most $p$ vertices, and
		\item every member of $\Y$ contains at most $q$ vertices.
	\end{itemize}
Note that if there exists a $(p,q)$-model for a $K_{s,t}$-minor in $G$, then $G$ contains $K_{s,t}$ as a $\max\{p,q\}$-shallow minor.
For collections $\C_1,\C_2$ of subgraphs of a graph $G$, {\it a $(p,q,\C_1,\C_2)$-model for a $K_{s,t}$-minor in $G$} is a $(p,q)$-model $(\X,\Y)$ for a $K_{s,t}$-minor such that $\X \subseteq \C_1$ and $\Y \subseteq \C_2$. 

For a set $S$ of nonnegative integers, we say a graph $G$ is an {\it $S$-subdivision} of a graph $H$ if for every $e \in E(H)$, there exists $s_e \in S$ such that $G$ can be obtained from $H$ by subdividing each edge $e$ of $H$ exactly $s_e$ times.
The vertices in $H$ are called the {\it branch vertices} of the $S$-subdivision.

For a graph $G$, a subset $Y$ of $V(G)$, and an integer $r$, we say a subgraph $H$ of $G$ is {\it $r$-adherent to $Y$} if $V(H) \cap Y = \emptyset$ and $\lvert N_G(V(H)) \cap Y \rvert \geq r$.

The proof of the following lemma is essentially the same as the proof of \cite[Lemma 4.2]{lw} but with more detailed analysis.

\begin{lemma} \label{r-adherent size}
For any positive integers $r,t$ and positive real number $k'$, there exists a real number $\alpha=\alpha(r,t,k')>0$ such that the following holds.
If $\ell \in {\mathbb N} \cup \{0,\infty\}$, $b \in {\mathbb N}$, $\epsilon' \geq 0$ is a real number, $G$ is a graph, $Y \subseteq V(G)$, and $\C$ is a collection of disjoint connected subgraphs of $G-Y$ on at most $b$ vertices such that each member of $\C$ is $r$-adherent to $Y$ and has radius at most $\ell$, then either
	\begin{enumerate}
		\item there exists a graph $H$ with $\lvert E(H) \rvert > k'\lvert V(H) \rvert^{1+\epsilon'}$ such that $G$ contains a subgraph isomorphic to a $[2\ell+1]$-subdivision of $H$ with all branch vertices in $Y$, 
		\item there exist a $(1,b, \{G[\{y\}]: y \in Y\}, \C)$-model for an $\ell$-shallow $K_{r,t}$-minor in $G$, or 
		\item $\lvert \C \rvert \leq \alpha \lvert Y \rvert^{1+\epsilon'(r-1)}$. 
	\end{enumerate}
\end{lemma}

\begin{pf}
Let $r \geq 1,t \geq 1$ be integers, and let $k'>0$ be a real number.
Define $\alpha=(t-1)(2k')^{r-1} + k'+t-1$.

Let $\ell,b,\epsilon',G,Y,\C$ be as stated in the lemma.
Let $G'$ be the graph obtained from $G$ by contracting each member of $\C$ into a vertex.
Note that $Y \subseteq V(G')$ since each member of $\C$ is disjoint from $Y$. 
Let $Z=V(G')-V(G)$.
Note that $Z$ is the set of the vertices of $G'$ obtained by contracting members of $\C$.
Define $G''$ to be the graph obtained from $G' - E(G'[Y])$ by repeatedly choosing a vertex $v$ in $Z$ adjacent in $G'$ to a pair of nonadjacent vertices $u,w$ in $(G' - E(G'[Y]))[Y]$, deleting $v$, and adding an edge $uw$, until for any remaining vertex in $Z$, its neighbors in $Y$ form a clique.

Let $H=G''[Y]$.
Note that some subgraph $H'$ of $G'$ is isomorphic to a $\{1\}$-subdivision of $H$ with all branch vertices in $Y$, and $V(H')-V(H) \subseteq Z$.
Since each vertex in $Z$ comes from contracting a member of $\C$, and each member of $\C$ has radius at most $\ell$, there exists a subgraph $H''$ of $G$ isomorphic to a $[2\ell+1]$-subdivision of $H$ with all branch vertices in $Y$.
So for every subgraph $H_0$ of $H$, some subgraph of $G$ is isomorphic to a $[2\ell+1]$-subdivision of $H_0$ with all branch vertices in $Y$.
Hence we may assume that for every subgraph $H_0$ of $H$, $\lvert E(H_0) \rvert \leq k'\lvert V(H_0) \rvert^{1+\epsilon'}$, for otherwise Statement 1 holds and we are done.
So for every subgraph $H_0$ of $H$, there exists a vertex of degree at most $2k'\lvert V(H_0) \rvert^{\epsilon'} \leq 2k'\lvert V(H) \rvert^{\epsilon'}$ in $H_0$.
By \cite[Lemma 18]{Wood-count-clique}, there are at most at most ${2k'\lvert V(H) \rvert^{\epsilon'} \choose r-1} \lvert V(H) \rvert \leq (2k')^{r-1}\lvert V(H) \rvert^{\epsilon'(r-1)+1} = (2k')^{r-1}\lvert Y \rvert^{\epsilon'(r-1)+1}$ cliques of size $r$ in $H$. 

We first assume that there exists a subgraph $L$ of $G''$ isomorphic to $K_{r,t}$ with a partition $\{R,T\}$ of $V(L)$, where $\lvert R \rvert=r$ and $\lvert T \rvert=t$ such that $R$ and $T$ are stable sets in $L$, $R \subseteq Y$ and $T \subseteq Z$. 
Since each member of $\C$ contains at most $b$ vertices and has radius at most $\ell$, there exists a $(1,b,\{G[\{y\}]: y\in Y\},\C)$-model for an $\ell$-shallow $K_{r,t}$-minor in $G$.
Therefore, Statement 2 holds.

Hence we may assume that there does not exist a subgraph $L$ of $G''$ isomorphic to $K_{r,t}$ with a partition $\{R,T\}$ of $V(L)$, where $\lvert R \rvert=r$ and $\lvert T \rvert=t$ such that $R$ and $T$ are stable sets in $L$, $R \subseteq Y$ and $T \subseteq Z$. 
This implies that for each clique $K$ in $G''[Y]$ of size $r$, $\lvert \{z \in Z \cap V(G''): K \subseteq N_{G''}(z)\} \rvert \leq t-1$.
In addition, for every $z \in Z \cap V(G'')$, since every member of $\C$ is $r$-adherent to $Y$, $N_{G''}(z) \cap Y$ is a clique consisting of at least $r$ vertices in $H$.
Hence $\lvert Z \cap V(G'') \rvert \leq (t-1)(2k')^{r-1}\lvert Y \rvert^{\epsilon'(r-1)+1}$. 

When $r=1$, if $\lvert Z \rvert \geq (t-1)\lvert Y \rvert+1$, then some vertex in $Y$ is adjacent in $G'$ to at least $t$ vertices in $Z$, so Statement 2 holds.
So if $r=1$, then we may assume that $\lvert Z \rvert \leq (t-1)\lvert Y \rvert \leq \alpha\lvert Y \rvert^{1+\epsilon'(r-1)}$, so Statement 3 holds since $\lvert \C \rvert = \lvert Z \rvert$.
Hence we may assume that $r \geq 2$.

By the definition of $G''$, the vertices in $Z$ but not in $V(G'')$ are the ones being deleted while adding edges between vertices in $Y$. 
Hence $\lvert Z-V(G'') \rvert \leq \lvert E(G''[Y]) \rvert = \lvert E(H) \rvert \leq k'\lvert V(H) \rvert^{1+\epsilon'} = k'\lvert Y \rvert^{1+\epsilon'} \leq k'\lvert Y \rvert^{1+\epsilon'(r-1)}$, since $r \geq 2$. 
So 
$$\lvert \C \rvert = \lvert Z \rvert = \lvert Z \cap V(G'') \rvert + \lvert Z-V(G'') \rvert \leq (t-1)(2k')^{r-1}\lvert Y \rvert^{\epsilon'(r-1)+1} + k'\lvert Y \rvert^{1+\epsilon'(r-1)} \leq \alpha \lvert Y \rvert^{1+\epsilon'(r-1)}.$$
Hence Statement 3 holds.
\end{pf}

\bigskip

Let $G$ be a graph.
For a subset $Y$ of $V(G)$, $v \in V(G)-Y$ and integers $k$ and $r$, we define a {\it $(v,Y,k,r)$-span} to be a connected subgraph $H$ of $G-Y$ containing $v$ such that $\lvert Y \cap N_G(V(H)) \rvert \geq r$, and for every vertex $u$ of $H$, there exists a path in $H$ from $v$ to $u$ of length at most $k$.
A $(v,Y,k,r)$-span $H$ is {\it minimal} if no proper subgraph of $H$ is a $(v,Y,k,r)$-span. 

The following lemma (Lemma \ref{island_1_prep}) is a strengthening of a result implicitly proved in earlier work of the author and Wei \cite[Lemma 4.4]{lw}.
A big portion of the proof of Lemma \ref{island_1_prep} is identical to the proof of \cite[Lemma 4.4]{lw}. 

\begin{lemma} \label{island_1_prep}
For any positive integers $r,t$, nonnegative integer $\ell$, real numbers $k'>0,\beta \geq 0$, there exists a real number $c=c(r,t,\ell,k',\beta)$ such that for any graph $G$, subset $Y_0$ of $V(G)$ and real number $1 \geq \epsilon' \geq 0$, either
	\begin{enumerate}
		\item there exists a graph $H$ with $\lvert E(H) \rvert > k'\lvert E(H) \rvert^{1+\epsilon'}$ such that some subgraph of $G$ is isomorphic to a $[2\ell+1]$-subdivision of $H$, 
		\item there exists a $(1,\ell r+1)$-model for an $\ell$-shallow $K_{r,t}$-minor in $G$, 
		\item $\lvert V(G) \rvert \leq c\lvert Y_0 \rvert^{\epsilon_0}$, where $\epsilon_0=(1+\epsilon'(r-1))^{(r+1)^\ell}$, or 
		\item there exist $X,Z \subseteq V(G)-Y_0$ with $Z \subseteq X$ and $\lvert Z \rvert > \beta\lvert V(G)-X \rvert^{1+\epsilon'(r-1)}$ such that
			\begin{enumerate}
				\item for any distinct $z,z' \in Z$, the distance in $G[X]$ between $z,z'$ is at least $\ell+1$, 
				\item for any $z \in Z$ and any $u \in X$ whose distance from $z$ in $G[X]$ is at most $\ell$, $\lvert N_G(u) -X \rvert  \leq r-1$, and 
				\item $\lvert N_G(N_{G[X]}^{\leq \ell-1}[z])-X \rvert \leq r-1$ for every $z \in Z$.
			\end{enumerate}
	\end{enumerate}
\end{lemma}

\begin{pf}
Let $r,t$ be positive integers, $\ell$ be a nonnegative integer, $k'$ be a positive real number, and $\beta$ be a nonnegative real number.
Let $\alpha=\alpha_{\ref{r-adherent size}}(r,t,k')$, where $\alpha_{\ref{r-adherent size}}$ is the number $\alpha$ mentioned in Lemma \ref{r-adherent size}. 
Let $\eta=\beta+1+(\ell r+1)\alpha$. 
Define $c=1+\eta^{(r+1)^{(r+1)^\ell}}$. 

Let $G$ be a graph.
Let $Y_0 \subseteq V(G)$.
Let $\epsilon'$ be a real number with $0 \leq \epsilon' \leq 1$.
Let $\epsilon_0=(1+\epsilon'(r-1))^{(r+1)^\ell}$.
We may assume that Statements 1 and 2 do not hold, for otherwise we are done.
We shall prove that Statements 3 or 4 hold.

Let $X_0=V(G)-Y_0$.
For each integer $i \geq 0$, we define the following.
	\begin{itemize}
		\item Define $\C_i$ to be a maximal collection of pairwise disjoint subgraphs of $G[X_i]$, where each member of $\C_i$ is a minimal $(v,Y_i,\ell,r)$-span for some vertex $v \in X_i$ satisfying that if $\ell \geq 1$, then $\lvert Y_i \cap N_G(N_{G-Y_i}^{\leq \ell-1}[v]) \rvert \geq r$. 
		\item $D_i = \bigcup_{L \in \C_i}V(L)$.
		\item $Z_i$ is a maximal subset of $X_i-D_i$ such that
			\begin{itemize}
				\item for any two distinct vertices in $Z_i$, the distance in $G[X_i]$ between them is at least $\ell+1$, and 
				\item for every $z \in Z_i$, $N_{G[X_i-D_i]}^{\leq \ell-1}[z] \cap N_G(D_i)=\emptyset$.
			\end{itemize}
		\item $X_{i+1} = X_i - (Z_i \cup D_i)$.
		\item $Y_{i+1} = Y_{i} \cup Z_i \cup D_i$.
	\end{itemize}

We remark that the definitions of $\C_i,D_i,Z_i,X_{i+1},Y_{i+1}$ for $i \geq 0$ are identical to their definitions in the proof of \cite[Lemma 4.4]{lw}, except the initial conditions for $X_0,Y_0$ are different. 
Hence the proof of following claim is identical to the proof of \cite[Claim 4.4.1 in Lemma 4.4]{lw}.

\medskip

\noindent{\bf Claim 1:} For any integer $i \geq 0$ and $z \in Z_i$, 
	\begin{itemize}
		\item $\lvert N_G(N_{G[X_i]}^{\leq \ell-1}[z])-X_i \rvert \leq r-1$, and 
		\item if $u \in N_{G[X_i]}^{\leq \ell}[z]$, then $\lvert N_G(u)-X_i \rvert \leq r-1$. 
	\end{itemize}

\medskip

If there exists a nonnegative integer $i^*$ such that $\lvert Z_{i^*} \rvert > \beta \lvert V(G)-X_{i^*} \rvert^{1+\epsilon'(r-1)}$, then by defining $X=X_{i^*}$ and $Z=Z_{i^*}$, we have $Z \subseteq X$ and $\lvert Z \rvert > \beta \lvert V(G)-X \rvert^{1+\epsilon'(r-1)}$ such that Statement 4 holds (Statements 4(a) follows from the definition of $Z_{i^*}$, and Statements 4(b) and 4(c) follow from Claim 1).

So we may assume that $\lvert Z_i \rvert \leq \beta \lvert V(G)-X_i \rvert^{1+\epsilon'(r-1)}=\beta \lvert Y_i \rvert^{1+\epsilon'(r-1)}$ for every nonnegative integer $i$. 

The proof of the following claim is identical to the proof of \cite[Claims 4.4.2-4.4.5 in Lemma 4.2]{lw}.

\medskip

\noindent{\bf Claim 2:} $X_0 \subseteq Y_{(r+1)^{\ell}}$.

\medskip

Now we are ready to complete the proof.
For each $i \geq 0$, every member of $\C_i$ is a minimal $(v,Y_i,\ell,r)$-span, so every member of $\C_i$ contains at most $\ell r+1$ vertices by \cite[Lemma 4.3]{lw}.
Hence for each $i \geq 0$, by Lemma \ref{r-adherent size}, $\lvert \C_i \rvert \leq \alpha\lvert Y_i \rvert^{1+\epsilon'(r-1)}$.
So for every nonnegative integer $i$,
$$
\lvert Y_{i+1}-Y_i \rvert = \lvert Z_i \rvert + \sum_{L \in \C_i} \lvert V(L) \rvert \leq \beta \lvert Y_i \rvert^{1+\epsilon'(r-1)} + \lvert \C_i \rvert \cdot (\ell r+1) \leq \beta \lvert Y_i \rvert^{1+\epsilon'(r-1)}+ (\ell r+1) \alpha \lvert Y_i \rvert^{1+\epsilon'(r-1)},
$$
so
$$\lvert Y_{i+1} \rvert \leq \lvert Y_i \rvert + (\beta+(\ell r+1)\alpha) \lvert Y_i \rvert^{1+\epsilon'(r-1)} \leq \eta\lvert Y_i \rvert^{1+\epsilon'(r-1)}.$$
Therefore, for every nonnegative integer $i$, $\lvert Y_i \rvert \leq \eta^{(2+\epsilon'(r-1))^i}\lvert Y_0 \rvert^{(1+\epsilon'(r-1))^i}$.

By Claim 2, 
$$\lvert X_0 \rvert \leq \lvert Y_{(r+1)^\ell} \rvert \leq \eta^{(2+\epsilon'(r-1))^{(r+1)^\ell}}\lvert Y_0 \rvert^{(1+\epsilon'(r-1))^{(r+1)^\ell}} \leq \eta^{(r+1)^{(r+1)^\ell}}\lvert Y_0 \rvert^{\epsilon_0},$$
since $\epsilon' \leq 1$.
Since $V(G)=X_0 \cup Y_0$ and $\epsilon_0 \geq 1$, 
$$\lvert V(G) \rvert = \lvert X_0 \rvert + \lvert Y_0 \rvert \leq \eta^{(r+1)^{(r+1)^\ell}}\lvert Y_0 \rvert^{\epsilon_0} + \lvert Y_0 \rvert \leq (1+\eta^{(r+1)^{(r+1)^\ell}})\lvert Y_0 \rvert^{\epsilon_0} = c\lvert Y_0 \rvert^{\epsilon_0}.$$
So Statement 3 holds.
\end{pf}

\begin{lemma} \label{island_1}
For any positive integers $r,t$, nonnegative integer $\ell$, real numbers $k>0,k'>0,\beta \geq 0$, there exists a real number $d=d(r,t,\ell,k,k',\beta)$ such that for any graph $G$ and real numbers $\epsilon \geq 0$ and $1 \geq \epsilon' \geq 0$, either
	\begin{enumerate}
		\item $\lvert E(G) \rvert > k\lvert V(G) \rvert^{1+\epsilon}$, 
		\item there exists a graph $H$ with $\lvert E(H) \rvert > k'\lvert E(H) \rvert^{1+\epsilon'}$ such that some subgraph of $G$ is isomorphic to a $[2\ell+1]$-subdivision of $H$, 
		\item there exists a $(1,\ell r+1)$-model for an $\ell$-shallow $K_{r,t}$-minor in $G$, or 
		\item there exist $X,Z \subseteq V(G)$ with $Z \subseteq X$ and $\lvert Z \rvert > \beta\lvert V(G)-X \rvert^{1+\epsilon'(r-1)}$ such that
			\begin{enumerate}
				\item every vertex in $X$ has degree at most $d \lvert V(G) \rvert^{1+\epsilon-\frac{1}{\epsilon_0}}$ in $G$, where $\epsilon_0=(1+\epsilon'(r-1))^{(r+1)^\ell}$, 
				\item for any distinct $z,z' \in Z$, the distance in $G[X]$ between $z,z'$ is at least $\ell+1$, 
				\item for any $z \in Z$ and any $u \in X$ whose distance from $z$ in $G[X]$ is at most $\ell$, $\lvert N_G(u) -X \rvert  \leq r-1$, and 
				\item $\lvert N_G(N_{G[X]}^{\leq \ell-1}[z])-X \rvert \leq r-1$ for every $z \in Z$.
			\end{enumerate}
	\end{enumerate}
\end{lemma}

\begin{pf}
Let $r,t$ be positive integers, $\ell$ be a nonnegative integer, $k,k'$ be positive real numbers, and $\beta$ be a nonnegative real number.
Let $\alpha=\alpha_{\ref{r-adherent size}}(r,t,k')$, where $\alpha_{\ref{r-adherent size}}$ is the number $\alpha$ mentioned in Lemma \ref{r-adherent size}. 
Let $c=c_{\ref{island_1_prep}}(r,t,\ell,k',\beta)$, where $c_{\ref{island_1_prep}}$ is the number $c$ mentioned in Lemma \ref{island_1_prep}.
Define $d=2kc$.

Let $G$ be a graph.
Let $\epsilon$ and $\epsilon'$ be nonnegative real numbers with $\epsilon' \leq 1$.
Let $\epsilon_0=(1+\epsilon'(r-1))^{(r+1)^\ell}$.
Let $Y_0$ be the set of all vertices of $G$ of degree greater than $d \lvert V(G) \rvert^{1+\epsilon-\frac{1}{\epsilon_0}}$.
If $\lvert V(G) \rvert \leq c\lvert Y_0 \rvert^{\epsilon_0}$, then since every vertex in $Y_0$ has degree greater than $d\lvert V(G) \rvert^{1+\epsilon-\frac{1}{\epsilon_0}}$ and $\epsilon_0 \geq 1$, 
$$2\lvert E(G) \rvert > d\lvert V(G) \rvert^{1+\epsilon-\frac{1}{\epsilon_0}} \cdot \lvert Y_0 \rvert \geq d\lvert V(G) \rvert^{1+\epsilon-\frac{1}{\epsilon_0}} \cdot (\frac{\lvert V(G) \rvert}{c})^{\frac{1}{\epsilon_0}} \geq \frac{d}{c} \lvert V(G) \rvert^{1+\epsilon} = 2k \lvert V(G) \rvert^{1+\epsilon},$$
so Statement 1 holds and we are done.
Hence we may assume that $\lvert V(G) \rvert > c\lvert Y_0 \rvert^{\epsilon_0}$ and Statements 1-3 do not hold.
So by Lemma \ref{island_1_prep}, there exist $X,Z \subseteq V(G)-Y_0$ with $Z \subseteq X$ and $\lvert Z \rvert > \beta\lvert V(G)-X \rvert^{1+\epsilon'(r-1)}$ such that Statements 4(b)-4(d) hold.
Since $X \subseteq V(G)-Y_0$, Statement 4(a) holds.
This proves the lemma.
\end{pf}

\begin{lemma} \label{island_2}
For any positive integers $r,t,t'$, nonnegative integer $\ell$, and positive real numbers $k,k'$, there exist real numbers $d=d(r,t,t',\ell,k,k')$ and $N=N(r,t,t',\ell,k,k')$ and a positive integer $r_0=r_0(r,\ell)$ such that for any graph $G$ and real numbers $\epsilon \geq 0$ and $1 \geq \epsilon' \geq 0$, either
	\begin{enumerate}
		\item $\lvert E(G) \rvert > k\lvert V(G) \rvert^{1+\epsilon}$, 
		\item there exists a graph $H$ with $\lvert E(H) \rvert > k'\lvert V(H) \rvert^{1+\epsilon'}$ such that some subgraph of $G$ is isomorphic to a $[2\max\{\ell,2\ell-2\}+1]$-subdivision of $H$, 
		\item there exists a $(1,\max\{\ell,2\ell-2\} r+1)$-model for a $\max\{\ell,2\ell-2\}$-shallow $K_{r,t}$-minor in $G$,
		\item there exist $X,Z,W \subseteq V(G)$ with $Z \subseteq X$, $\lvert Z \rvert =t'$, $W \subseteq V(G)-X$ and $\lvert W \rvert \leq r-1$ such that
			\begin{enumerate}
				\item every vertex in $X$ has degree at most $d \lvert V(G) \rvert^{1+\epsilon-\frac{1}{\epsilon_0}}$ in $G$, where $\epsilon_0=(1+(r-1)\epsilon')^{r_0}$, 
				\item for any distinct $z,z' \in Z$, the distance in $G[X]$ between $z,z'$ is at least $\max\{\ell,2\ell-2\}+1$, and
				\item $N_G(N_{G[X]}^{\leq \max\{\ell-1,0\}}[z])-X =W$ for every $z \in Z$,
			\end{enumerate}
		\item $(1+\epsilon-\frac{1}{\epsilon_0})(2\max\{\ell,2\ell-2\}+1) > \frac{1}{2}$, or 
		\item $\lvert V(G) \rvert \leq N$. 
	\end{enumerate}
\end{lemma}

\begin{pf}
Let $r,t,t'$ be positive integers, $\ell$ be a nonnegative integer, and let $k,k'$ be positive real numbers.
Let $\alpha=\max_{p \in [r]}\alpha_{\ref{r-adherent size}}(p,t',k')$, where $\alpha_{\ref{r-adherent size}}$ is the number $\alpha$ mentioned in Lemma \ref{r-adherent size}. 
Define $d=d_{\ref{island_1}}(r,t,\max\{\ell,2\ell-2\},k,k',r\alpha)$, where $d_{\ref{island_1}}$ is the real number $d$ mentioned in Lemma \ref{island_1}. 
Let $\alpha_0=\alpha_{\ref{r-adherent size}}(r,t,k')$.
Let $\gamma=\frac{t'}{\alpha}$.
Define $N=(\alpha_0\gamma+rt'+\gamma)^2 d^{4\max\{\ell,2\ell-2\}+2}$.
Define $r_0=(r+1)^{\max\{\ell,2\ell-2\}}$.

Let $G$ be a graph.
Let $\epsilon$ and $\epsilon'$ be nonnegative real numbers with $\epsilon' \leq 1$.
We may assume that Statements 1-3 do not hold, for otherwise we are done.
By Lemma \ref{island_1}, there exist $X,Z_0 \subseteq V(G)$ with $Z_0 \subseteq X$ and $\lvert Z_0 \rvert > r\alpha \cdot \lvert V(G)-X \rvert^{1+\epsilon'(r-1)}$ such that
	\begin{itemize}
		\item[(i)] every vertex in $X$ has degree at most $d\lvert V(G) \rvert^{1+\epsilon-\frac{1}{\epsilon_0}}$ in $G$, where $\epsilon_0=(1+\epsilon'(r-1))^{(r+1)^{\max\{\ell,2\ell-2\}}}=(1+\epsilon'(r-1))^{r_0}$,
		\item[(ii)] for any distinct $z,z' \in Z_0$, the distance in $G[X]$ between $z,z'$ is at least $\max\{\ell,2\ell-2\}+1$, 
		\item[(iii)] for every $z \in Z_0$ and $u \in X$ whose distance from $z$ in $G[X]$ is at most $\max\{\ell,2\ell-2\}$, $\lvert N_G(u) -X \rvert  \leq r-1$, and 
		\item[(iv)] $\lvert N_G(N_{G[X]}^{\leq \max\{\ell,2\ell-2\}-1}[z])-X \rvert \leq r-1$ for every $z \in Z_0$.
	\end{itemize}
We assume that $Z_0$ is maximal subject to the above properties. 

\medskip

\noindent{\bf Claim 1:} $\lvert N_G(N_{G[X]}^{\leq \max\{\ell-1,0\}}[z])-X \rvert \leq r-1$ for every $z \in Z_0$.

\noindent{\bf Proof of Claim 1:}
If $\ell-1 \geq 0$, then $\max\{\ell-1,0\}=\ell-1 \leq \max\{\ell,2\ell-2\}-1$, so $\lvert N_G(N_{G[X]}^{\leq \max\{\ell-1,0\}}[z])-X \rvert \leq \lvert N_G(N_{G[X]}^{\leq \max\{\ell,2\ell-2\}-1}[z])-X \rvert \leq r-1$ for every $z \in Z_0$ by (iv).
If $\ell-1<0$, then for every $z \in Z_0$, $N_{G[X]}^{\leq \max\{\ell-1,0\}}[z] = N_{G[X]}^{\leq 0}[z]=\{z\}$, so $\lvert N_G(N_{G[X]}^{\leq \max\{\ell-1,0\}}[z])-X \rvert = \lvert N_{G}(z)-X \rvert \leq r-1$ by (iii).
$\Box$

\medskip

By Claim 1 and piegeonhole principle, there exist $p \in [r-1] \cup \{0\}$ and $Z_1 \subseteq Z_0$ with $\lvert Z_1 \rvert \geq \frac{1}{r}\lvert Z_0 \rvert > \alpha\lvert V(G)-X \rvert^{1+\epsilon'(r-1)}$ such that $\lvert N_G(N_{G[X]}^{\leq \max\{\ell-1,0\}}[z])-X \rvert=p$ for every $z \in Z_1$.

We first assume that $p \in [r-1]$.
In particular, $r-1 \geq p \geq 1$ and $\alpha \geq \alpha_{\ref{r-adherent size}}(p,t',k')$.
Let $Y=V(G)-X$.
Let $\C = \{G[N_{G[X]}^{\leq \max\{\ell-1,0\}}[z]]: z \in Z_1\}$.
By (ii), $\C$ is a collection of pairwise disjoint connected subgraphs of $G-Y$ of radius at most $\max\{\ell-1,0\}$.
By (i), each member of $\C$ has at most $d^\ell\lvert V(G) \rvert^{(1+\epsilon-\frac{1}{\epsilon_0})\ell}$ vertices.
Furthermore, each member of $\C$ is $p$-adherent to $Y$.
Note that $\lvert \C \rvert=\lvert Z_1 \rvert>\alpha\lvert Y \rvert^{1+\epsilon'(r-1)}$.
Since Statement 2 of this lemma does not hold, by Lemma \ref{r-adherent size}, there exists a $(1,d^\ell\lvert V(G) \rvert^{(1+\epsilon-\frac{1}{\epsilon_0})\ell},\{G[\{y\}]: y \in Y\},\C)$-model for a $\max\{\ell-1,0\}$-shallow $K_{p,t'}$-minor in $G$.
So there exist $W \subseteq Y$ with $\lvert W \rvert=p$ and a subset $\C'$ of $\C$ with $\lvert \C' \rvert=t'$ and a $(1,d^\ell\lvert V(G) \rvert^{(1+\epsilon-\frac{1}{\epsilon_0})\ell},\{G[\{y\}]: y \in W\},\C')$-model for a $\max\{\ell-1,0\}$-shallow $K_{p,t'}$-minor in $G$.
Let $Z$ be the subset of $Z_1$ such that $\C' = \{G[N_{G[X]}^{\leq \max\{\ell-1,0\}}[z]]: z \in Z\}$. 
Then Statement 4 holds.

So we may assume that $p=0$.
If $\lvert Z_1 \rvert \geq t'$, then Statement 4 holds by choosing $W=\emptyset$ and $Z$ to be a subset of $Z_1$ with size $t'$.
So we may assume that $\lvert Z_1 \rvert < t'$.
Hence $\lvert Z_0 \rvert \leq r\lvert Z_1 \rvert<rt'$.
Since $\lvert Z_0 \rvert > r\alpha\lvert V(G)-X \rvert^{1+\epsilon'(r-1)}$, $\lvert V(G)-X \rvert^{1+\epsilon'(r-1)} \leq \frac{t'}{\alpha} = \gamma$.

Let $X_1 = X-\bigcup_{z \in Z_0}N_{G[X]}^{\leq \max\{\ell,2\ell-2\}}[z]$.
By the maximality of $Z_0$, for every $x \in X_1$, $\lvert N_G(N_{G[X]}^{\leq \max\{\ell,2\ell-2\}}[x])-X \rvert \geq r$.
So for each $x \in X_1$, there exists a subgraph $L_x$ of $G[N_{G[X]}^{\leq \max\{\ell,2\ell-2\}}[x]]$ consisting of a union of at most $r$ paths of length at most $\max\{\ell,2\ell-2\}$ from $x$ in $G[N_{G[X]}^{\leq \max\{\ell,2\ell-2\}}[x]]$ such that $L_x$ is $r$-adherent to $Y$. 
Note that each $L_x$ contains at most $r \cdot \max\{\ell,2\ell-2\}+1$ vertices.
Let $X_2$ be a maximal subset of $X_1$ such that the distance in $G[X]$ between any two distinct elements in $X_2$ is at least $2 \cdot \max\{\ell,2\ell-2\}+1$.
Let $\C_2 = \{L_x: x \in X_2\}$.
Note that members of $\C_2$ are pairwise disjoint $r$-adherent subgraphs of $G-Y$ on at most $\max\{\ell,2\ell-2\}r+1$ vertices of radius at most $\max\{\ell,2\ell-2\}$.
Since Statements 2 and 3 of this lemma do not hold, by Lemma \ref{r-adherent size}, $\lvert \C_2 \rvert \leq \alpha_0\lvert Y \rvert^{1+\epsilon'(r-1)} = \alpha_0\lvert V(G)-X \rvert^{1+\epsilon'(r-1)} \leq \alpha_0\gamma$.

Let $\ell_0=\max\{\ell,2\ell-2\}$.
By the maximality of $X_2$, $X_1 \subseteq \bigcup_{x \in X_2}N_{G[X]}^{\leq 2\max\{\ell,2\ell-2\}}[x] = \bigcup_{x \in X_2}N_{G[X]}^{\leq 2\ell_0}[x]$.
By (i), 
\begin{align*}
	\lvert X_1 \rvert \leq & \lvert X_2 \rvert \cdot \sum_{i=0}^{2\ell_0}(d\lvert V(G) \rvert^{1+\epsilon-\frac{1}{\epsilon_0}})^i \\
	\leq & \lvert X_2 \rvert \cdot d^{2\ell_0+1}\lvert V(G) \rvert^{(1+\epsilon-\frac{1}{\epsilon_0})(2\ell_0+1)} \\
	= & \lvert \C_2 \rvert \cdot d^{2\ell_0+1}\lvert V(G) \rvert^{(1+\epsilon-\frac{1}{\epsilon_0})(2\ell_0+1)} \\
	\leq & \alpha_0\gamma d^{2\ell_0+1}\lvert V(G) \rvert^{(1+\epsilon-\frac{1}{\epsilon_0})(2\ell_0+1)}.
\end{align*}
Since $X_1 = X-\bigcup_{z \in Z_0}N_{G[X]}^{\leq \max\{\ell,2\ell-2\}}[z] = X-\bigcup_{z \in Z_0}N_{G[X]}^{\leq \ell_0}[z]$, 
\begin{align*}
	\lvert X \rvert \leq & \lvert X_1 \rvert + \lvert \bigcup_{z \in Z_0}N_{G[X]}^{\leq \ell_0}[z] \rvert \\
	\leq & \alpha_0\gamma d^{2\ell_0+1}\lvert V(G) \rvert^{(1+\epsilon-\frac{1}{\epsilon_0})(2\ell_0+1)} + \lvert Z_0 \rvert d^{\ell_0+1}\lvert V(G) \rvert^{(1+\epsilon-\frac{1}{\epsilon_0})(\ell_0+1)} \\
	 \leq & (\alpha_0\gamma+\lvert Z_0 \rvert) d^{2\ell_0+1}\lvert V(G) \rvert^{(1+\epsilon-\frac{1}{\epsilon_0})(2\ell_0+1)} \\
	< & (\alpha_0\gamma+rt') d^{2\ell_0+1}\lvert V(G) \rvert^{(1+\epsilon-\frac{1}{\epsilon_0})(2\ell_0+1)}.
\end{align*}
Let $\epsilon_1 = (1+\epsilon-\frac{1}{\epsilon_0})(2\ell_0+1)$.
Hence $\lvert X \rvert <(\alpha_0\gamma+rt') d^{2\ell_0+1}\lvert V(G) \rvert^{\epsilon_1}$.

We may assume that $\epsilon_1 \leq \frac{1}{2}$, for otherwise Statement 5 holds.
Since $\epsilon_0 \geq 1$, $0 \leq \epsilon_1 \leq \frac{1}{2}$.
Hence 
\begin{align*}
	\lvert V(G) \rvert = & \lvert X \rvert + \lvert V(G)-X \rvert \\
	\leq & \lvert X \rvert + \lvert V(G)-X \rvert^{1+\epsilon'(r-1)} \\
	\leq & (\alpha_0\gamma+rt') d^{2\ell_0+1}\lvert V(G) \rvert^{\epsilon_1} + \gamma \\
	\leq & (\alpha_0\gamma+rt'+\gamma) d^{2\ell_0+1}\lvert V(G) \rvert^{\epsilon_1} \\
	\leq & \sqrt{N}\lvert V(G) \rvert^{\frac{1}{2}}.
\end{align*}
Therefore, $\lvert V(G) \rvert \leq N$, so Statement 6 holds.
\end{pf}

\bigskip

Recall that a set of homomorphisms is consistent if it satisfies (CON1) and (CON2) as stated in Section \ref{subsec:consistent}.

\begin{lemma} \label{island_3}
For any positive integers $r,t,t',a$, nonnegative integer $\ell$, and positive real numbers $k,k'$, there exist real numbers $d=d(r,t,t',a,\ell,k,k'),N=N(r,t,t',a,\ell,k,k')$ and a positive integer $r_0=r_0(r,\ell)$ such that for any graph $G$ and real numbers $\epsilon \geq 0$ and $1 \geq \epsilon' \geq 0$, either
	\begin{enumerate}
		\item $\lvert E(G) \rvert > k\lvert V(G) \rvert^{1+\epsilon}$, 
		\item there exists a graph $H$ with $\lvert E(H) \rvert>k'\lvert V(H) \rvert^{1+\epsilon'}$ such that some subgraph of $G$ is isomorphic to a $[2\max\{\ell,2\ell-2\}+1]$-subdivision of $H$,
		\item there exists a $(1,\max\{\ell,2\ell-2\} r+1)$-model for a $\max\{\ell,2\ell-2\}$-shallow $K_{r,t}$-minor in $G$, 
		\item $(1+\epsilon-\frac{1}{\epsilon_0})(2\max\{\ell,2\ell-2\}+1) > \frac{1}{2}$, where $\epsilon_0=(1+\epsilon'(r-1))^{r_0}$, 
		\item $\lvert V(G) \rvert \leq N$, or 
		\item there exist $X,Z,W \subseteq V(G)$ with $Z \subseteq X$, $\lvert Z \rvert =t'$, $W \subseteq V(G)-X$ with $\lvert W \rvert \leq r-1$ such that
			\begin{enumerate}
				\item every vertex in $X$ has degree at most $d\lvert V(G) \rvert^{1+\epsilon-\frac{1}{\epsilon_0}}$ in $G$, 
				\item for any distinct $z,z' \in Z$, the distance in $G[X]$ between $z,z'$ is at least $\max\{\ell,2\ell-2\}+1$,
				\item $N_G(N_{G[X]}^{\leq \max\{\ell-1,0\}}[z])-X =W$ for every $z \in Z$, and
				\item for any quasi-order $Q=(U,\preceq)$ with $\lvert U \rvert \leq a$, legal $Q$-labelling $f_G$ of $G$, set $\HH$ of homomorphisms from an element in $\{(L,f_L): (L,f_L)$ is a legal $Q$-labelled graph with $V(L) \subseteq [\ell]\}$ to the set of all induced subgraphs of $(G,f_G)$ satisfying (CON1) and (CON2), if there exist $z_0 \in Z$ and $\phi \in \HH$ from some $(L,f_L)$ to $(G[N_{G[X]}^{\leq \max\{\ell-1,0\}}[z_0] \cup W],f_G|_{N_{G[X]}^{\leq \max\{\ell-1,0\}}[z_0] \cup W})$ with $z_0 \in \phi(V(L))$, then for every $z \in Z$, there exists $\phi_z \in \HH$ from $(L,f_L)$ to $(G[N_{G[X]}^{\leq \max\{\ell-1,0\}}[z] \cup W], f_G|_{N_{G[X]}^{\leq \max\{\ell-1,0\}}[z] \cup W})$ such that 
					\begin{enumerate}
						\item $z \in \phi_z(V(L))$,  
						\item for every $y \in W$ and $u \in V(L)$, $y=\phi(u)$ if and only if $y=\phi_z(u)$, 
						\item for every $u \in V(L)$, $\phi(u)=z_0$ if and only if $\phi_z(u)=z$, and
						\item there exists an isomorphism $\iota_z$ from $(G[\phi(V(L)) \cup W],f|_{\phi(V(L)) \cup W})$ to \linebreak $(G[\phi_z(V(L)) \cup W],f|_{\phi_z(V(L)) \cup W})$ such that $\iota_z(\phi(u))=\phi_z(u)$ for every $u \in V(L)$, and $\iota_z(y)=y$ for every $y \in W$.
					\end{enumerate}
			\end{enumerate}
	\end{enumerate}
\end{lemma}

\begin{pf}
Let $r,t,t',a$ be positive integers, $\ell$ be a nonnegative integer, and let $k,k'$ be positive real numbers.
Let $\ell_1=2^{\ell+{\ell \choose 2}}$, $\ell_2=\ell_1 \cdot 2^{a^2} \cdot a^{2^\ell \cdot \ell !}$, and $\ell_3=\ell_2^4 \cdot \ell^\ell \cdot 2^{(2r-1)\ell}$. 
Define $d=d_{\ref{island_2}}(r,t,t'\ell_3,\ell,k,k')$ and $N=N_{\ref{island_2}}(r,t,t'\ell_3,\ell,k,k')$, where $d_{\ref{island_2}}$ and $N_{\ref{island_2}}$ are the real numbers $d$ and $N$ mentioned in Lemma \ref{island_2}, respectively.
Define $r_0=r_{\ref{island_2}}(r,\ell)$, where $r_{\ref{island_2}}$ is the integer $r_0$ mentioned in Lemma \ref{island_2}.

Let $G$ be a graph.
Let $\epsilon$ and $\epsilon'$ be nonnegative real numbers with $\epsilon' \leq 1$.
Let $\epsilon_0=(1+\epsilon'(r-1))^{r_0}$.
We may assume that Statements 1-5 do not hold, for otherwise we are done.
By Lemma \ref{island_2}, there exist $X,Z_0,W \subseteq V(G)$ with $Z_0 \subseteq X$, $\lvert Z_0 \rvert =t'\ell_3$, $W \subseteq V(G)-X$ and $\lvert W \rvert \leq r-1$ such that
	\begin{itemize}
		\item every vertex in $X$ has degree at most $d\lvert V(G) \rvert^{1+\epsilon-\frac{1}{\epsilon_0}}$ in $G$, 
		\item for any distinct $z,z' \in Z_0$, the distance in $G[X]$ between $z,z'$ is at least $\max\{\ell,2\ell-2\}+1$, and
		\item $N_G(N_{G[X]}^{\leq \max\{\ell-1,0\}}[z])-X =W$ for every $z \in Z_0$.
	\end{itemize}

Let $\F_1$ be the set of all graphs whose vertex-sets are contained in $[\ell]$.
So $\lvert \F_1 \rvert \leq \ell_1$.
Let $\Q$ be the set of all isomorphism classes of quasi-orders whose underlying set have size at most $a$.
Note that $\lvert \Q \rvert \leq 2^{a+2{a \choose 2}}=2^{a^2}$.
Let $\F_2 = \{(L,f_L): L \in \F_1, (L,f_L)$ is a $Q$-labelled graph for some $Q \in \Q\}$. 
So $\lvert \F_2 \rvert \leq \lvert \F_1 \rvert \cdot \lvert \Q \rvert \cdot a^{2^\ell \cdot \ell!} \leq \ell_2$.
Let $\HH' = \{\phi: \phi$ is a homomorphism from an element $(L,f_L)$ in $\F_2$ to an element $(L',f_{L'})$ in $\F_2$ such that there exists an isomorphism $\iota$ from $(L',f_{L'})$ to an induced subgraph of $(G,f_G)$ such that $\iota \circ \phi \in \HH\}$.
Since every graph in $\F_2$ contains at most $\ell$ vertices, $\lvert \HH' \rvert \leq \ell_2^2 \cdot \ell^\ell$.

Denote $W$ by $\{v_1,v_2,...,v_{\lvert W \rvert}\}$.
For each $z \in Z_0$, let $S_z$ be the set of all tuples $((L,f_L), \allowbreak (L',f_{L'}), \allowbreak \phi,A_0,A_1,A_2, \allowbreak ..., \allowbreak A_{\lvert W \rvert},C_1,C_2,...,C_{\lvert W \rvert})$ such that
	\begin{itemize}
		\item $(L,f_L),(L',f_{L'}) \in \F_2$, 
		\item $\phi \in \HH'$ is a homomorphism from $(L,f_L)$ to $(L',f_{L'})$, and 
		\item there exists an isomorphism $\iota$ from $(L',f_{L'})$ to an induced subgraph of $(G[N_{G[X]}^{\leq \max\{\ell-1,0\}}[z] \cup W],f_G|_{N_{G[X]}^{\leq \max\{\ell-1,0\}}[z] \cup W})$ containing $z$ such that 
			\begin{itemize}
				\item for each $j \in [\lvert W \rvert] \cup \{0\}$, $A_j=\{u \in V(L) \subseteq [\ell]: \iota(\phi(u))=v_j\}$, where $v_0=z$, and
				\item for each $j \in [\lvert W \rvert]$, $C_j = \{u \in V(L) \subseteq [\ell]: \iota(\phi(u))v_j \in E(G)\}$.
			\end{itemize}
	\end{itemize}
Note that the definition of $A_0,A_1,...,A_{\lvert W \rvert},C_1,...,C_{\lvert W \rvert}$ depends on the choice of $\iota$, but each $A_j$ or $C_j$ is a subset of $[\ell]$, so there are at most $2^\ell$ possibilities for $A_j$ or $C_j$ for each fixed $j \in [\lvert W \rvert] \cup \{0\}$.
Hence there are at most $\lvert \F_2 \rvert \cdot \lvert \F_2 \rvert \cdot \lvert \HH' \rvert \cdot ({2^\ell})^{2\lvert W \rvert+1} \leq \ell_2^2 \cdot \ell_2^2\ell^\ell \cdot 2^{(2r-1)\ell}=\ell_3$ possibilities for the sets $S_z$ among all $z \in Z_0$.
So there exists $Z \subseteq Z_0$ with $\lvert Z \rvert \geq \frac{1}{\ell_3}\lvert Z_0 \rvert \geq t'$ such that for any $z_1,z_2 \in Z$, $S_{z_1}=S_{z_2}$.
Then Statement 6 holds, since $\HH$ satisfies (CON1) and (CON2).
\end{pf}

\begin{lemma} \label{homogen_0}
For any positive integers $r,t,c,m,p$, there exists a positive integer $s=s(r,t,c, \allowbreak m,p)$ such that the following holds.
Let $G$ be a bipartite graph with a bipartition $\{A,B\}$ such that every edge of $G$ is colored with one element in $[c]$.
Assume that there exists no subgraph of $G$ isomorphic to $K_{r,t}$ such that the part consisting of $t$ vertices is a subset of $A$.
If $\lvert A \rvert \geq s$, then there exists a subset $A'$ of $A$ with $\lvert A' \rvert \geq m$ such that for every subset $T$ of $B$ with $\lvert T \rvert \leq p$, there exists $M_T \subseteq A'$ with $\lvert M_T \rvert \geq m$ such that for every $v \in T$, either there exists $c_v \in [c]$ such that $v$ is adjacent in $G$ to all vertices in $A'$ by using edges with color $c_v$, or $v$ is non-adjacent in $G$ to all vertices in $M_T$.
\end{lemma}

\begin{pf}
Let $r,t,c,m,p$ be positive integers.
Let $s_0=t$.
For every positive integer $i$, let $s_{i}=mpcs_{i-1}$.
Define $s=s_{r}$.

We shall prove this lemma by induction on $r$.

Let $G$ be a graph stated in this lemma.
We may assume that there exists a subset $T_0$ of $B$ with $\lvert T_0 \rvert \leq p$ such that for every subset $S$ of $A$ with $\lvert S \rvert=m$, there exists an edge $e_{T_0,S}$ of $G$ between $T_0$ and $S$, for otherwise we are done by choosing $A'=A$.
Since there are ${\lvert A \rvert \choose m}$ subsets of $A$ with size $m$, and for each edge $e$ of $G$ incident with $T_0$, there are at most ${\lvert A \rvert-1 \choose m-1}$ subsets $S$ of $A$ with $\lvert S \rvert=m$ such that $e_{T_0,S}=e$, we know that there are at least ${\lvert A \rvert \choose m}/{\lvert A \rvert-1 \choose m-1} = \frac{\lvert A \rvert}{m} \geq pcs_{r-1} \geq \lvert T_0 \rvert cs_{r-1}$ edges incident with $T_0$.
So there exist $v^* \in T_0$, $c_{v^*} \in [c]$ and $A_1 \subseteq A$ with $\lvert A_1 \rvert=s_{r-1}$ such that $v^*$ is adjacent to all vertices in $A_1$ by using edges with color $c_{v^*}$.
Since $s_{r-1} \geq s_0 \geq t$, there exists a $K_{1,t}$ subgraph whose part consisting of $t$ vertices is contained in $A$, so $r \geq 2$, and the lemma holds when $r=1$.

Let $G'=G[A_1 \cup (B-\{v^*\})]$.
Note that $G'$ has no $K_{r-1,t}$ subgraph whose part consisting of $t$ vertices is contained in $A_1$.
Since $\lvert A_1 \rvert \geq s_{r-1}$, by the induction hypothesis, there exists a subset $A'$ of $A_1$ with $\lvert A_1' \rvert \geq m$ such that for every subset $T$ of $B-\{v^*\}$ with $\lvert T \rvert \leq p$, there exists $M_T \subseteq A'$ with $\lvert M_T \rvert \geq m$ such that for every $v \in T$, either there exists $c_v \in [c]$ such that $v$ is adjacent in $G' \subseteq G$ to all vertices in $A'$ by using edges with color $c_v$, or $v$ is non-adjacent in $G'$ (and hence in $G$) to all vertices in $M_T$.
For every subset $T$ of $B$ with $\lvert T \rvert \leq p$, if $v^* \not \in T$, then $T \subseteq B-\{v^*\}$, so $M_T$ is desired; if $v^* \in T$, then let $M_T=M_{T-\{v^*\}}$, and $M_T$ is desired since $v^*$ is adjacent to all vertices in $A_1 \supseteq A'$ by using edges with color $c_{v^*}$.
This proves the lemma.
\end{pf}

\begin{lemma} \label{island_4}
For any positive integers $r,t,t',a,m,h$, nonnegative integer $\ell$, and positive real numbers $k,k'$, there exist positive real numbers $d=d(r,t,t',a,m,h,\ell,k,k'),N=N(r,t,t',a, \allowbreak m, h, \ell,k,k'), \allowbreak b=b(r,t,t',a,m,h,\ell,k,k'),r_0=r_0(r,\ell),r_1=r_1(r,\ell)$ such that if $G$ is a graph and $\epsilon,\epsilon'$ are nonnegative real numbers with $\epsilon \leq \frac{1}{8\max\{\ell+1,2\ell\}+4}$ and $\epsilon' \leq r_1$ such that
	\begin{enumerate}
		\item for every induced subgraph $L$ of $G$, $\lvert E(L) \rvert \leq k\lvert V(L) \rvert^{1+\epsilon}$, 
		\item for every graph $L$ for which some subgraph of $G$ is isomorphic to a $[2\max\{\ell+1,2\ell\}+1]$-subdivision of $L$, $\lvert E(L) \rvert \leq k'\lvert V(L) \rvert^{1+\epsilon'}$, and 
		\item there exists no $(1,\max\{\ell+1,2\ell\} r+1)$-model for a $\max\{\ell+1,2\ell\}$-shallow $K_{r,t}$-minor in $G$, 
	\end{enumerate}
then for any quasi-order $Q=(U,\preceq)$ with $\lvert U \rvert \leq a$ and legal $Q$-labelling $f$ of $G$, there exist an ordering $\sigma$ of $G$ and a sequence $(S_i: i \in [\lvert V(G) \rvert])$ such that for every $i \in [\lvert V(G) \rvert]$,
	\begin{enumerate}
		\item $S_i$ is a subset of $\{v \in V(G): \sigma(v) \geq i+1\}$ with size at most $r-1$, and
		\item if $i \in [\lvert V(G) \rvert-N]$, then the $\max\{\ell-1,0\}$-basin $B_i$ with respect to $\sigma$ and $(S_j: j \in [\lvert V(G) \rvert])$ at $i$ satisfies the following: 
			\begin{enumerate}
				\item $\lvert B_i \rvert \leq b \lvert V(G) \rvert^{(1+\epsilon-\frac{1}{\epsilon_0})\ell}$, where $\epsilon_0=(1+\epsilon'(r-1))^{r_0}$.
				\item Every vertex in $B_i$ has degree at most $d\lvert V(G) \rvert^{1+\epsilon-\frac{1}{\epsilon_0}}$ in $G_{\sigma,\geq i}$.
				\item There exists a set $Z^*=\{z_1,z_2,...,z_{\lvert Z^* \rvert}\}$ of at least $t'$ vertices in $\{v \in V(G): \sigma(v) \geq i+1\}$.
				\item There exist $\lvert Z^* \rvert$ pairwise disjoint subsets $Z'_1,Z'_2,...,Z'_{\lvert Z^* \rvert}$ of $\{v \in V(G): \sigma(v) \geq i+1\}-S_i$ disjoint from $N^{\leq \min\{\ell,1\}}_{G_{\sigma, \geq i}}[B_i]$ such that for any set $\HH$ of homomorphisms from an element in $\{(L,f_L): (L,f_L)$ is a legal $Q$-labelled graph with $V(L) \subseteq [\ell+1]\}$ to the set of all induced subgraphs of $(G_{\sigma, \geq i},f|_{V(G_{\sigma, \geq i})})$ satisfying (CON1) and (CON2), if $(L,f_L)$ is a legal $Q$-labelled graph such that there exists $\phi \in \HH$ from $(L,f_L)$ to $(G_{\sigma, \geq i}[B_i \cup S_i],f|_{B_i \cup S_i})$ with $\sigma^{-1}(i) \in \phi(V(L))$, then 
					\begin{enumerate}
						\item for every $j \in [\lvert Z^* \rvert]$, there exists $\phi_j \in \HH$ from $(L,f_L)$ to $(G_{\sigma,\geq i}[Z'_j \cup S_i],f|_{Z'_j \cup S_i})$ such that 
							\begin{itemize}
								\item $z_j \in \phi_j(V(L))$,  
								\item for every $y \in S_i$ and $u \in V(L)$, $y=\phi(u)$ if and only if $y=\phi_j(u)$,
								\item for every $u \in V(L)$, $\phi(u)=\sigma^{-1}(i)$ if and only if $\phi_j(u)=z_j$, and
								\item there exists an isomorphism $\iota_j$ from $(G_{\sigma, \geq i}[\phi(V(L)) \cup S_i],f|_{\phi(V(L)) \cup S_i})$ to \linebreak $(G_{\sigma, \geq i}[\phi_j(V(L)) \cup S_i], \allowbreak f|_{\phi_j(V(L)) \cup S_i})$ such that $\iota_j(\phi(u))=\phi_j(u)$ for every $u \in V(L)$, and $\iota_j(y)=y$ for every $y \in S_i$.
							\end{itemize}
						\item for any $I \subseteq [i-1]$ with $\lvert I \rvert \leq h$, there exists $M_{I} \subseteq Z^*$ with $\lvert M_{I} \rvert \geq m$ such that for every $v \in \sigma^{-1}(I)$ and $z_j \in M_I$, $\{u \in V(L): v\phi(u) \in E(G)\} \supseteq \{u \in V(L): v\phi_j(u) \in E(G)\}$.
					\end{enumerate}
				\item There exist $\lvert Z^* \rvert$ sets $Z_1,Z_2,...,Z_{\lvert Z^* \rvert}$ such that for every $j \in [\lvert Z^* \rvert]$, 
					\begin{itemize}
						\item $N_{G_{\sigma, \geq i}-S_i}^{\leq \max\{\ell-1,0\}}[z_j] \subseteq Z_j \subseteq Z'_j$, and 
						\item if $\ell \geq 1$, then $N_{G_{\sigma, \geq i}}[Z_j] \subseteq Z'_j \cup S_i$. 
					\end{itemize}
			\end{enumerate}
	\end{enumerate}
\end{lemma}

\begin{pf}
Let $r,t,t',a,m,h$ be positive integers, $\ell$ be a nonnegative integer, and let $k,k'$ be positive real numbers.
Let $\L$ be the collection of $Q$-labelled graphs with vertex-set contained in $[\ell+1]$, among all isomorphism classes of quasi-orders $Q$ with ground set size at most $a$.
Note that $\lvert \L \rvert$ is finite and only depends on $\ell$ and $a$. 
So there exists an integer $\ell^*$ such that there are at most $\ell^*$ homomorphisms from a member of $\L$ to a member of $\L$.
Let $s= s_{\ref{homogen_0}}(r,t,\lvert \L \rvert\ell^*(\ell+1),m+t'+1,h)$, where $s_{\ref{homogen_0}}$ is the integer $s$ mentioned in Lemma \ref{homogen_0}.
Define $d=d_{\ref{island_3}}(r,t,s,a,\ell+1,k,k')$, $N=N_{\ref{island_3}}(r,t,s,a,\ell+1,k,k')$, $r_0=r_{\ref{island_3}}(r,\ell+1)$, where $d_{\ref{island_3}}, N_{\ref{island_3}}, r_{\ref{island_3}}$ are the numbers $d,N,r_0$ mentioned in Lemma \ref{island_3}, respectively.
Define $b=d^\ell$ and $r_1 = \frac{(\frac{1}{1-\frac{1}{8\max\{\ell+1,2\ell\}+4}})^{\frac{1}{r_0}}-1}{r}$.

Let $G$ be a graph and $\epsilon,\epsilon'$ be real numbers as stated in the lemma.
Let $\epsilon_0=(1+\epsilon'(r-1))^{r_0}$ and $\epsilon_1 = (1+\epsilon-\frac{1}{\epsilon_0})(2\max\{\ell+1,2\ell\}+1)$.
Since $\epsilon' \leq r_1$, $\epsilon_0 \leq \frac{1}{1-\frac{1}{8\max\{\ell+1,2\ell\}+4}}$, so $\epsilon_1 \leq (\frac{1}{8\max\{\ell+1,2\ell\}+4}+\epsilon)(2\max\{\ell+1,2\ell\}+1)$.
Since $\epsilon \leq \frac{1}{8\max\{\ell+1,2\ell\}+4}$, $\epsilon_1 \leq \frac{1}{2}$.

Let $Q$ be a quasi-order and $f$ a legal $Q$-labelling as stated in the lemma.

Now we define the ordering $\sigma$ and the sequence $(S_j: j \in [\lvert V(G) \rvert])$ and show that they satisfy the conclusions of this lemma.
To define $\sigma$, it suffices to define its inverse function from $[\lvert V(G) \rvert]$ to $V(G)$.

Let $i \in [\lvert V(G) \rvert]$.
Assume that $\sigma^{-1}(j)$ and $S_j$ are defined for every $j \in [i-1]$.
For each $j \in [i-1]$, we denote the vertex $v$ with $\sigma(v)=j$ by $v_j$.
Let $G_i = G-\{v_j: j \in [i-1]\}$.
So $G_{\sigma,\geq i}=G_i$ no matter how we further define $\sigma^{-1}(i),\sigma^{-1}(i+1),...,\sigma^{-1}(\lvert V(G) \rvert)$.
If $i \geq \lvert V(G) \rvert-N+1$, then define $\sigma^{-1}(i)$ to be an arbitrary vertex in $G_i$ and define $S_i=\emptyset$.
If $i \in [\lvert V(G) \rvert-N]$, then $\lvert V(G_i) \rvert \geq N+1$, and since $\epsilon_1 \leq \frac{1}{2}$ and $G_i$ is an induced subgraph of $G$, by applying Lemma \ref{island_3} to $G_i$, we know that there exist $X,Z,W \subseteq V(G_i)$ with $Z \subseteq X$, $\lvert Z \rvert =s$, $W \subseteq V(G_i)-X$ and $\lvert W \rvert \leq r-1$ such that
			\begin{itemize}
				\item[(i)] every vertex in $X$ has degree at most $d \lvert V(G) \rvert^{1+\epsilon-\frac{1}{\epsilon_0}}$ in $G_i$, 
				\item[(ii)] for any distinct $z,z' \in Z$, the distance in $G_i[X]$ between $z,z'$ is at least $\max\{\ell+1,2\ell\}+1$,
				\item[(iii)] $N_{G_i}(N_{G_i[X]}^{\leq \ell}[z])-X =W$ for every $z \in Z$, and 
				\item[(iv)] for any set $\HH$ of homomorphisms from an element in $\{(L,f_L): (L,f_L)$ is a legal $Q$-labelled graph with $V(L) \subseteq [\ell+1]\}$ to the set of all induced subgraphs of $(G_i,f|_{V(G_i)})$ satisfying (CON1) and (CON2), if there exist $z_0 \in Z$ and $\phi \in \HH$ from $(L,f_L)$ to $(G_i[N_{G_i[X]}^{\leq \ell}[z_0] \cup W],f|_{N_{G_i[X]}^{\leq \ell}[z_0] \cup W})$ with $z_0 \in \phi(V(L))$, then for every $z \in Z$, there exists $\phi_{(L,f_L),z,\HH} \in \HH$ from $(L,f_L)$ to $(G_i[N_{G_i[X]}^{\leq \ell}[z] \cup W], f|_{N_{G_i[X]}^{\leq \ell}[z] \cup W})$ such that 
					\begin{itemize}
						\item $z \in \phi_{(L,f_L),z,\HH}(V(L))$,  
						\item for every $y \in W$ and $u \in V(L)$, $y=\phi(u)$ if and only if $y=\phi_{(L,f_L),z,\HH}(u)$, 
						\item for every $u \in V(L)$, $\phi(u)=z_0$ if and only if $\phi_{(L,f_L),z,\HH}(u)=z$, and 
						\item there exists an isomorphism $\iota_{(L,f_L),z,\HH}$ from $(G_i[\phi(V(L)) \cup W],f|_{\phi(V(L)) \cup W})$ to \linebreak $(G_i[\phi_{(L,f_L),z,\HH}(V(L)) \cup W], \allowbreak f|_{\phi_{(L,f_L),z,\HH}(V(L)) \cup W})$ such that $\iota_{(L,f_L),z,\HH}(\phi(u))= \linebreak \phi_{(L,f_L),z,\HH}(u)$ for every $u \in V(L)$, and $\iota_{(L,f_L),z,\HH}(y)=y$ for every $y \in W$.
					\end{itemize}
			\end{itemize}
Note that since every $\HH$ stated in (iv) satisfies (CON1) and (CON2), there are at most $\ell^*$ different $\HH$'s stated in (iv).
For every $z \in Z$, let $\L_{i,z} = \{((L,f_L),\HH): (L,f_L) \in \L, \phi_{(L,f_L),z,\HH}$ is defined for some $\HH$ stated in (iv)$\}$.
Let $F_i$ be the bipartite graph with $V(F_i)=Z \cup \sigma^{-1}([i-1])$, and $E(F_i)=\{z\sigma^{-1}(i'): z \in Z, i' \in [i-1], \sigma^{-1}(i')$ is adjacent in $G$ to some vertex in $\bigcup_{((L,f_L),\HH) \in \L_{i,z}}\phi_{(L,f_L),z,\HH}(V(L))-W\}$.
For each edge $z\sigma^{-1}(i')$ of $F_i$, color it with $\{((L,f_L),\HH,u): ((L,f_L),\HH) \in \L_{i,z}, u \in V(L), \allowbreak \sigma^{-1}(i')\phi_{(L,f_L),z,\HH}(u) \in E(G)\}$.
So we use at most $\lvert \L \rvert\ell^* (\ell+1)$ colors.
Since there exists no $(1,\max\{\ell+1,2\ell\} r+1)$-model for a $\max\{\ell+1,2\ell\}$-shallow $K_{r,t}$-minor in $G$, $F_i$ has no $K_{r,t}$ subgraph whose part consisting of $t$ vertices is contained $Z$.
Since $\lvert Z \rvert \geq s$, by Lemma \ref{homogen_0}, there exists $Z' \subseteq Z$ with $\lvert Z' \rvert \geq m+t'+1$ such that 
	\begin{itemize}
		\item[(v)] for every subset $I \subseteq [i-1]$ with $\lvert I \rvert \leq h$, there exists $M'_I \subseteq Z'$ with $\lvert M'_I \rvert \geq m+t'+1$ such that for every $i' \in I$, either $\sigma^{-1}(i')$ is adjacent in $F_i$ to all vertices in $Z'$ using the same color, or $\sigma^{-1}(i')$ is non-adjacent in $F_i$ to all vertices in $M'_I$.
			(Note that the latter implies that $\sigma^{-1}(i')$ is not adjacent in $G$ to any vertex in $\bigcup_{((L,f_L),\HH) \in \L_{i,z}}\phi_{(L,f_L),z,\HH}(V(L))-W$.)
	\end{itemize}
If $i \in [\lvert V(G) \rvert-N]$, then let $v_i$ be a vertex in $Z'$, and define $\sigma^{-1}(i)=v_i$ and $S_i=W$.

When $i \geq \lvert V(G) \rvert-N+1$, $S_i=\emptyset$; when $i \in [\lvert V(G) \rvert-N]$, $v_i \in Z' \subseteq X$, so $S_i=W \subseteq V(G_i)-X \subseteq \{v \in V(G): \sigma(v) \geq i+1\}$ and $\lvert S_i \rvert = \lvert W \rvert \leq r-1$.
So Conclusion 1 of this lemma holds.

Now we show that Conclusion 2 of this lemma holds for $i$.
So we may assume $i \in [\lvert V(G) \rvert-N]$.
Let $p=\max\{\ell-1,0\}$.

\medskip

\noindent{\bf Claim 1:} For every $z \in Z'$, $N_{G_i[X]}^{\leq p}[z]=N_{G_i-S_i}^{\leq p}[z]$. 

\noindent{\bf Proof of Claim 1:}
Suppose to the contrary that $N_{G_i[X]}^{\leq p}[z] \neq N_{G_i-S_i}^{\leq p}[z]$.
Since $G_i-S_i \supseteq G_i[X]$, $N_{G_i-S_i}^{\leq p}[z] \supset N_{G_i[X]}^{\leq p}[z]$.
So there exists $u \in N_{G_i-S_i}^{\leq p}[z] - N_{G_i[X]}^{\leq p}[z]$.
Hence there exists a path $P$ in $G_i-S_i$ of length at most $p$ from $z$ to $u$.
We may assume that $u$ is chosen so that $P$ is as short as possible.
So $P-u$ is a path in $G_i[X]$ and $u \in N_{G_i-S_i}^{\leq p}[z]-X$. 
Since $u \not \in X$, $u \neq z$, so $V(P)-\{u\} \neq \emptyset$.
Note that $V(P)-\{u\} \subseteq N_{G_i[X]}^{\leq p}[z]$.
Since $u \in N_{G_i}(V(P)-\{u\}) \subseteq N_{G_i}(N_{G_i[X]}^{\leq p}[z])$, we have $u \in X \cup W$ by (iii).
Since $u \not \in X$, $u \in W=S_i$, so $u \not \in N_{G_i-S_i}^{\leq p}[z]$, a contradiction.
$\Box$

\medskip

Let $B_i$ be the $\max\{\ell-1,0\}$-basin with respect to $\sigma$ and $(S_j: j \in [\lvert V(G) \rvert])$ at $i$.
Note that $B_i = N_{G_i-S_i}^{\leq p}[v_i]$.
Since $v_i \in Z'$, by Claim 1, $B_i = N_{G_i[X]}^{\leq p}[v_i]$.
In particular, $B_i \subseteq X$.
Hence Conclusion 2(b) follows immediately from (i).

If $p=0$, then $\lvert B_i \rvert = 1 \leq b \lvert V(G) \rvert^{(1+\epsilon-\frac{1}{\epsilon_0})\ell}$; if $p \geq 1$, then $p=\ell-1$, and since $B_i = N_{G_i[X]}^{\leq p}[v_i]$, by (i), $\lvert B_i \rvert \leq \sum_{j=0}^{p}(d \lvert V(G) \rvert^{1+\epsilon-\frac{1}{\epsilon_0}})^j \leq (d \lvert V(G) \rvert^{1+\epsilon-\frac{1}{\epsilon_0}})^{p+1} \leq (d \lvert V(G) \rvert^{1+\epsilon-\frac{1}{\epsilon_0}})^\ell=b \lvert V(G) \rvert^{(1+\epsilon-\frac{1}{\epsilon_0})\ell}$.
So Conclusion 2(a) holds.

Let $Z^*=Z'-\{v_i\}$.
So $\lvert Z^* \rvert \geq m+ t'$.
Let $z_1,z_2,...,z_{\lvert Z^* \rvert}$ be the vertices in $Z^*$.
For every $j \in [\lvert Z^* \rvert]$, let $Z'_j=N_{G_i[X]}^{\leq \ell}[z_j]$ and $Z_j=N_{G_i[X]}^{\leq p}[z_j]$.
By Claim 1, for each $j \in [\lvert Z^* \rvert]$, $N_{G_i-S_i}^{\leq p}[z_j] = Z_j \subseteq Z'_j \subseteq \{v \in V(G): \sigma(v) \geq i+1\}-S_i$.
Furthermore, for each $j \in [\lvert Z^* \rvert]$, if $\ell \geq 1$, then $\ell = p+1$, so $N_{G_i}[Z_j]=N_{G_i}[N_{G_i[X]}^{\leq p}[z_j]] \subseteq N_{G_i[X]}^{\leq p+1}[z_j] \cup (N_{G_i}[N_{G_i[X]}^{\leq p}[z_j]]-X) \subseteq N_{G_i[X]}^{\leq \ell}[z_j] \cup S_i = Z_j' \cup S_i$ by (iii).
So Conclusions 2(c) and 2(e) hold.

When $\ell=0$, each $Z'_j$ is disjoint from $B_i=N^{\leq \min\{\ell,1\}}_{G_{\sigma, \geq i}}[B_i]$, and $Z'_1,...,Z'_{t'}$ are pairwise disjoint; when $\ell \geq 1$, $\ell+1 \leq 2\ell$, so by (ii), each $Z'_j$ is disjoint from $N_{G_{\sigma, \geq i}}[B_i]=N^{\leq \min\{\ell,1\}}_{G_{\sigma, \geq i}}[B_i]$, and $Z'_1,...,Z'_{t'}$ are pairwise disjoint.
So Conclusion 2(d)(i) immediately follows from (iv).

For every $I \in [i-1]$ with $\lvert I \rvert \leq h$, let $M_I=M'_I \cap Z^*$, so $\lvert M_I \rvert \geq \lvert M'_I \rvert-1 \geq m+t'$. 
For $v \in \sigma^{-1}(I)$ and $z_j \in M_I$, if $u \in V(L)$ such that $v\phi_{(L,f_L),\HH,z_j}(u) \in E(G)$, then either $\phi_{(L,f_L),\HH,z_j}(u) \in W=S_i$ (so $\phi(u)=\phi_{(L,f_L),\HH,z_j}(u)$ and $v\phi(u) \in E(G)$), or $\phi_j(u) \in \bigcup_{((L,f_L),\HH) \in \L_{i,z}}\phi_{(L,f_L),z,\HH}(V(L))-W$, so by (v), $v$ is adjacent in $F_i$ to all vertices in $Z'$ by using edges with the same color, and hence $v\phi(u) \in E(G)$ by the definition of the edge-coloring of $F_i$.
Hence Conclusion 2(d)(ii) holds.
This proves the lemma.
\end{pf}

\section{Homomorphism counts and shallow minors} \label{sec:homo_density}

Before considering homomorphisms, we first prove Lemma \ref{parts_collection} that will be convenient for cleaning up structures.
We need the following well-known bipartite Ramsey theorem.

\begin{lemma} \label{bi_ramsey} 
For any positive integers $a,b,c$, there exists a positive integer $R=R(a,b,c)$ such that any $c$-coloring of the edges of $K_{R,R}$ results in a monochromatic $K_{a,b}$.
\end{lemma}

\begin{lemma} \label{2_part_collection}
For any positive integers $r,h,w$, there exists a positive integer $s=s(r,h,w)$ such that the following holds.
Let $G$ be a graph with no $K_{r,r}$-subgraph.
For every $i \in [2]$, let $\C_i$ be a collection of pairwise disjoint subsets of $V(G)$ with each member having size at most $h$.
If $\lvert \C_i \rvert \geq s$ for every $i \in [2]$, then for every $i \in [2]$, there exists a subset $\C_i'$ of $\C_i$ with $\lvert \C_i' \rvert=w$ such that every member of $\C_1'$ is disjoint and non-adjacent in $G$ to every member of $\C_2'$.
\end{lemma}

\begin{pf}
Let $r,h,w$ be positive integers.
Let $m=\max\{r,h+1,w\}$.
Let $s=R_{\ref{bi_ramsey}}(m,m,h^2+2)$, where $R_{\ref{bi_ramsey}}$ is the integer $R$ mentioned in Lemma \ref{bi_ramsey}.

Let $G$, $\C_1$ and $\C_2$ be as stated in this lemma.
By taking subsets, we may assume that $\lvert \C_1 \rvert = \lvert \C_2 \rvert =s$.
Denote the members of $\C_1$ by $A_1,A_2,...,A_s$, and denote the members of $\C_2$ by $B_1,B_2,...,B_s$.
For each $i \in [s]$, we denote $A_i = \{v_{i,j}: j \in [\lvert A_i \rvert]\}$ and $B_i=\{u_{i,j}: j \in [\lvert B_i \rvert]\}$.
Let $H$ be a graph isomorphic to $K_{s,s}$ with a bipartition $\{A,B\}$, where $A=\{a_i: i \in [s]\}$ and $B=\{b_i: i \in [s]\}$.
For every $i \in [s]$ and $j \in [s]$, define $f(a_ib_j)$ to be
	\begin{itemize}
		\item $0$, if $A_i$ is disjoint and non-adjacent in $G$ to $B_j$,
		\item $1$, if $A_i \cap B_j \neq \emptyset$,
		\item $(x,y)$, if $A_i \cap B_j = \emptyset$, $v_{i,x}$ is adjacent to $u_{j,y}$ in $G$, and subject to this, the lexicographic order of $(x,y)$ is minimal.
	\end{itemize}
So $f$ is an $(h^2+2)$-coloring of $E(H)$.
By Lemma \ref{bi_ramsey}, there exists a monochromatic $K_{m,m}$, say with color $c$.

If $c=1$, then some member of $\C_1$ intersects at least $m \geq h+1$ members of $\C_2$; but it is impossible since every member of $\C_1$ has size at most $h$ and members of $\C_2$ are pairwise disjoint, a contradiction.
If $c=(x,y)$ for some $x \in [h]$ and $y \in [h]$, then there exists a $K_{m,m}$-subgraph of $G$ consisting of the $x$-th vertex of each of some $m$ members of $\C_1$ and the $y$-th vertex of each of some $m$ members of $\C_2$, contradicting that $G$ has no $K_{r,r}$-subgraph.
So $c=0$.
Hence for each $i \in [2]$, there exists $\C_i' \subseteq \C_i$ with $\lvert \C_i' \rvert = m \geq w$ such that every member of $\C_1'$ is disjoint and non-adjacent in $G$ to every member of $\C_2'$.
\end{pf}

\begin{lemma} \label{parts_collection}
For any positive integers $r,h,w,p$, there exists a positive integer $s=s(r,h,w,p)$ such that the following holds.
Let $G$ be a graph with no $K_{r,r}$-subgraph.
For every $i \in [p]$, let $\C_i$ be a collection of pairwise disjoint subsets of $V(G)$ with each member having size at most $h$.
If $\lvert \C_i \rvert \geq s$ for every $i \in [p]$, then for every $i \in [p]$, there exists a subset $\C_i'$ of $\C_i$ with $\lvert \C_i' \rvert=w$ such that for any $i_1,i_2$ with $1 \leq i_1 < i_2 \leq p$, every member of $\C_{i_1}'$ is disjoint and non-adjacent in $G$ to every member of $\C_{i_2}'$.
\end{lemma}

\begin{pf}
Let $r,h,w,p$ be positive integers.
Let $f(1,1)=w$. 
For any positive integers $x$ and $y$, let $f(x+1,1)=f(x,x)$ and let $f(x+1,y+1)=s_{\ref{2_part_collection}}(r,h,f(x+1,y))$.
Define $s=f(p,p)$.

We shall prove this lemma by induction on $p$.
Since $f(1,1)=w$, the case $p=1$ holds.
So we may assume that $p \geq 2$ and the lemma holds when $p$ is smaller.

Let $G$, $\C_1,\C_2, ... \C_p$ be as stated in this lemma.
Let $\C_{1,1}=\C_1$.
By induction, for every positive integer $x$ with $2 \leq x \leq p$, applying Lemma \ref{2_part_collection} with taking $(\C_1,\C_2)=(\C_{1,x-1},\C_x)$, we know that there exist a subset $\C_{1,x}$ of $\C_{1,x-1}$ of size $f(p,p-x+1)$ and a subset $\D_x$ of $\C_x$ with $\lvert \D_x \rvert = f(p,p-x+1) \geq f(p,1) = f(p-1,p-1)$ such that every member of $\C_{1,x}$ is disjoint and non-adjacent in $G$ to every member of $\D_x$.
So every member of $\C_{1,p}$ is disjoint and non-adjacent to every member of $\D_x$ for every $2 \leq x \leq p$.
By the induction hypothesis, for every $x$ with $2 \leq x \leq p$, since $\lvert \D_x \rvert \geq f(p-1,p-1)$, there exists $\C_x' \subseteq \D_x$ with $\lvert \C_x' \rvert = w$ such that for any $i_1,i_2$ with $2 \leq i_1 < i_2 \leq p$, every member of $\C_{i_1}'$ is disjoint and non-adjacent in $G$ to every member of $\C_{i_2}'$.
Hence the lemma follows by taking $\C_1'=\C_{1,p}$.
\end{pf}

\bigskip

Now we are ready to consider homomorphisms.

\begin{lemma} \label{shallow_density}
For any positive integers $r,t,w,h,a$ and positive real numbers $k,k'$, there exist positive real numbers $c=c(r,t,w,h,a,k,k')$, $r_0=r_0(r,h)$ and $r_1=r_1(r,h)$ such that if $G$ is a graph and $\epsilon,\epsilon'$ are nonnegative real numbers with $\epsilon \leq \frac{1}{16h+4}$ and $\epsilon' \leq r_1$ such that
	\begin{enumerate}
		\item for every induced subgraph $L$ of $G$, $\lvert E(L) \rvert \leq k\lvert V(L) \rvert^{1+\epsilon}$, 
		\item for every graph $L$ for which some subgraph of $G$ is isomorphic to a $[4h+1]$-subdivision of $L$, $\lvert E(L) \rvert \leq k'\lvert V(L) \rvert^{1+\epsilon'}$, and 
		\item there exists no $(1,2hr+1)$-model for a $2h$-shallow $K_{r,t}$-minor in $G$,
	\end{enumerate}
then for any quasi-order $Q=(U,\preceq)$ with $\lvert U \rvert \leq a$, legal $Q$-labelled graph $(H,f_H)$ on $h$ vertices, legal $Q$-labelling $f_G$ of $G$, and consistent set $\HH$ of homomorphisms from $(H,f_H)$ to the set of all induced subgraphs of $(G,f_G)$, there exists a nonnegative integer $q$ such that the following hold. 
	\begin{enumerate}
		\item The number of members of $\HH$ with codomain $V(G)$ is at most $c\lvert V(G) \rvert^{q+(1+\epsilon-\frac{1}{\epsilon_0})\lvert V(H) \rvert^2}$, where $\epsilon_0=(1+\epsilon'(r-1))^{r_0}$.  
		\item There exists an independent collection $\L$ of separations of $(H,f_H)$ with $\lvert \L \rvert = q$ such that for every $(X,Y) \in \L$,
			\begin{enumerate}
				\item every vertex in $X \cap Y$ is adjacent in $H$ to a vertex in $X-Y$, and
				\item if every member of $\HH$ is injective, then $H[X-Y]$ is connected.
			\end{enumerate}
		\item There exists an induced subgraph $(H',f_{H'})$ of $(G,f_G)$ such that 
			\begin{enumerate}
				\item there exists an onto homomorphism $\phi_H \in \HH$ from $(H,f_H)$ to $(H',f_{H'})$ such that 
					\begin{itemize}
						\item $\phi_H(A-B) \cap \phi_H(B)=\emptyset$ for every $(A,B) \in \L$,
						\item $\{(N_{H'}[\phi_H(A-B)], V(H')-\phi_H(A-B)): (A,B) \in \L\}$, denoted by $\L_{H'}$, is an independent collection of separations of $(H',f_{H'})$, and
						\item $\lvert \L_{H'} \rvert  = \lvert \L \rvert$,
					\end{itemize}
				\item $(H',f_{H'}) \wedge_w \L_{H'}$ is isomorphic to an induced subgraph of $(G,f_G)$. 
			\end{enumerate}
		\item There exists a spanning legal $Q$-labelled subgraph $(H'',f_{H''})$ of $(H',f_{H'})$ such that
			\begin{enumerate}
				\item $\phi_H$ is an onto homomorphism from $(H,f_H)$ to $(H'',f_{H''})$ such that $\{(N_{H''}[\phi_H(A-B)], V(H'')-\phi_H(A-B)): (A,B) \in \L\}$, denoted by $\L_{H''}$, is an independent collection of separations of $(H'',f_{H''})$ of order at most $r-1$ with $\lvert \L_{H''} \rvert  = \lvert \L \rvert$, and
				\item $(H'',f|_{H''}) \wedge_w \L_{H''}$ is isomorphic to a legal $Q$-labelled subgraph of $(G,f_G)$. 
			\end{enumerate}
	\end{enumerate}
\end{lemma}

\begin{pf}
Let $r,t,w,h,a$ be positive integers, and let $k,k'$ be positive real numbers. 
We define the following.
	\begin{itemize}
		\item Define $r_0=r_{\ref{island_4}}(r,h)$, where $r_{\ref{island_4}}$ is the real number $r_0$ mentioned in Lemma \ref{island_4}.
		\item Define $r_1=r'_{\ref{island_4}}(r,h)$, where $r'_{\ref{island_4}}$ is the real number $r_1$ mentioned in Lemma \ref{island_4}.
		\item Let $w_1=\max_{1 \leq i \leq h}s_{\ref{parts_collection}}(r+t,h,w,i)$, where $s_{\ref{parts_collection}}$ is the integer $s$ mentioned in Lemma \ref{parts_collection}.
		\item Let $\ell^*$ be the number of isomorphism classes of $Q$-labelled graphs with vertex-set contained in $[h]$, among all isomorphism classes of quasi-orders $Q$ with ground set size at most $a$.
		\item Let $t'=rh+w_1\ell^*$. 
		\item Let $N_1=\lceil N_{\ref{island_4}}(r,t,t',a,t',h,h,k,k') \rceil$, where $N_{\ref{island_4}}$ is the real number $N$ mentioned in Lemma \ref{island_4}.
		\item Let $b_1=b_{\ref{island_4}}(r,t,t',a,t',h,h,k,k')$, where $b_{\ref{island_4}}$ is the real number $b$ mentioned in Lemma \ref{island_4}.
		\item Define $c= c_{\ref{container}}(r-1,h,N_1) \cdot b_1^h$, where $c_{\ref{container}}$ is the number $c$ mentioned in Lemma \ref{container}. 
	\end{itemize}

Let $G$, $\epsilon,\epsilon'$, $Q=(U,\preceq)$, $(H,f_H)$, $f_G$, and $\HH$ be as stated in the lemma.
Let $\epsilon_0=(1+\epsilon'(r-1))^{r_0}$.
Note that $\epsilon_0 \geq 1$, so $1+\epsilon-\frac{1}{\epsilon_0} \geq 0$.

Since $h \geq 1$, $\max\{h+1,2h\}=2h$.
By Lemma \ref{island_4}, there exist an ordering $\sigma$ of $G$ and a sequence $(S_i: i \in [\lvert V(G) \rvert])$ such that for every $i \in [\lvert  V(G) \rvert]$,
	\begin{itemize}
		\item[(i)] $S_i$ is a subset of $\{v \in V(G): \sigma(v) \geq i+1\}$ with size at most $r-1$, and
		\item[(ii)] if $i \in [\lvert V(G) \rvert-N_1]$, then the $p$-basin $B_i$ with respect to $\sigma$ and $(S_j: j \in [\lvert V(G) \rvert])$ at $i$, where $p=\max\{h-1,0\}=h-1$, satisfies the following:
			\begin{itemize}
				\item[(iia)] $\lvert B_i \rvert \leq b_1 \lvert V(G) \rvert^{(1+\epsilon-\frac{1}{\epsilon_0})h}$. 
				\item[(iib)] There exists a subset $Z^*_i=\{z_{i,1},z_{i,2},...,z_{i,\lvert Z^*_i \rvert}\}$ of $\{v \in V(G): \sigma(v) \geq i+1\}$ with $\lvert Z^*_i \rvert \geq t'$.
				\item[(iic)] There exist $\lvert Z^*_i \rvert$ pairwise disjoint subsets $Z'_{i,1},Z'_{i,2},..., \allowbreak Z'_{i,\lvert Z^*_i \rvert}$ of $\{v \in V(G): \sigma(v) \geq i+1\}-S_i$ disjoint from $N_{G_{\sigma, \geq i}}[B_i]$.
				\item[(iid)] For any set $\HH'$ of homomorphisms from an element in $\{(L,f_L): (L,f_L)$ is a legal $Q$-labelled graph with $V(L) \subseteq [h]\}$ to the set of all induced subgraphs of $(G_{\sigma, \geq i},f_G|_{V(G_{\sigma, \geq i})})$ satisfying (CON1) and (CON2), if $(L,f_L)$ is a legal $Q$-labelled graph such that there exists $\phi \in \HH'$ from $(L,f_L)$ to $(G_{\sigma, \geq i}[B_i \cup S_i],f_G|_{B_i \cup S_i})$ with $\sigma^{-1}(i) \in \phi(V(L))$, then 
					\begin{itemize}
						\item for every $j \in [\lvert Z^*_i \rvert]$, there exists $\phi_{i,j} \in \HH'$ from $(L,f_L)$ to $(G_{\sigma,\geq i}[Z'_{i,j} \cup S_i],f_G|_{Z'_{i,j} \cup S_i})$ such that
							\begin{itemize}
								\item $z_{i,j} \in \phi_{i,j}(V(L))$, 
								\item for every $y \in S_i$ and $u \in V(L)$, $y=\phi(u)$ if and only if $y=\phi_{i,j}(u)$, 
								\item for every $u \in V(L)$, $\phi(u)=\sigma^{-1}(i)$ if and only if $\phi_{i,j}(u)=z_{i,j}$, and 
								\item there exists an isomorphism $\iota_{i,j}$ from $(G_{\sigma, \geq i}[\phi(V(L)) \cup S_i],f|_{\phi(V(L)) \cup S_i})$ to \linebreak $(G_{\sigma, \geq i}[\phi_{i,j}(V(L)) \cup S_i], \allowbreak f|_{\phi_{i,j}(V(L)) \cup S_i})$ such that $\iota_{i,j}(\phi(u))=\phi_{i,j}(u)$ for every $u \in V(L)$, and $\iota_{i,j}(y)=y$ for every $y \in S_i$.
							\end{itemize}
						\item for any $I \subseteq [i-1]$ with $\lvert I \rvert \leq h$, there exists $M_{i,I} \subseteq Z_i^*$ with $\lvert M_{i,I} \rvert \geq t'$ such that for every $v \in \sigma^{-1}(I)$ and $z_{i,j} \in M_{i,I}$, $\{u \in V(L): v\phi(u) \in E(G)\} \supseteq \{u \in V(L): v\phi_{i,j}(u) \in E(G)\}$. 
					\end{itemize}
				\item[(iie)] There exist $\lvert Z^*_i \rvert$ sets $Z_{i,1},Z_{i,2},..., Z_{i,\lvert Z^*_i \rvert}$ such that for every $j \in [\lvert Z^*_i \rvert]$, $N_{G_{\sigma, \geq i}-S_i}^{\leq h-1}[z_{i,j}] \allowbreak \subseteq Z_{i,j} \subseteq Z'_{i,j}$ and $N_{G_{\sigma, \geq i}}[Z_{i,j}] \subseteq Z_{i,j}' \cup S_i$.
			\end{itemize}
	\end{itemize}
So by (i), (iia) and Lemma \ref{container}, there exists an integer $q$ such that 
	\begin{itemize}
		\item[(iii)] The number of members of $\HH$ with codomain $V(G)$ is at most $c_{\ref{container}}(r-1,h,N_1) \cdot (b_1\lvert V(G) \rvert^{(1+\epsilon-\frac{1}{\epsilon_0})h})^h \cdot \lvert V(G) \rvert^{q} \leq c\lvert V(G) \rvert^{q+(1+\epsilon-\frac{1}{\epsilon_0})h^2}$, and 
		\item[(iv)] there exist a homomorphism $\phi_0 \in \HH$ from $(H,f_H)$ to $(G,f_G)$, an independent collection $\L$ of separations of $H$ with $\lvert \L \rvert =q$, and an injection $\iota_{\phi_0}: \L \rightarrow [\lvert V(G) \rvert-N_1]$ such that for every $(X,Y) \in \L$, there exists $i_X \in [\lvert V(G) \rvert-N_1]$ with $i_X=\iota_{\phi_0}((X,Y))$ such that 
			\begin{itemize}
				\item[(iva)] $\sigma^{-1}(i_X) \in \phi_0(X-Y) \subseteq B_{i_X}$,
				\item[(ivb)] $\phi_0(X \cap Y) \subseteq S_{i_X}$, 
				\item[(ivc)] for every component $C$ of $H[X-Y]$, $\sigma^{-1}(i_X) \in \phi_0(V(C))$, and
				\item[(ivd)] every vertex in $X \cap Y$ is adjacent in $H$ to some vertex in $X-Y$.
			\end{itemize}
	\end{itemize}
Note that (iv) implies that $q \geq 0$.
So Conclusion 1 of this lemma follows from (iii).
And Conclusion 2(a) follows from (ivd).
Moreover, if every member of $\HH$ is injective, then $\phi_0$ is injective, so $H[X-Y]$ is connected by (ivc).
So Conclusion 2 of this lemma holds.

For every $i \in [\lvert V(G) \rvert-N_1]$, by permuting indices, we may assume that $M_{i,\phi_0(V(H)) \cap \sigma^{-1}([i-1])} \allowbreak \supseteq \{z_{i,j}: j \in [t']\}$.

Since $\lvert V(H) \rvert=h$, $\lvert \phi_0(V(H)) \rvert \leq h$, so $\lvert \phi_0(V(H)) \cup \bigcup_{v \in \phi_0(V(H))}S_{\sigma(v)} \rvert \leq rh$ by (i).
For each $i \in [\lvert V(G) \rvert-N_1]$, $Z'_{i,j}$'s are pairwise disjoint by (iic), so by permuting indices, we may assume that $Z'_{i,j} \cap (\phi_0(V(H)) \cup \bigcup_{v \in \phi_0(V(H))}S_{\sigma(v)})=\emptyset$ for every $j \in [t'-rh]$.

For every $(X,Y) \in \L$, let $(L_X,f_{L_X})$ be a legal $Q$-labelled graph with $V(L_X) \subseteq [h]$ isomorphic to $(H[X \cup \phi_0^{-1}(S_{i_X} \cap \phi_0(V(H)))],f_H|_{X \cup \phi_0^{-1}(S_{i_X} \cap \phi_0(V(H)))})$ with an isomorphism $\iota_X$, and let $\HH'_X = \{\phi|_{X \cup \phi_0^{-1}(S_{i_X} \cap \phi_0(V(H)))} \circ \iota_X: \phi \in \HH,\phi(X \cup \phi_0^{-1}(S_{i_X} \cap \phi_0(V(H)))) \subseteq V(G_{\sigma, \geq i_X})\}$.
So $\HH'_X$ is a set of homomorphisms from an element in $\{(L,f_L): (L,f_L)$ is a legal $Q$-labelled graph with $V(L) \subseteq [h]\}$ to the set of all induced subgraphs of $(G_{\sigma, \geq i_X},f_G|_{V(G_{\sigma, \geq i_X})})$. 
By (iva) and (ivb), $\phi_0|_{X \cup \phi_0^{-1}(S_{i_X} \cap \phi_0(V(H)))} \circ \iota_X \in \HH'_X$.
Let $\HH''_X$ be the minimal set of homomorphisms from an element in $\{(L,f_L): (L,f_L)$ is a legal $Q$-labelled graph with $V(L) \subseteq [h]\}$ to the set of all induced subgraphs of $(G_{\sigma, \geq i_X},f_G|_{V(G_{\sigma, \geq i_X})})$ containing $\HH_X'$ and satisfying (CON1) and (CON2).

For every $(X,Y) \in \L$, since $\sigma^{-1}(i_X) \in \phi_0(X) \subseteq B_{i_X} \cup S_{i_X}$ by (iva) and (ivb), we know that by (iid) (with taking $i=i_X$, $\HH'=\HH''_X$, $(L,f_L)=(L_X,f_{L_X})$, $\phi=\phi_0|_{X \cup \phi_0^{-1}(S_{i_X} \cap \phi_0(V(H)))} \circ \iota_X$, and $I=\phi_0(V(H)) \cap \sigma^{-1}([i_X-1])$), there exists $J_X \subseteq [t'-rh]$ with $\lvert J_X \rvert \geq (t'-rh)/\ell^* = w_1$ such that
	\begin{itemize}
		\item[(v)] there exists a legal $Q$-labelled graph $(R_X,f_{R_X})$ such that for every $j \in J_X$, there exists a homomorphism $\phi_{i_X,j}$ from $(H[X \cup \phi_0^{-1}(S_{i_X} \cap \phi_0(V(H)))],f_H|_{X \cup \phi_0^{-1}(S_{i_X} \cap \phi_0(V(H)))})$ to $(G_{\sigma, \geq i_X}[Z'_{i_X,j} \cup S_{i_X}],f_G|_{Z'_{i_X,j} \cup S_{i_X}})$ and an isomorphism $\iota_{X,j}$ from $(R_X,f_{R_X})$ to $(G[\phi_{i_X,j}(X) \allowbreak \cup S_{i_X} \cup \{v \in \phi_0(V(H)): \sigma(v) \leq i_X-1\}],f|_{\phi_{i_X,j}(X) \cup S_{i_X} \cup \{v \in \phi_0(V(H)): \sigma(v) \leq i_X-1\}})$ such that  
			\begin{itemize}
				\item[(va)] $z_{i_X,j} \in \phi_{i_X,j}(X)$,
				\item[(vb)] for every $y \in S_{i_X}$ and $u \in X \cup \phi_0^{-1}(S_{i_X} \cap \phi_0(V(H)))$, $y = \phi_0(u)$ if and only if $y = \phi_{i_X,j}(u)$, 
				\item[(vc)] for every $u \in X$, $\phi_0(u)=\sigma^{-1}(i_X)$ if and only if $\phi_{i_X,j}(u)=z_{i_X,j}$, 
				\item[(vd)] there exists an isomorphism $\iota'_{i_X,j}$ from $(G[\phi_0(X) \cup S_{i_X}], f|_{\phi_0(X) \cup S_{i_X}})$ to $(G[\phi_{i_X,j}(X) \cup S_{i_X}], f|_{\phi_{i_X,j}(X) \cup S_{i_X}})$ such that $\iota'_{i_X,j} \circ \phi_0|_{X \cup \phi_0^{-1}(S_{i_X} \cap \phi_0(V(H)))} = \phi_{i_X,j}|_{X \cup \phi_0^{-1}(S_{i_X} \cap \phi_0(V(H)))}$ and $\iota'_{i_X,j}|_{S_{i_X}}$ is the identity function, 
				\item[(ve)] for every $y \in S_{i_X} \cup (\phi_0(V(H)) \cap \sigma^{-1}([i_X-1]))$, $\iota_{X,j}^{-1}(y)=\iota_{X,j'}^{-1}(y)$ for any $j' \in J_X$,  
				\item[(vf)] for every $x \in X \cup \phi_0^{-1}(S_{i_X} \cap \phi_0(V(H)))$, $\iota_{X,j}^{-1}(\phi_{i_X,j}(x))=\iota_{X,j'}^{-1}(\phi_{i_X,j'}(x))$ for any $j' \in J_X$, and
				\item[(vg)]for every $x \in X \cup \phi_0^{-1}(S_{i_X} \cap \phi_0(V(H)))$ and $y \in \phi_0(V(H)) \cap \sigma^{-1}([i_X-1])$, if $y\phi_{i_X,j}(x) \in E(G)$, then $y\phi_0(x) \in E(G)$.
			\end{itemize}
	\end{itemize}

\medskip

\noindent{\bf Claim 1:} For every $(X,Y) \in \L$ and $j \in J_X$, the following hold.
	\begin{itemize}
		\item For every $u \in X \cap Y$, $\phi_0(u) = \phi_{i_X,j}(u)$.
		\item $\phi_{i_X,j}(X) \cap S_{i_X} = \phi_0(X \cap Y) = \phi_{i_X,j}(X \cap Y)$.
		\item $\phi_{i_X,j}(X-Y) \subseteq Z_{i_X,j}$.
	\end{itemize}

\noindent{\bf Proof of Claim 1:}
Let $(X,Y) \in \L$.
Let $j \in J_X$.
Let $u \in X \cap Y$.
Then $\phi_0(u) \in S_{i_X}$ by (ivb).
Since $u \in X \cap Y \subseteq X$, $\phi_0(u)=\phi_{i_X,j}(u)$ by (vb).
This proves the first statement of this claim.
Note that this implies that $\phi_0(X \cap Y) = \phi_{i_X,j}(X \cap Y)$.

If $v \in \phi_{i_X,j}(X) \cap S_{i_X}$, then $v \in \phi_0(X) \cap S_{i_X}$ by (vb), so $v \in \phi_0(X \cap Y)$ by (iva).
Hence $\phi_{i_X,j}(X) \cap S_{i_X} \subseteq \phi_0(X \cap Y) = \phi_{i_X,j}(X \cap Y)$.
Since $\phi_0(X \cap Y) \subseteq S_{i_X}$ by (ivb), $\phi_{i_X,j}(X) \cap S_{i_X} \subseteq \phi_0(X \cap Y) = \phi_0(X \cap Y) \cap S_{i_X} = \phi_{i_X,j}(X \cap Y) \cap S_{i_X} \subseteq \phi_{i_X,j}(X) \cap S_{i_X}$.
So all inclusions are equalities.
This proves the second statement of this claim.

Since $\phi_{i_X,j}(X) \subseteq S_{i_X} \cup Z'_{i_X,j}$ by (v), and $\phi_{i_X,j}(X) \cap S_{i_X}= \phi_{i_X,j}(X \cap Y)$ by Statement 2 of this claim, we know $\phi_{i_X,j}(X-Y) \subseteq Z'_{i_X,j}$.
By (ivc), for every component $C$ of $H[X-Y]$, there exists $v_C \in V(C)$ such that $\sigma^{-1}(i_X)=\phi_0(v_C)$, so by (vc), $\phi_{i_X,j}(v_C)=z_{i_X,j}$.
For any component $C$ of $H[X-Y]$ and vertex $v$ of $C$, since there exists a path $P_C$ in $C \subseteq H[X-Y]$ from $v_C$ to $v$, there exists a path in $G_{\sigma, \geq i_X}$ from $\phi_{i_X,j}(v_C)=z_{i_X,j}$ to $\phi_{i_X,j}(v)$ whose vertex-set is contained in $\phi_{i_X,j}(V(P_C)) \subseteq \phi_{i_X,j}(X-Y) \subseteq Z'_{i_X,j}$, so $\phi_{i_X,j}(v) \in N_{G_{\sigma, \geq i_X}[Z'_{i_X,j}]}^{\leq \lvert E(P_C) \rvert}[z_{i_X,j}] \subseteq N_{G_{\sigma, \geq i_X}-S_{i_X}}^{\leq h-1}[z_{i_X,j}] \subseteq Z_{i_X,j}$ by (iic) and (iie).
Hence $\phi_{i_X,j}(V(C)) \subseteq Z_{i_X,j}$ for every component $C$ of $H[X-Y]$.
Therefore, $\phi_{i_X,j}(X-Y) \subseteq Z_{i_X,j}$.
This proves the third statement of this claim.
$\Box$

\medskip

Denote the elements of $\L$ by $(X_1,Y_1),(X_2,Y_2),...,(X_{\lvert \L \rvert},Y_{\lvert \L \rvert})$.
For every $j \in [\lvert \L \rvert]$, let $i_j=i_{X_j}$.

By Statement 3 of Claim 1, (iic) and (iie), we know that for every $j \in [\lvert \L \rvert]$, $\{\phi_{i_j,\ell}(X_j-Y_j): \ell \in J_{X_j}\}$ is a collection of pairwise disjoint subsets of $V(G)$, and each member has size at most $h$.
Since there exists no $(1,2hr+1)$-model for a $2h$-shallow $K_{r,t}$-minor in $G$, $G$ has no $K_{r+t,r+t}$-subgraph.
So by Lemma \ref{parts_collection} (with taking $\C_j=\{\phi_{i_j,\ell}(X_j-Y_j): \ell \in J_{X_j}\}$ for every $j \in [\lvert \L \rvert]$), we know that for every $(X,Y) \in \L$, there exists $J_X' \subseteq J_X$ with $\lvert J_X' \rvert =w$ such that for any distinct $j_1,j_2 \in [\lvert \L \rvert]$, every member of $\{\phi_{i_{j_1},\ell}(X_{j_1}-Y_{j_1}): \ell \in J_{X_{j_1}}'\}$ is disjoint and non-adjacent in $G$ to every member of $\{\phi_{i_{j_2},\ell}(X_{j_2}-Y_{j_2}): \ell \in J'_{X_{j_2}}\}$.
And by Statement 3 in Claim 1, (iic) and (iie), for every $j \in [\lvert \L \rvert]$, members of $\{\phi_{i_j,\ell}(X_j-Y_j): \ell \in J_{X_j}'\}$ are pairwise disjoint and non-adjacent in $G$.
Since $J_X' \subseteq J_X \subseteq [t'-rh]$, by permuting indices, we may assume that $J_X'=[w]$ for every $(X,Y) \in \L$.
Hence we have
	\begin{itemize}
		\item[(vi)] $\{\phi_{i_X,\ell}(X-Y): (X,Y) \in \L, \ell \in [w]\}$ is a collection of pairwise disjoint and pairwise non-adjacent in $G$ subsets of $V(G)$, and this collection has size $\lvert \L \rvert w$. 
	\end{itemize}

Define $G^*=G[\phi_0(\bigcap_{(A,B) \in \L}B) \cup \bigcup_{(X,Y) \in \L}\bigcup_{\ell \in [w]}\phi_{i_X,\ell}(X-Y)]$.

Since $Z_{i_X,j} \cap \phi_0(V(H))=\emptyset$ for all $(X,Y) \in \L$ and $j \in [w]=J_X' \subseteq [t'-rh]$, by Statement 3 in Claim 1 and (vi), we have
	\begin{itemize}
		\item[(vii)] $\{\phi_{i_j,\ell}(X_j-Y_j), \phi_0(\bigcap_{(A,B) \in \L}B): j \in [\lvert \L \rvert], \ell \in [w]\}$ is a partition of $V(G^*)$ with size $\lvert \L \rvert w+1$, where $\phi_0(\bigcap_{(A,B) \in \L}B)$ is the only possibly empty member and is the only possible member intersecting $\phi_0(V(H))$. 
	\end{itemize}

Define $H' = G[\phi_0(\bigcap_{(A,B) \in \L}B) \cup \bigcup_{(X,Y) \in \L}\phi_{i_X,1}(X-Y)]$.
Define $\phi_H: V(H) \rightarrow V(H')$ such that for every $x \in \bigcap_{(A,B) \in \L}B$, $\phi_H(x)=\phi_0(x)$, and for every $x \in A-B$ for some $(A,B) \in \L$, $\phi_H(x)=\phi_{i_A,1}(x)$.
Note that $\phi_H$ is well-defined since $\L$ is an independent collection of separations of $H$.
By Statement 1 of Claim 1, we have 
	\begin{itemize}
		\item[(viii)] for any $(X,Y) \in \L$ and $v \in X$, $\phi_H(v)=\phi_{i_X,1}(v)$.
	\end{itemize}

Define $H''$ to be the spanning subgraph of $H'$ obtained from $H'$ by, for every $(X,Y) \in \L$, deleting all edges of $H'$ with one end in $\phi_{i_X,1}(X-Y)$ and the other end in $V(H') \cap \sigma^{-1}[i_X-1]$. 
Let $f_{H''}$ be the function whose domain consisting of the marches $\pi$ in the domain of $f_G|_{V(H')}$ such that the entries of $\pi$ form a clique in $H''$, and for each $x$ in the domain of $f_{H''}$, $f_{H''}(x)=f_G(x)$.
So $f_{H''}$ is a legal labelling of $H''$, and $(H'',f_{H''})$ is a spanning subgraph of $(H',f_G|_{V(H')})$ and hence a subgraph of $(G,f_G)$.

\medskip

\noindent{\bf Claim 2:} $\phi_H$ is an onto homomorphism from $(H,f_H)$ to $(H'',f_{H''})$.

\noindent{\bf Proof of Claim 2:}
Note that $V(H')=V(H'')$, so $\phi_H$ is a function from $V(H)$ to $V(H'')$.
Clearly, $\phi_H$ is onto.
So it suffices to show that $\phi_H$ is a homomorphism from $(H,f_H)$ to $(H'',f_{H''})$.

Suppose to the contrary that there exists $uv \in E(H)$ such that $\phi_H(u)\phi_H(v) \not \in E(H'')$.
If both $u,v$ are in $\bigcap_{(A,B) \in \L}B$, then $\phi_H(u)\phi_H(v)=\phi_0(u)\phi_0(v) \in E(G)$ since $\phi_0$ is a homomorphism; since both $\phi_0(u)$ and $\phi_0(v)$ are not in $\bigcup_{j \in [\lvert \L \rvert]}\bigcup_{\ell \in [w]}\phi_{i_j,\ell}(X_j-Y_j)$ by (vii), $\phi_0(u)\phi_0(v) \in E(H'')$, a contradiction.
So there exists $(X,Y) \in \L$ such that both $u,v$ are in $X$.
Hence by (viii), $\phi_H(u)\phi_H(v)=\phi_{i_X,1}(u)\phi_{i_X,1}(v) \in E(G)$, since $\phi_{i_X,1}$ is a homomorphism.
By (v) and (vii), $\phi_{i_X,1}(u)\phi_{i_X,1}(v) \in E(H'')$, a contradiction.

So $\phi_H$ preserves adjacency.
Let $\pi$ be a march in the domain of $f_H$.
Since $f_H$ is legal, either all entries of $\pi$ are in $X$ for some $(X,Y) \in \L$, or all entries of $\pi$ are in $\bigcap_{(A,B) \in \L}B$.
For the former, $f_H(\pi) \preceq f_G(\phi_{i_X,1}(\pi))=f_G(\phi_H(\pi))$ by (viii); for the latter, $f_H(\pi) \preceq f_G(\phi_0(\pi)) = f_G(\phi_H(\pi))$.
Since $f_H$ is legal and $\phi_H$ preserves adjacency, $\phi_H(\pi)$ is a clique in $H''$.
So $\phi_H(\pi)$ is in the domain of $f_{H''}$, and $f_{H''}(\phi_H(\pi))=f_G(\phi_H(\pi))$.
Hence $f_H(\pi) \preceq f_{H''}(\phi_H(\pi))$.
Therefore, $\phi_H$ is an onto homomorphism from $(H,f_H)$ to $(H'',f_{H''})$.
$\Box$

\medskip

\noindent{\bf Claim 3:} $\phi_H$ is an onto homomorphism from $(H,f_H)$ to $(H',f_G|_{V(H')})$ and $\phi_H \in \HH$.

\noindent{\bf Proof of Claim 3:}
Since $(H'',f_{H''})$ is a spanning subgraph of $(H',f_G|_{V(H')})$, $\phi_H$ is an onto homomorphism from $(H,f_H)$ to $(H',f_G|_{V(H')})$ by Claim 2.
Let $G_0$ be the spanning subgraph of $G[\phi_0(V(H))]$ obtained from $G[\phi_0(V(H))]$ by, for every $(X,Y) \in \L$, deleting all edges of $G[\phi_0(V(H))]$ satisfying that
	\begin{itemize}
		\item one end is in $\phi_0(X-Y)$ and the other end is in $N_{G_{\sigma, \geq i_X}}[B_{i_X}]-(S_{i_X} \cup \phi_0(X))$, or
		\item one end is in $\phi_0(u)$ for some $u \in X-Y$ and the other end $v$ is in $\phi_0(V(H)) \cap \sigma^{-1}[i_X-1]$ with $v\phi_{i_X,1}(u) \not \in E(G)$.
	\end{itemize}
Let $f_0$ be the function whose domain consists of the marches $\pi$ in the domain of $f_G|_{\phi_0(V(H))}$ such that the entries of $\pi$ form a clique in $G_0$, and for each $x$ in the domain of $f_0$, $f_0(x)=f_G(x)$.
Note that $(G_0,f_0)$ is a spanning subgraph of $(G[\phi_0(V(H))],f_G|_{\phi_0(V(H))})$, and $\phi_0$ is also a homomorphism from $(H,f_H)$ to $(G_0,f_0)$ by (iva) and (ivb).
Since $\HH$ satisfies (CON1) and $(G[\phi_0(V(H))],f_G|_{\phi_0(V(H))})$ is an induced subgraph of $(G,f_G)$, we know that the function $\phi_0'$ obtained from $\phi_0$ by restricting the codomain of $\phi_0$ to $\phi_0(V(H))$ is a member in $\HH$ and an onto homomorphism from $(H,f_H)$ to $(G_0,f_0)$. 

By (iic), (iie), (vd), (vg) and (vi), $(H',f_G|_{V(H')})$ is isomorphic to $(G_0,f_0)$ with an isomorphism $\iota_H^*$ such that $(\iota^*_H)^{-1} \circ \phi'_0=\phi_H$. 
Since $\HH$ satisfies (CON3), $\phi_H=(\iota^*_H)^{-1} \circ \phi'_0 \in \HH$.
$\Box$

\medskip

By (vii) and (viii), $\phi_H(A-B) \cap \phi_H(B)=\emptyset$ for every $(A,B) \in \L$.
Define $\L_{H'} = \{(N_{H'}[\phi_H(X-Y)], V(H')-\phi_H(X-Y)): (X,Y) \in \L\}$.
By (vi) and (viii), $\L_{H'}$ is an independent collection of separations of $H'$, and $\lvert \L_{H'} \rvert = \lvert \L \rvert$.
These together with Claim 3 imply Conclusion 3(a).

Recall that $G^*=G[\phi_0(\bigcap_{(A,B) \in \L}B) \cup \bigcup_{(X,Y) \in \L}\bigcup_{\ell \in [w]}\phi_{i_X,\ell}(X-Y)]$.
Note that if there exists $v \in N_{G^*}[\phi_{i_X,\ell}(X-Y)]-\phi_{i_X,\ell}(X)$ for some $(X,Y) \in \L$ and $\ell \in [w]$, then $v \in \phi_0(\bigcap_{(A,B) \in \L}B)$, and either $\sigma^{-1}(v) \leq i_X-1$ or $v \in S_{i_X}$ by Statement 3 of Claim 1 and (iie), so $v \in N_{G^*}[\phi_{i_X,\ell'}(X-Y)]-\phi_{i_X,\ell'}(X)$ for every $\ell' \in [w]$ by (vd), (ve) and (vf). 
Hence by (vb), (ve) and (vf), $(H',f_G|_{V(H')}) \wedge_w \L_{H'}$ is isomorphic to $(G^*,f_G|_{V(G^*)})$.
Therefore Conclusion 3(b) of this lemma holds.

Let $\L_{H''} = \{(N_{H''}[\phi_H(X-Y)], V(H'')-\phi_H(X-Y)): (X,Y) \in \L\}$.
Since $(H'',f_{H''})$ is a spanning legal $Q$-labelled subgraph of $(H',f_G|_{V(H')})$, Conclusion 4(b) follows from Conclusion 3(b).
By Claim 3 and Conclusion 3(a) of this lemma, to prove this lemma, it suffices to prove that every separation in $\L_{H''}$ has order at most $r-1$.

Suppose to the contrary that there exists $(X,Y) \in \L$ such that the order of $(N_{H''}[\phi_H(X-Y)], V(H'')-\phi_H(X-Y))$ is at least $r$.
Note that $(N_{H''}[\phi_H(X-Y)], V(H'')-\phi_H(X-Y)) = (N_{H''}[\phi_{i_X,1}(X-Y)], V(H'')-\phi_{i_X,1}(X-Y))$ by (viii).
And the order of $(N_{H''}[\phi_{i_X,1}(X-Y)], V(H'')-\phi_{i_X,1}(X-Y))$ is $\lvert N_{H''}(\phi_{i_X,1}(X-Y)) \rvert$, so $N_{H''}(\phi_{i_X,1}(X-Y)) \not \subseteq S_{i_X}$ by (i).
Hence there exists $u \in \phi_{i_X,1}(X-Y)$ and $v \in N_{H''}(\phi_{i_X,1}(X-Y))-S_{i_X}$.
By Statement 3 in Claim 1 and (iie), either $\sigma(v) \leq i_X-1$, or $\sigma(v) \geq i_X$ and $v \in Z'_{i_X,1}$.
Since $u \in \phi_{i_X,1}(X-Y)$ and $uv \in E(H'')$, $\sigma(v) \geq i_X$ by the definition of $H''$.
So $v \in Z'_{i_X,1}$.
Since $v \in V(H'')=V(H')$ and $uv \in E(H'') \subseteq E(H')$, by (vi) and (vii), $v \in \phi_0(\bigcap_{(A,B) \in \L}B) \subseteq \phi_0(V(H))$.
But $Z'_{i_X,1} \cap \phi_0(V(H)) = \emptyset$, a contradiction.

Therefore, every separation in $\L_{H''}$ has order at most $r-1$.
This proves Conclusion 4(a) of this lemma and completes the proof.
\end{pf}

\begin{lemma} \label{shallow_density_dup}
For any quasi-order $Q$ with finite ground set, positive integers $r,t,h$, hereditary class $\G$ of legal $Q$-labelled graphs, legal $Q$-labelled graph $(H,f_H)$ on $h$ vertices, consistent set $\HH$ of homomorphisms from $(H,f_H)$ to $\G$, and positive real numbers $k,k'$, there exist positive real numbers $c=c(Q,r,t,h,\G,(H,f_H),\HH,k,k')$, $r_0=r_0(r,h)$ and $r_1=r_1(r,h)$ such that if $(G,f_G)$ is a member in $\G$, and $\epsilon,\epsilon'$ are nonnegative real numbers with $\epsilon \leq \frac{1}{16h+4}$ and $\epsilon' \leq r_1$ such that
	\begin{enumerate}
		\item for every induced subgraph $L$ of $G$, $\lvert E(L) \rvert \leq k\lvert V(L) \rvert^{1+\epsilon}$, 
		\item for every graph $L$ in which some subgraph of $G$ is isomorphic to a $[4h+1]$-subdivision of $L$, $\lvert E(L) \rvert \leq k'\lvert V(L) \rvert^{1+\epsilon'}$, and
		\item there exists no $(1,2hr+1)$-model for a $2h$-shallow $K_{r,t}$-minor in $G$,
	\end{enumerate}
then there exists a nonnegative integer $q$ such that 
	\begin{enumerate}
		\item there are at most $c\lvert V(G) \rvert^{q+(1+\epsilon-\frac{1}{\epsilon_0})\lvert V(H) \rvert^2}$ homomorphisms in $\HH$ from $(H,f_H)$ to $(G,f_G)$, where $\epsilon_0=(1+\epsilon'(r-1))^{r_0}$, 
		\item there exists a $(\G,\HH)$-duplicable independent collection $\L$ of separations of $(H,f_G)$ with $\lvert \L \rvert = q$, and 
		\item if every member of $\HH$ is injective, then for every separation $(X,Y)$ in $\L$, it has order at most $r-1$ and $H[X-Y]$ is connected. 
	\end{enumerate}
\end{lemma}

\begin{pf}
Let $Q$ be a quasi-order with finite ground set.
Let $r,t,h$ be positive integers.
Let $\G$ be a hereditary class of legal $Q$-labelled graphs.
Let $(H,f_H)$ be a legal $Q$-labelled graph on $h$ vertices.
Let $\HH$ be a consistent set of homomorphisms from $(H,f_H)$ to $\G$. 
Let $k,k'$ be positive real numbers.

For any non-$(\G,\HH)$-duplicable independent collection $\L$ of separations of $H$, let $s_{\L}$ be the minimum positive integer $s$ such that for every $Q$-labelled graph $(H',f_{H'})$ and $\phi \in \HH$ from $(H,f_H)$ to $(H',f_{H'})$ satisfying that $\phi$ is onto, $\phi(A-B) \cap \phi(B)=\emptyset$ for every $(A,B) \in \L$, $\W_{\L,\phi}$ is an independent collection of separations of $(H',f_{H'})$ with size $\lvert \L \rvert$, we have $(H',f_{H'}) \wedge_{s'} \W_{\L,\phi} \not \in \G$ for every $s' \geq s$, where $\W_{\L,\phi}=\{(N_{H'}[\phi(A-B)], V(H')-\phi(A-B)): (A,B) \in \L\}$. 
Note that such $s_\L$ exists since there are only finitely many legal $Q$-labelled graphs $(H',f_{H'})$ such that there exists an onto function $\phi \in \HH$ from $(H,f_H)$ to $(H',f_{H'})$.
If there exists no non-$(\G,\HH)$-duplicable independent collections $\L$ of separations of $H$, then let $s=1$; otherwise, let $s=\max_{\L}s_\L$, where the maximum is over all non-$(\G,\HH)$-duplicable independent collections $\L$ of separations of $H$. 
Note that this maximum exists since there are only finitely many independent collections of separations of $H$.
Let $a = \lvert U \rvert$, where $U$ is the ground set of $Q$.
Define $c=c_{\ref{shallow_density}}(r,t,s+t,h,a,k,k')$, where $c_{\ref{shallow_density}}$ is the real number $c$ mentioned in Lemma \ref{shallow_density}.
Define $r_0=r_{\ref{shallow_density}}(r,h)$ and $r_1=r'_{\ref{shallow_density}}(r,h)$, where $r_{\ref{shallow_density}}$ and $r'_{\ref{shallow_density}}$ are the real numbers $r_0$ and $r_1$ in Lemma \ref{shallow_density}.

Let $(G,f_G) \in \G$ and let $\epsilon$ and $\epsilon'$ be the real numbers as stated in the lemma.
Let $\HH'$ be the subset of $\HH$ consisting of the members of $\HH$ from $(H,f_H)$ to some induced subgraphs of $(G,f_G)$.
Since $\HH$ is consistent, $\HH'$ is consistent.
Applying Lemma \ref{shallow_density} (with taking $\HH=\HH'$), there exists a nonnegative integer $q$ such that the following hold.
	\begin{itemize}
		\item[(i)] There are at most $c\lvert V(G) \rvert^{q+(1+\epsilon-\frac{1}{\epsilon_0})\lvert V(H) \rvert^2}$ homomorphisms in $\HH'$ from $(H,f_H)$ to $(G,f_G)$, where $\epsilon_0=(1+\epsilon'(r-1))^{r_0}$.
		\item[(ii)] There exist an independent collection $\L$ of separations of $(H,f_H)$ with $\lvert \L \rvert = q$ such that for every $(X,Y) \in \L$,
			\begin{itemize}
				\item[(iia)] every vertex in $X \cap Y$ is adjacent in $H$ to a vertex in $X-Y$, and
				\item[(iib)] if every member of $\HH'$ is injective, then $H[X-Y]$ is connected.
			\end{itemize}
		\item[(iii)] There exists an induced subgraph $(H',f_{H'})$ of $(G,f_G)$ such that 
			\begin{itemize}
				\item[(iiia)] there exists an onto homomorphism $\phi_H \in \HH' \subseteq \HH$ from $(H,f_H)$ to $(H',f_{H'})$ such that 
					\begin{itemize}
						\item $\phi_H(A-B) \cap \phi_H(B)=\emptyset$ for every $(A,B) \in \L$,
						\item $\{(N_{H'}[\phi_H(A-B)], V(H')-\phi_H(A-B)): (A,B) \in \L\}$, denoted by $\L_{H'}$, is an independent collection of separations of $(H',f_{H'})$, and
						\item $\lvert \L_{H'} \rvert  = \lvert \L \rvert$,
					\end{itemize}
				\item[(iiib)] $(H',f|_{H'}) \wedge_{s+t} \L_{H'}$ is isomorphic to an induced subgraph of $(G,f_G)$.
			\end{itemize}
		\item[(iv)] There exists a legal spanning subgraph $(H'',f_{H''})$ of $(H',f_{H'})$ such that
			\begin{itemize}
				\item $\phi_H$ is an onto homomorphism from $(H,f_H)$ to $(H'',f_{H''})$ such that $\{(N_{H''}[\phi_H(A-B)], V(H'')-\phi_H(A-B)): (A,B) \in \L\}$, denoted by $\L_{H''}$, is an independent collection of separations of $(H'',f_{H''})$ of order at most $r-1$ with $\lvert \L_{H''} \rvert  = \lvert \L \rvert$, and
				\item $(H'',f|_{H''}) \wedge_{s+t} \L_{H''}$ is isomorphic to a legal subgraph $(G',f_{G'})$ of $(G,f_G)$. 
			\end{itemize}
	\end{itemize}

Since every member of $\HH$ with codomain $V(G)$ belongs to $\HH'$, Conclusion 1 follows from (i).

Suppose that $\L$ is non-$(\G,\HH)$-duplicable.
By (iiia) and the definition of $s_\L$, $(H',f|_{H'}) \wedge_{s+t} \L_{H'} \not \in \G$.
But $\G$ is hereditary, so it contradicts (iiib).
Hence Conclusion 2 follows from (ii).

Now we prove Conclusion 3.
Assume that every member of $\HH$ is injective. 
By (iib), it suffices to show that every member of $\L$ has order at most $r-1$.
Suppose to the contrary that there exists $(X,Y) \in \L$ of order at least $r$.
Let $(A^*,B^*) = (N_{H''}[\phi_H(X-Y)], V(H'')-\phi_H(X-Y))$.
Note that $(A^*,B^*)$ is a member of $\L_{H''}$, so it is a separation of $(H'',f_{H''})$ of order at most $r-1$ by (iv).
Since $A^* \cap B^* = N_{H''}(\phi_H(X-Y))$, $\lvert N_{H''}(\phi_H(X-Y)) \rvert = r-1<r \leq \lvert X \cap Y \rvert = \lvert N_H(X-Y) \rvert$, where the last equality follows from (iia).
Since $\phi_H \in \HH$, $\phi_H$ is injective.
So $\lvert \phi_H(N_H(X-Y)) \rvert \geq \lvert N_H(X-Y) \rvert > \lvert N_{H''}(\phi_H(X-Y)) \rvert$.
Since $\phi_H$ is an injective homomorphism from $(H,f_H)$ to $(H'',f_{H''})$ by (iv), $\phi_H(N_H(X-Y)) \subseteq N_{H''}(\phi_H(X-Y))$.
Hence $\lvert \phi_H(N_H(X-Y)) \rvert \leq \lvert N_{H''}(\phi_H(X-Y)) \rvert$, a contradiction.
This proves Conclusion 3.
\end{pf}

\bigskip

Now we are ready to prove the main theorem of this paper.
Recall that $\ex(H,f_H,\HH,\G,n)$ is the maximum number of homomorphisms in $\HH$ from $(H,f_H)$ to $(G,f_G)$, among all members $(G,f_G)$ of $\G$ on $n$ vertices; and $\dup_\HH(H,f_H,\G)$ is the maximum size of a $(\G,\HH)$-duplicable independent collection of separations of $(H,f_H)$.

\begin{theorem} \label{main}
For any positive integers $r,h$, there exist positive real numbers $r_0$ and $r_1$ such that for any quasi-order $Q$ with finite ground set, positive integer $t$, hereditary class $\G$ of legal $Q$-labelled graphs, legal $Q$-labelled graph $(H,f_H)$ on $h$ vertices, consistent set $\HH$ of homomorphisms from $(H,f_H)$ to $\G$, and positive real numbers $k,k'$, there exist positive real numbers $c_1$ and $c_2$ such that if there exist nonnegative real numbers $\epsilon,\epsilon'$ with $\epsilon \leq \frac{1}{20h^2}$ and $\epsilon' \leq r_1$ such that any member $(G,f_G)$ of $\G$ satisfies 
	\begin{enumerate}
		\item for every induced subgraph $L$ of $G$, $\lvert E(L) \rvert \leq k\lvert V(L) \rvert^{1+\epsilon}$,
		\item for every graph $L$ in which some subgraph of $G$ is isomorphic to a $[4h+1]$-subdivision of $L$, $\lvert E(L) \rvert \leq k'\lvert V(L) \rvert^{1+\epsilon'}$, and
		\item there exists no $(1,2hr+1)$-model for a $2h$-shallow $K_{r,t}$-minor in $G$. 
	\end{enumerate}
then
	\begin{enumerate}
		\item $\ex(H,f_H,\HH,\G,n) \geq c_1n^{\dup_\HH(H,f_H,\G)}$ for infinitely many positive integers $n$,
		\item $\ex(H,f_H,\HH,\G,n) \leq c_2n^{\dup_\HH(H,f_H,\G)+(1+\epsilon-\frac{1}{\epsilon_0})h^2}$ for all positive integers $n$, where $\epsilon_0=(1+\epsilon'(r-1))^{r_0}$, and
		\item if every member of $\HH$ is injective, then there exists a $(\G,\HH)$-duplicable independent collection $\L$ of separations of $H$ with $\lvert \L \rvert = \dup_\HH(H,f_H,\G)$ such that for every $(X,Y) \in \L$, $\lvert X \cap Y \rvert \leq r-1$ and $H[X-Y]$ is connected.
	\end{enumerate}
\end{theorem}

\begin{pf}
Let $r$ and $h$ be positive integers.
Define $r_0=r_{\ref{shallow_density_dup}}(r,h)$ and $r'_1=r_{\ref{shallow_density_dup}}'(r,h)$, where $r_{\ref{shallow_density_dup}}$ and $r'_{\ref{shallow_density_dup}}$ are the real numbers $r_0$ and $r_1$ mentioned in Lemma \ref{shallow_density_dup}, respectively.
Define $r_1$ to be $\min\{r_1', \frac{(\frac{1}{1-\frac{1}{2h^2}})^{\frac{1}{r_0}}-1}{r-1}\}$.
Let $Q$ be a quasi-order with finite ground set.
Let $t$ be a positive integer.
Let $\G$ be a hereditary class of legal $Q$-labelled graphs.
Let $(H,f_H)$ be a legal $Q$-labelled graph in $h$ vertices.
Let $\HH$ be a consistent set of homomorphisms from $(H,f_H)$ to $\G$.
Let $k$ and $k'$ be positive real numbers.
Define $c_1 = c_{\ref{lower_bound}}(H)$, where $c_{\ref{lower_bound}}$ is the real number mentioned in Proposition \ref{lower_bound}.
Define $c_2= c_{\ref{shallow_density_dup}}(Q,r,t,h,\G,(H,f_H),\HH,k,k')$, where $c_{\ref{shallow_density_dup}}$ is the real number $c$ mentioned in Lemma \ref{shallow_density_dup}.

Let $\epsilon$ and $\epsilon'$ be nonnegative real numbers with $\epsilon \leq \frac{1}{20h^2}$ and $\epsilon' \leq r_1$ such that every member of $\G$ satisfies the conditions stated in the lemma.

Then Conclusion 1 of this lemma immediately follows from Proposition \ref{lower_bound}.

Let $n$ be a positive integer.
For every member $(G,f_G)$ of $\G$ on $n$ vertices, by Lemma \ref{shallow_density_dup}, there exists a nonnegative integer $q_G$ such that 
	\begin{enumerate}
		\item[(i)] there are at most $c_2n^{q_G+(1+\epsilon-\frac{1}{\epsilon_0})h^2}$ homomorphisms in $\HH$ from $(H,f_H)$ to $(G,f_G)$, where $\epsilon_0=(1+\epsilon'(r-1))^{r_0}$, 
		\item[(ii)] there exists a $(\G,\HH)$-duplicable independent collection $\L_G$ of separations of $(H,f_H)$ with $\lvert \L_G \rvert = q_G$, and 
		\item[(iii)] if every member of $\HH$ is injective, then for every separation $(X,Y)$ in $\L_G$, it has order at most $r-1$ and $H[X-Y]$ is connected. 
	\end{enumerate}
By (ii), $q_G \leq \dup_\HH(H,f_H,\G)$ for every $(G,f_G) \in \G$. 
So by (i), $\ex(H,f_H,\HH,\G,n) \leq c_2n^{\dup_\HH(H,f_H,\G)+(1+\epsilon-\frac{1}{\epsilon_0})h^2}$ homomorphisms in $\HH$ from $(H,f_H)$ to $(G,f_G)$.
Hence Conclusion 2 holds.

Since $\epsilon' \leq r_1$, $1-\frac{1}{\epsilon_0} \leq \frac{1}{2h^2}$.
So $(1+\epsilon-\frac{1}{\epsilon_0})h^2 \leq \frac{11}{20}<1$.
Note that $\dup_\HH(H,f_H,\G)$ and $q_G$ are integers for any $(G,f_G) \in \G$.
So by Conclusion 1 and 2 of this lemma, there exist a large integer $N$ with $\ex(H,f_H,\HH,\G,N) \geq c_1N^{\dup_\HH(H,f_H,\G)}$ and $(G^*,f^*) \in \G$ on $N$ vertices such that $q_{G^*}=\dup_\HH(H,f_H,\G)$. 
Define $\L=\L_{G^*}$.
Then Conclusion 3 follows from (iii) by taking $\L=\L_{G^*}$.
\end{pf}

\bigskip

Now we are ready to deduce Theorem \ref{bdd_expan_intro} from Theorem \ref{main}.
The following corollary immediately implies Theorem \ref{bdd_expan_intro}.

\begin{corollary} \label{cor_bdd_expan}
For any quasi-order $Q$ with finite ground set, hereditary class $\G$ of legal $Q$-labelled graphs with bounded expansion, legal $Q$-labelled graph $(H,f_H)$, consistent set $\HH$ of homomorphisms from $(H,f_H)$ to $\G$, there exist positive integers $c_1$ and $c_2$ such that
	\begin{enumerate}
		\item $\ex(H,f_H,\HH,\G,n) \geq c_1n^{\dup_\HH(H,f_H,\G)}$ for infinitely many positive integers $n$,
		\item $\ex(H,f_H,\HH,\G,n) \leq c_2n^{\dup_\HH(H,f_H,\G)}$ for every positive integer $n$, and
		\item if every member of $\HH$ is injective and there exist positive integers $r$ and $t$ such that every graph in $\G$ has no $(1,2r\lvert V(H) \rvert+1)$-model for a $2\lvert V(H) \rvert$-shallow $K_{r,t}$-minor, then there exists a $(\G,\HH)$-duplicable independent collection $\L$ of separations of $(H,f_H)$ of order at most $r-1$ with $\lvert \L \rvert =\dup_\HH(H,f_H,\G)$ such that $H[X-Y]$ is connected for every $(X,Y) \in \L$.
	\end{enumerate}
\end{corollary}

\begin{pf}
Let $Q$ be a quasi-order with finite ground set.
Let $\G$ be a hereditary class of legal $Q$-labelled graphs with bounded expansion.
So there exists a function $f: {\mathbb N} \cup \{0\} \rightarrow {\mathbb N} \cup \{0\}$ such that $\G$ has expansion at most $f$.
Let $(H,f_H)$ be a legal $Q$-labelled graph.
Let $\HH$ be a consistent set of homomorphisms from $(H,f_H)$ to $\G$.
Let $h$ be the number of vertices of $H$.
Define $c_1$ and $c_2$ to be the real numbers $c_1$ and $c_2$, respectively, mentioned in Theorem \ref{main} by taking $(r,h,Q,t,\G,(H,f_H),\HH,k,k')=(f(2h)+1,h,Q,f(2h)+1, \G, (H,f_H), \HH,f(0),f(2h+1))$.

If there exists $(G,f_G) \in \G$ such that there exists a $(1,2h \cdot (f(2h)+1) +1)$-model for a $2h$-shallow $K_{f(2h)+1,f(2h)+1}$-minor in $G$, then $G$ contains some graph with average degree greater than $f(2h)$ as a $2h$-shallow minor, a contradiction.
Then applying Theorem \ref{main} by taking $\epsilon=\epsilon'=0$, Statements 1 and 2 immediately follows from Conclusions 1 and 2 of Theorem \ref{main}.
If there exist positive integers $r$ and $t$ such that every graph in $\G$ has no $(1,2r\lvert V(H) \rvert+1)$-model for a $2\lvert V(H) \rvert$-shallow $K_{r,t}$-minor, then we can apply Theorem \ref{main} by taking $(r,t)=(r,t)$, so Statement 3 of this corollary follows from Conclusion 3 of Theorem \ref{main}.
\end{pf}

\bigskip

Another consequence of Theorem \ref{main} is Theorem \ref{nowhere_dense_intro}.
We will use the following characterization for the classes of nowhere dense graphs.

\begin{theorem}[{\cite[Theorem 4.1]{no_nowhere_dense}}] \label{nowhere_chara}
Let $\C$ be a hereditary class of graphs.
For every nonnegative integer $\ell$, let $M_\ell$ be the set consisting of all graphs $L$ such that some graph in $\C$ contains $L$ as an $\ell$-shallow minors, and let $S_\ell$ be the set consisting of all graphs $L$ such that some graph in $\C$ contains a $([2\ell] \cup \{0\})$-subdivision of $L$.
Then following are equivalent.
	\begin{enumerate}
		\item $\C$ is nowhere dense.
		\item For every positive integer $k$, there exists an integer $g(k)$ such that the size of every clique in any graph in $M_k$ is at most $g(k)$.
		\item $\lim_{\ell \to \infty}\limsup_{L \in M_\ell} \frac{\log\lvert E(L) \rvert}{\log\lvert V(L) \rvert} \leq 1$.
		\item For every nonnegative integer $\ell$, $\limsup_{L \in M_\ell} \frac{\log\lvert E(L) \rvert}{\log\lvert V(L) \rvert} \leq 1$.
		\item $\lim_{\ell \to \infty}\limsup_{L \in S_\ell} \frac{\log\lvert E(L) \rvert}{\log\lvert V(L) \rvert} \leq 1$.
		\item For every nonnegative integer $\ell$, $\limsup_{L \in S_\ell} \frac{\log\lvert E(L) \rvert}{\log\lvert V(L) \rvert} \leq 1$.
	\end{enumerate}
\end{theorem}

Note that Statements 4 and 6 are not stated in \cite[Theorem 4.1]{no_nowhere_dense}, but they are equivalent to Statements 3 and 5, respectively, since $M_0 \subseteq M_1 \subseteq M_2 \subseteq ...$ and $S_0 \subseteq S_1 \subseteq S_2 \subseteq ...$ by the definition of $M_\ell$ and $S_\ell$.

Now we are ready to prove Theorem \ref{nowhere_dense_intro}.
The following is a restatement of Theorem \ref{nowhere_dense_intro}.

\begin{corollary} \label{nowhere_dense}
Let $Q$ be a quasi-order $Q$ with finite ground set.
Let $(H,f_H)$ be a legal $Q$-labelled graph.
Let $\G$ be a hereditary nowhere dense class of legal $Q$-labelled graphs.
Let $\HH$ be a consistent set of homomorphisms from $(H,f_H)$ to $\G$.
Then the asymptotic logarithmic density for $(H,f_H,\G,\HH)$ equals $\dup_\HH(H,f_H,\G)$.
\end{corollary}

\begin{pf}
By Proposition \ref{lower_bound}, the asymptotic logarithmic density for $(H,f_H,\G,\HH)$ is at least $\dup_\HH(H,f_H,\G)$.
To prove that the asymptotic logarithmic density for $(H,f_H,\G,\HH)$ is at most $\dup_\HH(H,f_H,\G)$, it suffices to show that it is at most $\dup_\HH(H,f_H,\G)+\delta$ for every $\delta>0$. 

Let $h$ be the number of vertices of $H$.
By Theorem \ref{nowhere_chara}, there exists a positive integer $r$ such that no member of $\G$ contains $K_r$ as a $(6h+1)$-shallow minor.
Since $K_{r,r}$ contains $K_r$ as a 1-shallow minor, no member of $\G$ contains $K_{r,r}$ as a $2h$-shallow minor.
Let $r_0$ and $r_1$ be the positive real numbers $r_0$ and $r_1$, respectively, mentioned in Theorem \ref{main} by taking $(r,h)=(r,h)$.

Let $\delta$ be a positive real number with $\delta<1$.
Let $\epsilon'=\min\{r_1, \frac{1}{(r-1)(1-\frac{\delta}{2h^2})^{1/r_0}}-\frac{1}{r-1}\}$.
Let $\epsilon_0=(1+\epsilon'(r-1))^{r_0}$.
So $(1-\frac{1}{\epsilon_0})h^2 \leq \frac{\delta}{2}$.
Let $\epsilon = \min\{\frac{1}{20h^2}, \frac{\delta}{2h^2}\}$.
Hence $(1+\epsilon-\frac{1}{\epsilon_0})h^2 \leq \delta$.

For every nonnegative integer $\ell$, let $M_\ell$ be the set consisting of all graphs $L$ such that some graph in $\G$ contains $L$ as an $\ell$-shallow minors, and let $S_\ell$ be the set consisting of all graphs $L$ such that some graph in $\G$ contains a $([2\ell] \cup \{0\})$-subdivision of $L$.
By Theorem \ref{nowhere_chara}, there exists a positive integer $N$ such that $\frac{\log \lvert E(L) \rvert}{\log \lvert V(L) \rvert} \leq 1+\min\{\epsilon,\epsilon'\}$ for every $L \in M_0 \cup S_{2h+1}$ with $\lvert V(L) \rvert \geq N$.
Hence for every $L \in M_0 \cup S_{2h+1}$, $\lvert E(L) \rvert \leq N^2 \lvert V(L) \rvert^{1+\min\{\epsilon,\epsilon'\}}$.
Let $c$ be the real number $c_2$ mentioned in Theorem \ref{main} by taking $(r,h,Q,t,\G,(H,f_H),\HH,k,k')=(r,h,Q,r,\G,(H,f_H),\HH,N^2,N^2)$.

Let $(G,f_G)$ be a member of $\G$. 
Then for every induced subgraph $L$ of $G$, $L \in M_0$, so $\lvert E(L) \rvert \leq N^2\lvert V(L) \rvert^{1+\epsilon}$.
For every graph $L$ in which some subgraph of $G$ is isomorphic to a $[4h+1]$-subdivision of $L$, $L \in S_{2h+1}$, so $\lvert E(L) \rvert \leq N^2\lvert V(L) \rvert^{1+\epsilon'}$.
Recall that no member of $\G$ contains $K_{r,r}$ as a $2h$-shallow minor, so there exists no $(1,2hr+1)$-model for a $2h$-shallow $K_{r,r}$-minor in $G$.

So by Theorem \ref{main}, $\ex(H,f_H,\HH,\G,n) \leq cn^{\dup_{\HH}(H,f_H,\G)+(1+\epsilon-\frac{1}{\epsilon_0})h^2} \leq cn^{\dup_{\HH}(H,f_H,\G)+\delta}$ for every positive integer $n$.
Hence the asymptotic logarithmic density for $(H,f_H,\G,\HH)$ is at most $\dup_\HH(H,f_H,\G)+\delta$.

Since $\delta$ is an arbitrary small positive real number, the asymptotic logarithmic density for $(H,f_H,\G,\HH)$ is at most $\dup_\HH(H,f_H,\G)$.
This proves the corollary.
\end{pf}

\section{Concrete applications} \label{sec:concrete}

In this section, we apply Theorem \ref{bdd_expan_intro} (or precisely, Corollary \ref{cor_bdd_expan}) to solve open questions on some extensively studied graph classes of (unlabelled) graphs.
We will concentrate on determining the number of (induced or not) subgraphs isomorphic to $H$ and the number homomorphisms from $H$ (without other conditions).

We first show that Corollary \ref{cor_bdd_expan} implies that the subgraph counts determine these three quantities up to a constant factor for monotone classes with bounded expansion.

A class $\G$ of graphs is {\it monotone} if every subgraph of any member of $\G$ belongs to $\G$.
Note that every monotone class is hereditary.
For any graph $H$, class of graphs $\G$ and integer $n$, we define 
	\begin{itemize}
		\item $\ex(H,\G,n)$ to be the maximum number of subgraphs of $G$ isomorphic to $H$, among all $n$-vertices graphs $G$ in $\G$,
		\item $\ind(H,\G,n)$ to be the maximum number of induced subgraphs of $G$ isomorphic to $H$, among all $n$-vertices graphs $G$ in $\G$,
		\item $\hom(H,\G,n)$ to be the maximum number of homomorphisms from $H$ to $G$, among all $n$-vertices graphs $G$ in $\G$, and
	\end{itemize}

\begin{corollary} \label{3_sub_equiv}
Let $H$ be a graph.
Let $\G$ be a monotone class of graphs with bounded expansion.
Then 
	\begin{enumerate}
		\item $\ind(H,\G,n) = \Theta(\ex(H,\G,n))$.
		\item $\hom(H,\G,n) = \max_{H'} \Theta(\ex(H',\G,n))$, where the maximum is over all graphs $H'$ with $\lvert V(H') \rvert \leq \lvert V(H) \rvert$ such that there exists an onto homomorphism from $H$ to $H'$.
	\end{enumerate}
\end{corollary}

\begin{pf}
Let $\HH_1$ be the set of all injective homomorphisms from $H$ to $\G$.
Note that $\ex(H,\G,n) = \Theta(\max_{G \in \G,\lvert V(G) \rvert=n} \lvert \{\phi \in \HH_1: \phi$ is from $H$ to $G\} \rvert)$. 
Let $\HH_2 = \{\phi \in \HH_1: \phi$ maps any pair of non-adjacent vertices to a pair of non-adjacent vertices$\}$.
Note that $\ind(H,\G,n) = \Theta(\max_{G \in \G, \lvert V(G) \rvert=n} \lvert \{\phi \in \HH_2: \phi$ is from $H$ to $G\} \rvert)$. 
Since $\G$ is monotone, an independent collection $\L$ of separations of $H$ is $(\G,\HH_1)$-duplicable if and only if $\L$ is $(\G,\HH_2)$-duplicable.
Since $\HH_1$ and $\HH_2$ are consistent, and since every monotone class of graphs is hereditary, Statement 1 of this proposition follows from Statements 1 and 2 of Corollary \ref{cor_bdd_expan}.

Let $S$ be the set of all isomorphism classes of graphs $H'$ with $\lvert V(H') \rvert \leq \lvert V(H) \rvert$ such that there exists an onto homomorphism from $H$ to $H'$.
For every homomorphism $\phi$ from $H$ to a graph $G$, let $H_\phi$ be the subgraph of $G$ induced by $\phi(V(H))$, so $H_\phi \in S$.
So $\hom(H,\G,n) \leq h^h\sum_{H' \in S}\ind(H',\G,n)$.
Since $\lvert S \rvert$ is upper bounded by a function only depending on $H$, we have $\hom(H,\G,n) =\max_{H' \in S}O(\ind(H',\G,n))$.
On the other hand, $\hom(H,\G,n) \geq \ind(H',\G,n)$ for every $H' \in S$.
So $\hom(H,\G,n) =\max_{H' \in S}\Theta(\ind(H,\G,n))$.
Hence Statement 2 follows from Statement 1.
\end{pf}

\bigskip

By Corollary \ref{3_sub_equiv}, we can concentrate on $\ex(H,\G,n)$.
For a graph class $\G$ and a graph $H$, we say that an independent collection is {\it $\G$-duplicable} if it is $(\G,\HH)$-duplicable, where $\HH$ is the collection of all injective homomorphisms from $H$ to $\G$.
So it suffices to understand the maximum size of a $\G$-duplicable collection by Corollary \ref{cor_bdd_expan}.

\subsection{Preparation}

A collection $\L$ of separations of a graph $H$ is {\it essential} if for every $(A,B) \in \L$, every vertex in $A \cap B$ is adjacent in $H$ to some vertex in $A-B$, and $H[A-B]$ is connected.

\begin{prop} \label{make_essential}
Let $\G$ be a monotone class of graphs.
Let $H$ be a graph.
Let $k$ be a nonnegative integer.
If $\L$ is a $\G$-duplicable independent collection of separations of $H$ such that every separation in $\L$ has order at most $k$, then there exists an essential $\G$-duplicable independent collection $\L'$ of separations of $H$ with $\lvert \L' \rvert \geq \lvert \L \rvert$ such that every separation in $\L'$ has order at most $k$.
\end{prop}

\begin{pf}
For every $(A,B) \in \L$, let $\L_A = \{(N_H[V(C)], V(H)-V(C)): C$ is a component of $H[A-B]\}$. 
Clearly, $\L_A$ is an essential collection of separations of $H$ of order at most $\lvert A \cap B \rvert$.
Let $\L' = \bigcup_{(A,B) \in \L}\L_A$.
Since $\L$ is independent, $\lvert \L' \rvert \geq \lvert \L \rvert$ and $\L'$ is independent.
Since $\G$ is monotone and $\L$ is $\G$-duplicable, $\L'$ is $\G$-duplicable.
\end{pf}

\begin{prop} \label{easy_dup_subdiv}
Let $\G$ be a monotone class of graphs.
Let $H$ be a graph.
Let $\L$ be an essential $\G$-duplicable independent collection of separations of $H$.
Let $H'$ be the graph obtained from $H$ by adding edges such that $A \cap B$ is a clique for every $(A,B) \in \L$.
Then some $([\lvert V(H) \rvert-2] \cup \{0\})$-subdivision of $H'$ belongs to $\G$.
\end{prop}

\begin{pf}
Let $w={\lvert V(H) \rvert \choose 2}+1$.
Since $\L$ is $\G$-duplicable and $\G$ is monotone, $H \wedge_w \L \in \G$.
Since $\L$ is essential, for every $(A,B) \in \L$ and for every $u,v \in A \cap B$, there exists a path in $H[A]$ on at most $\lvert V(H) \rvert$ vertices from $u$ to $v$ whose all internal vertices are in $A-B$.
Since $w={\lvert V(H) \rvert \choose 2}+1$, $H \wedge_w \L$ contains a $([\lvert V(H) \rvert-2] \cup \{0\})$-subdivision of $H'$.
This proposition follows since $\G$ is monotone.
\end{pf}

\bigskip

Let $\L$ be an independent collection of separations of a graph $H$.
A {\it torso} of $\L$ is a graph that is either obtained from $H[\bigcap_{(A,B) \in \L}B]$ by adding edges such that $A \cap B$ is a clique for every $(A,B) \in \L$, or obtained from $H[X]$ for some $(X,Y) \in \L$ by adding edges such that $X \cap Y$ is a clique.

We say that a graph $H$ is a {\it topological minor} of another graph $G$ if some subgraph of $G$ is isomorphic to a subdivision of $H$.
A class $\G$ of graphs is {\it topological minor-closed} if every topological minor of any member of $\G$ is a member of $\G$.

\begin{prop} \label{topo_torso}
Let $\G$ be a topological minor-closed class of graphs.
Let $H$ be a graph.
Let $\L$ be an essential $\G$-duplicable independent collection of separations of $H$.
Then every torso of $\L$ belongs to $\G$.
\end{prop}

\begin{pf}
Let $H'$ be the graph obtained from $H$ by adding edges such that $A \cap B$ is a clique for every $(A,B) \in \L$.
By Proposition \ref{easy_dup_subdiv}, $H' \in \G$ since $\G$ is topological minor-closed.
This proposition follows since every torso of $\L$ is a subgraph of $H'$.
\end{pf}

\bigskip

Let $t$ be a nonnegative integer.
Let $G_1$ and $G_2$ be graphs.
If $G_1$ contains a clique $\{u_1,u_2,...,u_t\}$ of size $t$, and $G_2$ contains a clique $\{v_1,v_2,.,,,v_t\}$ of size $t$, then a {\it $t$-sum} of $G_1$ and $G_2$ is a graph obtained from the disjoint union of $G_1$ and $G_2$ by, for each $i \in [t]$, identifying $u_i$ and $v_i$ into a new vertex $w_i$, and then deleting any number of edges whose both ends are in $\{w_1,w_2,...,w_t\}$.
A {\it $(\leq t)$-sum} of $G_1$ and $G_2$ is a $k$-sum of $G_1$ and $G_2$ for some integer $k$ with $0 \leq k \leq t$.
Note that for every independent collection $\L$ of separations of a graph $H$, $H$ is a $(\leq t)$-sum of all torsos of $\L$, where $t = \max_{(A,B) \in \L}\lvert A \cap B \rvert$.
A class $\G$ of graphs is {\it closed under $(\leq t)$-sum} if any $(\leq t)$-sum of members of $\G$ belongs to $\G$.

\begin{prop} \label{clique_sum_dup}
Let $\G$ be a monotone class of graphs.
Let $H$ be a graph.
Let $t$ be a nonnegative integer.
Let $\L$ be an independent collection of separations of $H$ such that every separation in $\L$ has order at most $t$.
If $\G$ is closed under $(\leq t)$-sum, and every torso of $\L$ belongs to $\G$, then $\L$ is $\G$-duplicable.
\end{prop}

\begin{pf}
For every positive integer $w$, $H \wedge_w \L$ is a subgraph of a $(\leq t)$-sum of torsos of $\L$ and hence belongs to $\G$.
So $\L$ is $\G$-duplicable.
\end{pf}

\begin{corollary} \label{dup_torso}
Let $t$ be a positive integer.
Let $\G$ be a topological minor-closed class of graphs closed under $(\leq t-1)$-sum.
Let $\L$ be an essential independent collection of separations of a graph $H$ such that every member of $\L$ has order at most $t-1$.
Then $\L$ is $\G$-duplicable if and only if every torso of $\L$ belongs to $\G$.
\end{corollary}

\begin{pf}
It is an immediate corollary of Propositions \ref{topo_torso} and \ref{clique_sum_dup}.
\end{pf}

\begin{prop} \label{minor_sum}
Let $t$ be a positive integer.
Let $\preceq$ be the minor relation or the topological minor relation.
Let $\G$ be a $\preceq$-closed class of graphs.
Let $\F$ be the set of $\preceq$-minimal graphs that do not belong to $\G$.
If every member of $\F$ is $t$-connected, then $\G$ is closed under $(\leq t-1)$-sum.
\end{prop}

\begin{pf}
Let $G$ be a $(\leq t-1)$-sum of graphs $G_1,G_2 \in \G$.
So there exists a separation $(A,B)$ of $G$ with order at most $t-1$ such that $V(A)=V(G_1)$ and $V(B)=V(G_2)$.
Suppose to the contrary that $G \not \in \G$.
Then $H \preceq G$ for some graph $H \in \F$.
So $H$ is $t$-connected.
Since $\lvert A \cap B \rvert \leq t-1$ and $A \cap B$ is a clique in $G_1$ and a clique in $G_2$, either $H \preceq G_1$ or $H \preceq G_2$. 
So one of $G_1,G_2$ is not in $\G$, a contradiction.
\end{pf}

\begin{corollary} \label{minor_topo_ex}
Let $r$ be a positive integer.
Let $\preceq$ be the minor relation or the topological minor relation.
Let $\G$ be a $\preceq$-closed class of graphs.
Let $\F$ be the set of $\preceq$-minimal graphs that do not belong to $\G$.
If every member of $\F$ is $r$-connected, and there exists a positive integer $t$ such that every graph in $\G$ has no $(1,2r\lvert V(H) \rvert+1)$-model for a $2\lvert V(H) \rvert$-shallow $K_{r,t}$-minor, then $\ex(H,\G,n) = \Theta(n^k)$, where $k$ is the maximum size of an essential independent collection of separations of $H$ of order at most $r-1$ whose every torso is in $\G$.
\end{corollary}

\begin{pf}
By Corollary \ref{cor_bdd_expan}, $\ex(H,\G,n)=\Theta(n^k)$, where $k$ is the maximum size of a $\G$-duplicable independent collection of separations of $H$ of order at most $r-1$. 
By Proposition \ref{make_essential}, $k$ equals the maximum size of an essential $\G$-duplicable independent collection of separations of $H$ of order at most $r-1$.
By Proposition \ref{minor_sum}, $\G$ is closed under $(\leq r-1)$-sum. 
Since every minor-closed class is topological minor-closed, by Corollary \ref{dup_torso}, $k$ equals the maximum size of an essential independent collection of separations of $H$ of order at most $r-1$ whose every torso of is in $\G$. 
\end{pf}

\subsection{Minor-closed families}

We address minor-closed families, one of the most extensively studied sparse graph classes, in this subsection.

\subsubsection{Common minor-closed families} 

Huynh and Wood \cite{hw} proposed the following conjecture which was a potential answer of a question of Eppstein \cite{e}.
(Recall that $\flap_d(H)$ is the maximum size of an independent collection of separations of $H$ of order at most $d$.)

\begin{conj}[{\cite[Conjecture 15]{hw}}] \label{Kst_minor_conj}
Let $s,t$ be positive integers with $s \leq t$.
Let $\G$ be the class of graphs containing no $K_{s,t}$-minor.
Then for every $K_{s,t}$-minor free graph $H$, $\ex(H,\G,n)=\Theta(n^{\flap_{s-1}(H)})$.
\end{conj}

The following corollary disproves Conjecture \ref{Kst_minor_conj} by providing the correct bound, as it is easy to see that the bound mentioned in the following corollary is not equal to $\flap_{s-1}(H)$ for infinitely many graphs $H$. 

\begin{corollary} \label{cor_Kst_minor}
Let $s,t$ be positive integers with $s \leq t$.
Let $\G$ be the class of graphs containing no $K_{s,t}$-minor.
Then for every $K_{s,t}$-minor free graph $H$, $\ex(H,\G,n)=\Theta(n^k)$, where $k$ is the maximum size of an essential independent collection $\L$ of separations of $H$ of order at most $s-1$ such that every torso of $\L$ is $K_{s,t}$-minor free.
\end{corollary}

\begin{pf}
Since $\G$ is $K_{s,t}$-minor free, every graph in $\G$ has no $(1,2s\lvert V(H) \rvert+1)$-model for a $2\lvert V(H) \rvert$-shallow $K_{s,t}$-minor.
Since $K_{s,t}$ is $s$-connected, this corollary immediately follows from Corollary \ref{minor_topo_ex}.
\end{pf}

\bigskip

Corollary \ref{cor_Kst_minor} generalizes a result of Eppstein \cite{e} who proved the case that $H$ is $s$-connected, a result of Huynh, Joret and Wood \cite{hjw} who proved the case that $s \leq 3$, and a result of Huynh and Wood \cite{hw} who proved the case that $H$ is a tree.
Note that the result in \cite{hw} for trees is described by using $\alpha_{s-1}(H)$, but it can be easily shown that it matches the result in Corollary \ref{cor_Kst_minor} by counting the number of edges and considering $(\leq 1)$-sums.

A similar result for $K_t$-minor free graphs can be proved in a similar way.

\begin{corollary} \label{cor_Kt_minor}
Let $t$ be a positive integer.
Let $\G$ be the class of graphs containing no $K_t$-minor.
Then for every $K_t$-minor free graph $H$, $\ex(H,\G,n)=\Theta(n^k)$, where $k$ is the maximum size of an essential independent collection $\L$ of separations of $H$ of order at most $t-2$ such that every torso of $\L$ is $K_{t}$-minor free.
\end{corollary}

\begin{pf}
Since every graph in $\G$ is $K_{t}$-minor free, every graph in $\G$ has no $2\lvert V(H) \rvert$-shallow $K_{t-1,t-1}$-minor.
Since $K_{t}$ is $(t-1)$-connected, this corollary immediately follows from Corollary \ref{minor_topo_ex}.
\end{pf}

\bigskip

Similar arguments lead to results for the class of bounded tree-width graphs, one of the most extensively studied minor-closed families.
The {\it tree-width} of a graph $G$ is the minimum $k$ such that $G$ is a subgraph of a chordal graph with maximum clique size $k+1$.

\begin{corollary} \label{cor_tw}
Let $t$ be a positive integer.
Let $\G$ be the class of graphs of tree-width at most $t$.
Then for every graph $H$ of tree-width at most $t$, $\ex(H,\G,n)=\Theta(n^k)$, where $k$ is the maximum size of an essential independent collection $\L$ of separations of $H$ of order at most $t$ such that every torso of $\L$ has tree-width at most $t$.
\end{corollary}

\begin{pf}
Since $K_{t+1,t+1}$ has tree-width $t+1$, every graph in $\G$ has no $2\lvert V(H) \rvert$-shallow $K_{t+1,t+1}$-minor.
By Corollary \ref{cor_bdd_expan}, $\ex(H,\G,n)=\Theta(n^k)$, where $k$ is the maximum size of a $\G$-duplicable independent collection of separations of $H$ of order at most $t$.
By Proposition \ref{make_essential}, $k$ equals the maximum size of an essential $\G$-duplicable independent collection of separations of $H$ of order at most $t$.
By the definition of tree-width, $\G$ is closed under $(\leq t)$-sum. 
By Corollary \ref{dup_torso}, $k$ equals the maximum size of an essential independent collection of separations of $H$ of order at most $t$ whose every torso is in $\G$.
\end{pf}

\bigskip

The special case that $H$ is a tree in Corollaries \ref{cor_Kt_minor} and \ref{cor_tw} were proved by Huynh and Wood \cite{hw}.

Another common minor-closed family is the class of graphs of bounded path-width.
The {\it path-width} of a graph $G$ is the minimum $k$ such that there exists a function $f$ that maps each vertex of $G$ to an interval in ${\mathbb R}$ with positive integral endpoints such that $f(u) \cap f(v) \neq \emptyset$ for every edge $uv \in E(G)$, and $\lvert \{v \in V(G): x \in f(v)\} \rvert \leq k+1$ for every real number $x$; the function $f$ is called a {\it path-decomposition of $G$ of width $k$}.
Huynh and Wood \cite{hw} asked the following question.

\begin{question}[\cite{hw}]
If $H$ is a tree, what is the maximum number of $H$-subgraphs in an $n$-vertex graph of path-width at most $t$.
\end{question}

Corollary \ref{cor_bdd_expan} solves this question, up to a constant factor, even when $H$ is not a tree.

\begin{corollary}
Let $t$ be a positive integer.
Let $\G$ be the class of graphs of path-width at most $t$.
Then for every graph $H$ of path-width at most $t$, $\ex(H,\G,n)=\Theta(n^k)$, where $k$ is the maximum size of an essential independent collection $\L$ of separations of $H$ of order at most $t$ such that $H \wedge_{{t \choose 2}+2t+3} \L$ has path-width at most $t$. 
\end{corollary}

\begin{pf}
Since $K_{t+1,t+1}$ has tree-width $t+1$, it has path-width at least $t+1$.
So every graph in $\G$ has no $(1,2(t+1)\lvert V(H) \rvert+1)$-model for a $2\lvert V(H) \rvert$-shallow $K_{t+1,t+1}$-minor.
By Corollary \ref{cor_bdd_expan} and Proposition \ref{make_essential}, $\ex(H,\G,n)=\Theta(n^k)$, where $k$ is the maximum size of an essential $\G$-duplicable independent collection of separations of $H$ of order at most $t$.

It suffices to show that an essential independent collection $\L$ of separations of $H$ of order at most $t$ is $\G$-duplicable if and only if $H \wedge_{{t \choose 2}+2t+3} \L \in \G$.
Let $\L$ be an essential independent collection of separations of $H$ of order at most $t$.
It suffices to show that if $H \wedge_{{t \choose 2}+2t+3} \L \in \G$, then $\L$ is $\G$-duplicable by the definition of $\G$-duplicable collections.

Assume that $H \wedge_{{t \choose 2}+2t+3} \L \in \G$.
Let $H'$ be the graph obtained from $H$ by adding edges such that $A \cap B$ is a clique in $H'$ for every $(A,B) \in \L$.
Note that $\L$ is also an essential independent collection of separations of $H'$ of order at most $t$.
Since $\L$ is essential, $H' \wedge_{2t+3} \L$ is a minor of $H \wedge_{{t \choose 2}+2t+3} \L$.
Since $\G$ is minor-closed, $H' \wedge_{2t+3} \L \in \G$.
So there exists a path-decomposition $f$ of $H' \wedge_{2t+3} \L$ with width at most $t$.

For each $(A,B) \in \L$, let $Y_{A,1},Y_{A,2},...,Y_{A,2t+3}$ be the copies of $A-B$ in $H'[A] \wedge_{2t+3} (A \cap B)$.
For each $(A,B) \in \L$ and $i \in [2t+3]$, let $I_{A,i}=\bigcup_{v \in Y_{A,i}}f(v)$; let $\C_A$ be the multiset $\{I_{A,j}: j \in [2t+3]\}$.
Since $\L$ is essential, $H'[A-B]$ is connected, so $I_{A,i}$ is an interval with integral endpoints for every $(A,B) \in \L$ and $i \in [2t+3]$; we let $a_{A,i}$ and $b_{A,i}$ be the left and right endpoints of $I_{A,i}$, respectively.
Since the width of $f$ is at most $t$, for every $(A,B) \in \L$, there do not exist $t+2$ pairwise intersecting members of $\C_A$ by Helly's property.
For each $(A,B) \in \L$, since $\lvert \C_A \rvert = 2(t+1)+1$, by repeatedly greedily choosing an interval in $\C_A$ with smallest $b_i$ among all intervals in $\C_A$ disjoint from the intervals that have been chosen, there exist three pairwise disjoint members in $\C_A$; by symmetry, we may assume that $I_{A,1},I_{A,2},I_{A,3}$ are pairwise disjoint, and $b_{A,1}<a_{A,2} \leq b_{A,2}<a_{A,3}$.
Since $\L$ is essential, for every $(A,B) \in \L$, every vertex in $A \cap B$ is adjacent to some vertex in $A-B$, so $I_{A,2}$ is contained in the interior of $f(v)$ for every $v \in A \cap B$.
For every $(A,B) \in \L$, let $j_A$ be an integer in $I_{A,2}$ such that $\lvert \{v \in Y_{A,2}: j_A \in f(v)\} \rvert$ is as large as possible.
Note that $j_A \in I_{A,2}$ is contained in the interior of $f(v)$ for every $v \in A \cap B$.

To show that $\L$ is $\G$-duplicable, it suffices to show that $H' \wedge_w \L \in \G$ for every integer $w$.
Let $w$ be an arbitrary positive integer.
We shall construct a path-decomposition of $H' \wedge_w \L$ of width at most $t$.
Note that it suffices to define a function $g$ that maps each vertex of $H' \wedge_w \L$ to a closed interval in ${\mathbb R}$ such that $g(u) \cap g(v) \neq \emptyset$ for every $uv \in E(H' \wedge_w \L)$, and $\lvert \{v \in V(H' \wedge_w \L): x \in g(v)\} \rvert \leq t+1$ for every $x \in {\mathbb R}$, as we can stretch all intervals to make them having integral endpoints without changing the intersection patterns.

For each $u \in \bigcap_{(X,Y) \in \L}Y$, define $g(u)=f(u)$.
For each $(A,B) \in \L$ and $v \in A-B$, let $v'$ be the copy of $v$ in $Y_{A,2}$ and denote $f(v')=[a_v,b_v]$, and let $f'(v) = [(a_v-a_{A,2}) \cdot \frac{1}{(b_{A,2}-a_{A,2}+1)2w},(b_v-a_{A,2})\cdot \frac{1}{(b_{A,2}-a_{A,2}+1)2w}]$.
That is, $f'(v)$ is obtained from $f(v')$ by first shifting and then scaling.
For each $(A,B) \in \L$, by the definition of $I_{A,2}$, we know that $\bigcup_{v \in A-B}f'(v)$ is an interval with left endpoint $0$ and with length less than $\frac{1}{2w}$.
For each $(A,B) \in \L$ and $v \in A-B$, denote $f'(v)$ by $[a'_v,b'_v]$, and for each $i \in [w]$, let $g$ map the $i$-th copy of $v$ to $[j_A+\frac{i-1}{w}+a'_v,j_A+\frac{i-1}{w}+b'_v]$.
Clearly, for distinct $i_1,i_2 \in [w]$, $g$ maps the $i_1$-th copy of $A-B$ and the $i_2$-th copy of $A-B$ into disjoint intervals contained in $[j_A,j_A+1-\frac{1}{2w}]$.
Then it is straight forward to show that $g(u) \cap g(v) \neq \emptyset$ for every $uv \in E(H' \wedge_w \L)$. 
And for every $x \in {\mathbb R}$, $\lvert \{v \in V(H' \wedge_w \L): x \in g(v)\} \rvert \leq t+1$ by the definition of $j_A$ for each $(A,B) \in \L$.
This proves the corollary.
\end{pf}

\subsubsection{Minor-closed families with topological properties} 

Another extensively studied minor-closed family is the class of linkless embeddable graphs and flat embeddable graphs.
A graph is {\it linkless embeddable} if it can be embedded in ${\mathbb R}^3$ such that any two disjoint cycles are unlinked (in the topology sense). 
A graph is {\it flat embeddable} if it can be embedded in ${\mathbb R}^3$ such that every cycle is the boundary of an open disk disjoint from the embedding.
(See \cite{rst_survey} for a survey about these two notions.)
It is easy to see that the classes of linkless embeddable or flat embeddable graphs are minor-closed.
Robertson, Seymour and Thomas \cite{rst_flat} proved that linkless embeddable graphs are equivalent to flat embeddable graphs and proved that a graph is linkless embeddable if and only if it does not contain any graph in the Petersen family as a minor.
The Petersen family consists of the graphs that can be obtained from $K_6$ by repeatedly applying $\Delta Y$-operations or $Y\Delta$-operations.
There are seven graphs in the Petersen family, and the Petersen graph is one of them.

Linkless embeddable graphs are related to the Colin de Verdi\`{e}re parameter.
Colin de Verdi\`{e}re \cite{c_para,c_para2} introduced a parameter $\mu(G)$ for any graph $G$, motivated by the study of the maximum multiplicity of the second eigenvalue of certain Schr\"{o}dinger operators.
It is defined to be the largest corank of certain matrices associated with $G$. 
For any integer $k$, the class of graphs with $\mu \leq k$ is minor-closed \cite{c_para}.
Many graphs with certain topological properties can be characterized by using this algebraic parameter.
For example, $\mu(G) \leq 1$ if and only if $G$ is a disjoint union of paths \cite{c_para}; $\mu(G) \leq 2$ if and only if $G$ is outerplanar \cite{c_para}; $\mu(G) \leq 3$ if and only if $G$ is planar \cite{c_para}; $\mu(G) \leq 4$ if and only if $G$ is linkless embeddable \cite{ls,rst_survey,rst_flat}.
See \cite{hls} for a survey.

\begin{corollary} \label{cor_para}
Let $t$ be a positive integer.
Let $\G$ be the class of all graphs $G$ with $\mu(G) \leq t$.
Let $H$ be a graph with $\mu(H) \leq t$.
Then $\ex(H,\G,n) = \Theta(n^k)$, where $k$ is the maximum size of an essential independent collection of separations of $H$ of order at most $t-1$ such that every its torso $L$ satisfies $\mu(L) \leq t$.
\end{corollary}

\begin{pf}
It is shown in \cite{hls} that $\mu(K_{t,t+3})=t+1$.
Since $\G$ is minor-closed, no graph in $\G$ contains a $2\lvert V(H) \rvert$-shallow $K_{t,t+3}$-minor.
So $\ex(H,\G,n) = \Theta(n^k)$ by Corollary \ref{cor_bdd_expan}, where $k$ is the maximum size of a $\G$-duplicable independent collection of separations of order at most $t-1$.
By \cite[Theorem 2.10, Corollary 2.11]{hls}, $\G$ is closed under $(\leq t-1)$-sum.
By Proposition \ref{make_essential} and Corollary \ref{dup_torso}, $k$ equals the maximum size of an essential independent collection of separations of $H$ of order at most $t-1$ whose every torso is in $\G$.
\end{pf}

\bigskip

Recall that $\mu(G) \leq 4$ if and only if $G$ is linkless embeddable \cite{ls,rst_survey,rst_flat}.
So Corollary \ref{cor_para} immediately implies the following result for linkless embeddable graphs.

\begin{corollary} \label{cor_linkless}
Let $\G$ be the class of linkless embeddable graphs.
Let $H$ be a linkless embeddable graph.
Then $\ex(H,\G,n) = \Theta(n^k)$, where $k$ is the maximum size of an essential independent collection of separations of $H$ of order at most $3$ such that every its torso is linkless embeddable.
\end{corollary}

Note that Corollary \ref{cor_linkless} can also be derived by using the characterization for linkless embeddable graphs by Robertson, Seymour and Thomas \cite{rst_flat} in terms of Petersen family minors without considering $\mu$, but we omit this direct proof.
The characterization in \cite{rst_flat} is also a key ingredient for showing the equivalence between $\mu \leq 4$ and linkless embeddability.

Recall that $\mu(G) \leq 3$ if and only if $G$ is planar \cite{c_para}.
So the case $\mu \leq 3$ in Corollary \ref{cor_para} immediately implies the result for the class of planar graphs.
The case for $3$-connected $H$ in this result was conjectured by Perles (see \cite{ac}) and independently proved by Wormald \cite{w_3conn} and Eppstein \cite{e}.
On the other hand, this result for planar graphs is a special case of a recent result of Huynh, Joret and Wood \cite{hjw} for graphs of bounded Euler genus.
In Section \ref{subsec:topo}, we will show that Corollary \ref{cor_bdd_expan} implies an even stronger result that allows crossings in drawings in surfaces of bounded Euler genus. %stating that $\ex(H,\G,n)$ for the class $\G$ of bounded Euler genus

Since $\mu(G) \leq 2$ if and only if $G$ is outerplanar \cite{c_para}, the case $\mu \leq 2$ in Corollary \ref{cor_para} implies the following result about outerplanar graphs.
(A {\it cut-vertex} of a graph is a vertex whose deletion increases the number of components.
A {\it block} of a graph $H$ is the maximal connected subgraph $L$ of $H$ such that $L$ has no cut-vertex.
A {\it end-block} of $H$ is a block of $H$ containing at most one cut-vertex of $H$.)

\begin{corollary} \label{cor_outerplanar}
Let $\G$ be the class of outerplanar graphs.
Let $H$ be an outerplanar graph.
Then $\ex(H,\G,n) = \Theta(n^k)$, where $k$ is the number of end-blocks of $H$.
\end{corollary}

\begin{pf}
Since $\mu(G) \leq 2$ if and only if $G$ is outerplanar, Corollary \ref{cor_para} implies that $\ex(H,\G,n) \allowbreak = \Theta(n^k)$, where $k$ is the maximum size of an essential independent collection of separations of $H$ of order at most $1$ such that every its torso $L$ is outerplanar.
Note that the condition for torsos can be dropped since the maximum order of separations considered here is at most 1.
And it is easy to see that the maximum size of an essential independent collection of separations of $H$ of order at most $1$ equals the number of end-blocks of $H$ by seeing the block-tree.
\end{pf}

\bigskip

The special case that $H$ is 2-connected in Corollary \ref{cor_outerplanar} was proved by Eppstein \cite{e}.
Note that outerplanar graphs are exactly the graphs with stack number 1.
We will consider graphs with bounded stack number in Section \ref{subsubsec:stack_queue}.

A graph is {\it knotless embeddable} if it can be embedded in ${\mathbb R}^3$ such that every cycle forms a trivial knot.
Knotless embeddable graphs share some features with linkless embeddable graphs, but they are significantly more complicated.
For example, even though the class of knotless embeddable graphs is a minor-closed family, the complete list of minor-minimal non-knotless embeddable graphs remains unknown.
In contrast to the known complete list of 7 graphs in the Petersen family for linkless embeddable graphs, the current incomplete list of minor-minimal non-knotless embeddable graphs contains more than 200 graphs \cite{gmn}, including $K_7$ \cite{cg} and $K_{3,3,1,1}$ \cite{f_knot}. 
See \cite{fmmnn} for a survey.

Corollary \ref{cor_bdd_expan} can be applied to the class of knotless embeddable graphs as well.
However, since structures of knotless embeddable graphs are far from being well-understood, and in particular it seems unknown whether the class of knotless embeddable graphs is closed under $(\leq 4)$-sum, we only obtain less explicit result.

\begin{corollary} \label{cor_knotless}
Let $\G$ be the class of knotless embeddable graphs.
Let $H$ be a knotless embeddable graph.
Then the following statements hold. 
	\begin{enumerate}
		\item $\ex(H,\G,n) = \Theta(n^k)$, where $k$ is the maximum size of an essential $\G$-duplicable independent collection of separations of $H$ of order at most $4$.
		\item $\ex(H,\G,n) = O(n^k)$, where $k$ is the maximum size of an essential independent collection of separations of $H$ of order at most $4$ whose torsos are knotless embeddable.
		\item If $\G$ is closed under $(\leq 4)$-sum, then $\ex(H,\G,n) = \Theta(n^k)$, where $k$ is the maximum size of an essential independent collection of separations of $H$ of order at most $4$ whose torsos are knotless embeddable.
	\end{enumerate}
\end{corollary}

\begin{pf}
By \cite{s_knot}, $K_{5,5} \not \in \G$.
Since $\G$ is minor-closed, no graph in $\G$ contains a $2\lvert V(H) \rvert$-shallow $K_{5,5}$-minor.
So $\ex(H,\G,n) = \Theta(n^k)$ by Corollary \ref{cor_bdd_expan}, where $k$ is the maximum size of a $\G$-duplicable independent collection of separations of order at most $4$.
By Proposition \ref{make_essential}, $k$ equals the maximum size of an essential $\G$-duplicable independent collection of separations of $H$ of order at most $4$.
So Statement 1 holds.
Since $\G$ is minor-closed, if $\L$ is an essential $\G$-duplicable independent collection of separations of $H$, then every torso of $\L$ belongs to $\G$.
So Statement 2 holds. 
Statement 3 immediately follows from Corollary \ref{dup_torso}.
\end{pf}

\bigskip

We remark that the special cases that $H$ is a tree for Corollaries \ref{cor_para}, \ref{cor_linkless} and \ref{cor_knotless} were proved by Huynh and Wood \cite{hw}.

\subsection{Topological minor-closed families} \label{subsec:topo}

We consider topological minor-closed but not minor-closed classes in this subsection.
The class of all graphs with no $K_t$-topological minor is a typical example.

\begin{corollary}
Let $t$ be a positive integer.
Let $\G$ be the class of graphs containing no $K_t$-topological minor.
Then for every $K_t$-topological minor free graph $H$, $\ex(H,\G,n)=\Theta(n^k)$, where $k$ is the maximum size of an essential independent collection $\L$ of separations of $H$ of order at most $t-1$ such that every torso of $\L$ is $K_{t}$-topological minor free.
\end{corollary}

\begin{pf}
Since every graph in $\G$ is $K_{t}$-topological minor free, every graph in $\G$ has no $(1,2t\lvert V(H) \rvert+1)$-model for a $2\lvert V(H) \rvert$-shallow $K_{t,{t \choose 2}}$-minor.
So by Corollary \ref{cor_bdd_expan} and Proposition \ref{make_essential}, $\ex(H,\G,n) = \Theta(n^k)$, where $k$ equals the maximum size of an essential $\G$-duplicable independent collection of separations of $H$ of order at most $t-1$.
To prove this corollary, it suffices to prove that $\G$ is closed under $(\leq t-1)$-sum by Corollary \ref{dup_torso}.

Suppose to the contrary that there exist $G_1,G_2 \in \G$ and an integer $q$ with $0 \leq q \leq t-1$ such that the $q$-sum $G$ of $G_1$ and $G_2$ contains a subdivision of $K_t$.
For each $i \in [2]$, let $S_i$ be the subset of $V(G_i)$ consisting of the $q$ vertices identified with vertices in $G_{3-i}$. 
Since $S_1$ is a clique in $G_1$ and $S_2$ is a clique in $G_2$, some branch vertex $v_1$ of this subdivision is in $V(G)-V(G_2)$ and some branch vertex is in $V(G)-V(G_1)$.
Since $K_t$ is $(t-1)$-connected, $q=t-1$ and there exist $t-1$ disjoint paths in $G_1$ from $v_1$ to $S_1$.
Since $S_1$ is a clique of size $t-1$ in $G_1$, $G_1$ contains a subdivision of $K_t$, a contradiction.
\end{pf}

\bigskip

$K_{3,t}$-topological minor free graphs are also of particular interests since every graph embeddable in a surface of Euler genus at most $g$ is $K_{3,2g+3}$-topological minor free.

\begin{corollary}
Let $s,t$ be positive integers with $s \leq 3$ and $s \leq t$.
Let $\G$ be the class of graphs containing no $K_{s,t}$-topological minor.
Then for every $K_{s,t}$-topological minor free graph $H$, $\ex(H,\G,n)=\Theta(n^k)$, where $k$ is the maximum size of an essential independent collection $\L$ of separations of $H$ of order at most $s-1$ such that every torso of $\L$ is $K_{s,t}$-topological minor free.
\end{corollary}

\begin{pf}
If some graph $G$ in $\G$ has a $(1,2s\lvert V(H) \rvert+1)$-model for a $2\lvert V(H) \rvert$-shallow $K_{s,t}$-minor, then since $s \leq 3$, every vertex of $K_{s,t}$ of degree greater than 3 corresponds to one vertex in this shallow minor, so $G$ contains a $K_{s,t}$-topological minor, a contradiction.
Since $t \geq s$, $K_{s,t}$ is $s$-connected.
Then this corollary immediately follows from Corollary \ref{minor_topo_ex}.
\end{pf}

\bigskip

Arguably the most well-known topological minor-closed but non-minor-closed class is the class of graphs with bounded crossing number.

A {\it surface} is a non-null connected 2-dimensional compact manifold without boundary.
Let $\Sigma$ be a surface.
A {\it drawing} of a graph is a mapping that maps vertices to points $\Sigma$ injectively and maps each edge to a curve in $\Sigma$ connecting its ends internally disjoint from the points embedding the vertices such that no point in $\Sigma$ belongs to the interior of at least three of those curves; a {\it crossing} of this drawing is a point that belongs to the intersection of the interior of two of those curves.
The {\it $\Sigma$-crossing number} of a graph $G$ is the minimum nonnegative integer $k$ such that $G$ has a drawing in $\Sigma$ with at most $k$ crossings.
When $\Sigma$ is the sphere, the $\Sigma$-crossing number is called the {\it crossing number}.
Note that graphs with $\Sigma$-crossing number $0$ are exactly the graphs that can be drawn in $\Sigma$ with no crossing.
It is easy to see that the class of graphs with $\Sigma$-crossing number at most a fixed integer $k$ is topological minor-closed but not minor-closed.
See \cite{s_survey} for a survey about crossing number.

Recall that Wormald \cite{w_3conn} and Eppstein \cite{e} proved that $\ex(H,\G,n)=\Theta(n)$ when $\G$ is the set of planar graphs and $H$ is a 3-connected planar graph.
Huynh, Joret and Wood \cite{hjw} generalized this result by determining the value $k$ such that $\ex(H,\G,n)=\Theta(n^k)$ when $\G$ is the class of graphs of bounded Euler genus and $H$ is an arbitrary graph in $\G$.
Corollary \ref{cor_bdd_expan} leads to a more general result: it determines $\ex(H,\G,n)$ up to a constant factor for the class of graphs with $\Sigma$-crossing number at most $t$ for any surface $\Sigma$ and nonnegative integer $t$.

For an independent collection $\L$ of separations of a graph $H$, the {\it central torso} of $\L$ is the graph obtained from $H[\bigcap_{(A,B) \in \L}B]$ by adding edges such that for every $(A,B) \in \L$, $A \cap B$ is a clique; every edge in the central torso of $\L$ whose both ends are in $A \cap B$ for some $(A,B) \in \L$ is called a {\it peripheral edge}; every torso of $\L$ that is not the central torso is called a {\it peripheral torso}.

\begin{lemma} \label{total_cr_dup}
Let $\Sigma$ be a surface.
Let $t$ be a nonnegative integer.
Let $\G$ be the class of graphs with $\Sigma$-crossing number at most $t$.
Let $\L$ be an essential independent collection of $H$ of order at most $2$.
Then $\L$ is $\G$-duplicable if and only if every peripheral torso of $\L$ is planar and the central torso of $\L$ can be drawn in $\Sigma$ with at most $t$ crossings such that every peripheral edge does not contain any crossing.
\end{lemma}

\begin{pf}
It suffices to show that if $\L$ is $\G$-duplicable, then every peripheral torso of $\L$ is planar and the central torso of $\L$ can be drawn in $\Sigma$ with at most $t$ crossings such that every peripheral edge does not contain any crossing, since the converse statement is obvious.

Assume that $\L$ is $\G$-duplicable.
Hence for every positive integer $w$, $H \wedge_w \L$ can be drawn in $\Sigma$ with at most $t$ crossings, and we denote such a drawing by $D_w$.
Since $\L$ is essential, for every $(A,B) \in \L$ with $\lvert A \cap B \rvert=2$, there exists a path $P_A$ in $H[A]$ between the two vertices in $A \cap B$ with length at least two.
Since every crossing is contained in at most two edges, by using the drawing $D_{2t+1}$, $H$ can be drawn in $\Sigma$ with at most $t$ crossings such that for every $(A,B) \in \L$ with $\lvert A \cap B \rvert=2$, no edge in $P_A$ contains an crossing, so the central torso of $\L$ can be drawn in $\Sigma$ with at most $t$ crossings such that every peripheral edge does not contain any crossing.
In addition, by using the drawing $D_{2g+2t+3}$, where $g$ is the Euler genus of $\Sigma$, we know that for every $(A,B) \in \L$, $H[A] \wedge_{2g+2} (A \cap B)$ can be drawn in $\Sigma$ with no crossing, so by \cite[Theorem 1.1]{a_add_genus} (or \cite[Theorem 1]{m_add_genus}), $H[A] \wedge_2 (A \cap B)$ can be drawn in the sphere with no crossing, and hence the peripheral torso corresponding to $(A,B)$ is planar.
\end{pf}

\begin{corollary} \label{cor_total_cr}
Let $\Sigma$ be a surface.
Let $t$ be a nonnegative integer.
Let $\G$ be the class of graphs with $\Sigma$-crossing number at most $t$.
Let $H$ be a graph in $\G$.
Then $\ex(H,\G,n) = \Theta(n^k)$, where $k$ is the maximum size of an essential independent collection of separations of $H$ of order at most $2$ whose every peripheral torso is planar and whose central torso can be drawn in $\Sigma$ with at most $t$ crossings such that every peripheral edge does not contain any crossing.
\end{corollary}

\begin{pf}
If $L$ is a graph that can drawn in $\Sigma$ with at most $t$ crossings, then adding a vertex on each crossing leads to a new graph $L'$ with $\lvert V(L') \rvert \leq \lvert V(L) \rvert+t$ and $\lvert E(L') \rvert \geq \lvert E(L) \rvert$.
So Euler's formula implies that there exists a positive integer $q$ such that $K_{3,q} \not \in \G$.
If some graph $G$ in $\G$ contains a $(1,6\lvert V(H) \rvert+1)$-model for a $2\lvert V(H) \rvert$-shallow $K_{3,q}$-minor, then $G$ contains a subdivision of $K_{3,q}$ since each of the branch sets corresponding to vertices in $K_{3,q}$ with degree greater than 3 has size 1.
Since $\G$ is topological minor-closed, no graph in $\G$ contains a $(1,6\lvert V(H) \rvert+1)$-model for a $2\lvert V(H) \rvert$-shallow $K_{3,q}$-minor.

Since every topological minor-closed class has bounded expansion, by Corollary \ref{cor_bdd_expan} and Proposition \ref{make_essential}, $\ex(H,\G,n) = \Theta(n^k)$, where $k$ is the maximum size of an essential $\G$-duplicable independent collection of separations of $H$ of order at most $2$.
Then this corollary follows from Lemma \ref{total_cr_dup}.
\end{pf}

\bigskip

We remark that the case $t=0$ implies the result of Huynh, Joret and Wood \cite{hjw} for graphs with bounded Euler genus.

\subsection{Non-topological minor-closed classes}

In this subsection, we give examples for non-topological minor-closed classes but still with bounded expansion.

\subsubsection{With some subgraphs forbidden}

It is natural to study $\ex(H,\G,n)$, where $\G$ is a class of planar graphs with no subgraph isomorphic to certain fixed graph.
We can make it more general: what is $\ex(H,\G,n)$, where $\G$ is a class of graphs satisfying a topological minor-closed property and with no subgraph isomorphic to any member in a fixed set of graphs?
Such as a class is not topological minor-closed but still has bounded expansion.
Corollary \ref{cor_bdd_expan} can deal with such classes.
One example of such a question is the following conjecture proposed by Gy\H{o}ri, Paulos, Salia, Tompkins and Zamora \cite{gpstz_turan}, and they proved the case for $\ell=2$ and $p \geq 5$.

\begin{conj}[{\cite[Conjecture 2.10]{gpstz_turan}}] \label{conj_planar_subgraph}
For every integer $\ell \geq 2$ and for every sufficiently large integer $p$, if $H$ is the $p$-cycle and $\G$ is the class of all planar graphs with no even cycle of length between $4$ and $2\ell$, then $\ex(H,\G,n)=\Theta(n^{\lfloor \frac{p}{\ell+1} \rfloor})$.
\end{conj}

Here we provide the following example (Corollary \ref{cor_cycle}) to show how to use Corollary \ref{cor_bdd_expan} to solve questions in this kind; Conjecture \ref{conj_planar_subgraph} follows from its special case that $\Sigma$ is the sphere, $t=0$, and $H$ is a $p$-cycle, since it is easy to construct a collection of separations of the $p$-cycle with size $\lfloor \frac{p}{\ell+1} \rfloor$ satisfying the condition mentioned in Corollary \ref{cor_cycle} by considering edge-disjoint paths of length $\ell+1$; the special case that $\Sigma$ is the sphere, $t=0$ and $H$ is a tree was also proved in \cite{gpstz_turan}.

\begin{corollary} \label{cor_cycle}
Let $\Sigma$ be a surface.
Let $t$ be a nonnegative integer.
Let $\ell$ be a positive integer.
Let $\G$ be the class of all graphs with $\Sigma$-crossing number at most $t$ and with no even cycle of length between $4$ and $2\ell$.
Let $H$ be a graph in $\G$.
Then $\ex(H,\G,n)=\Theta(n^k)$, where $k$ is the maximum size of an essential independent collection $\L$ of separations of $H$ of order at most $2$ such that 
	\begin{enumerate}
		\item every peripheral torso is planar, 
		\item the central torso can be drawn in $\Sigma$ with at most $t$ crossings such that no peripheral edge contains a crossing, and 
		\item for every $(A,B) \in \L$ with $\lvert A \cap B \rvert=2$, there exists no path in $H[A]$ between the two vertices in $A \cap B$ with length between $2$ and $\ell$.
	\end{enumerate}
\end{corollary}

\begin{pf}
Let $\G'$ be the class of graphs with $\Sigma$-crossing number at most $t$.
Since $\G'$ is topological minor-closed, $\G'$ has bounded expansion.
Since $\G$ is a subset of $\G'$, for any nonnegative integer $r$, any $r$-shallow minor of a member of $\G$ is an $r$-shallow minor of a member of $\G'$.
So $\G$ has bounded expansion.
And as shown in the proof of Corollary \ref{cor_total_cr}, no graph in $\G'$ (and hence in $\G$) contains a $(1,6\lvert V(H) \rvert+1)$-model for a $2\lvert V(H) \rvert$-shallow $K_{3,q}$-minor for some integer $q$.
So by Corollary \ref{cor_bdd_expan} and Proposition \ref{make_essential}, $\ex(H,\G,n) = \Theta(n^k)$, where $k$ is the maximum size of an essential $\G$-duplicable independent collection of separations of $H$ of order at most $2$.
Since every $\G$-duplicable collection is $\G'$-duplicable, by Lemma \ref{total_cr_dup}, every essential $\G$-duplicable collection satisfies Statements 1 and 2. 

In addition, if $(A,B)$ is a member of a $\G$-duplicable independent collection of separations of $H$ with $\lvert A \cap B \rvert=2$, then there is no path in $H[A]$ between the two vertices in $A \cap B$ with length between $2$ and $\ell$, for otherwise $H[A] \wedge_2 (A \cap B)$ contains an even cycle of length between $4$ and $\ell$.
Hence every essential independent collection of separations of $H$ of order at most $2$ satisfies Statements 1-3. 
On the other hand, it is easy to see that every essential independent collection of separations of order at most $2$ satisfying Statements 1-3 is $\G$-duplicable. 
This proves the corollary.
\end{pf}

\subsubsection{Linearly many crossings}

For a surface $\Sigma$ and a nonnegative integer $t$, a graph is {\it $(\Sigma,t)$-planar} if it can be drawn in $\Sigma$ such that every edge contains at most $t$ crossings.
When $\Sigma$ is the sphere, $(\Sigma,t)$-planar graphs are called {\it $t$-planar} graphs in the literature.
The case $t=1$ is of particular interests in the literature.
See \cite{klm} for a survey.
$(\Sigma,t)$-planar graphs are almost equivalent to {\it $(g,t)$-planar graphs} which are the graphs that can be drawn in a surface of Euler genus at most $g$ such that every edge contains at most $t$ crossings.
Note that every graph with $\Sigma$-crossing number at most $t$ is $(\Sigma,t)$-planar.
But in general, the number of crossings of a $(\Sigma,t)$-planar graph can be linear in the number of its edges.

For a surface $\Sigma$ and a nonnegative real number $t$, a graph $G$ is {\it $t$-close to $\Sigma$} if for every subgraph $L$ of $G$, $L$ can be drawn in $\Sigma$ with at most $t\lvert E(L) \rvert$ crossings.
So $(\Sigma,t)$-planar graphs are $t$-close to $\Sigma$.
Ne\v{s}et\v{r}il, Ossona de Mendez and Wood \cite{now} showed that the class of $(\Sigma,t)$-planar graphs has bounded expansion when $\Sigma$ is the sphere; Dujmovi\'{c}, Eppstein and Wood \cite{dew} further showed that for every surface $\Sigma$, the class of $(\Sigma,t)$-planar graphs has bounded layered tree-width and hence has bounded expansion.

In fact, the class of graphs that are $t$-close to $\Sigma$ also has bounded expansion.
This fact seems not appeared in the literature but might be well-known by the community.
Here we include a proof of this fact for completeness; our proof is a slight modification of the proof of a result in \cite{oow}.

\begin{prop} \label{avg_cr_bdd_expan}
For any surface $\Sigma$ and nonnegative real number $t$, the class of graphs that are $t$-close to $\Sigma$ has bounded expansion.
\end{prop}

\begin{pf}
By \cite[Theorem 11]{d}, it suffices to show that there exists $f: {\mathbb N} \cup \{0\} \rightarrow {\mathbb R}$ such that if $L$ is a graph such that some graph $G$ that is $t$-close to $\Sigma$ contains a subgraph $H$ isomorphic to an $([\ell] \cup \{0\})$-subdivision of $L$ for some nonnegative integer $\ell$, then the average degree of $L$ is at most $f(\ell)$.
By definition, $H$ is $t$-close to $\Sigma$, so $H$ can be drawn in $\Sigma$ with at most $t\lvert E(H) \rvert$ crossings.
Since $H$ is an $([\ell] \cup \{0\})$-subdivision of $L$, $L$ can be drawn in $\Sigma$ with at most $t\lvert E(H) \rvert$ edge-corssings, and $\lvert E(H) \rvert \leq (\ell+1)\lvert E(L) \rvert$.
So this drawing of $L$ has at most $t(\ell+1)\lvert E(L) \rvert$ crossings.
By \cite[Lemma 4.6]{oow}, $\lvert E(L) \rvert \leq c_{t,\ell,\Sigma}\lvert V(L) \rvert$ for some number $c_{t,\ell,\Sigma}$ only dependent on $t$, $\ell$ and $\Sigma$.
So we are done by defining $f(\ell)=2c_{t,\ell,\Sigma}$ for every nonnegative integer $\ell$.
\end{pf}

\begin{lemma} \label{per_cr_dup}
Let $\Sigma$ be a surface.
Let $t$ be a nonnegative integer.
Let $\G$ be the class of $(\Sigma,t)$-planar graphs.
Let $\L$ be an essential independent collection of $H$ of order at most $2$.
Then $\L$ is $\G$-duplicable if and only if 
	\begin{enumerate}
		\item the central torso of $\L$ can be drawn in $\Sigma$ such that every edge contains at most $t$ crossings, and no peripheral edge contains a crossing, and 
		\item for every $(A,B) \in \L$, the peripheral torso of $\L$ corresponding to $(A,B)$ can be drawn in the sphere such that every edge contains at most $t$ crossings, and if $\lvert A \cap B \rvert=2$, then the edge with both ends in $A \cap B$ does not contain any crossing.
	\end{enumerate}
\end{lemma}

\begin{pf}
It is easy to see that if $\L$ satisfies Statements 1 and 2, then $\L$ is $\G$-duplicable.
So we may assume that $\L$ is $\G$-duplicable and show that $\L$ satisfies Statements 1 and 2.

Since $\L$ is $\G$-duplicable, $H \wedge_{t\lvert \L \rvert\lvert E(H) \rvert(2g+3)+t\lvert E(H) \rvert} \L$ can be drawn in $\Sigma$ such that every edge contains at most $t$ crossings.
Since every edge contains at most $t$ crossings, $H \wedge_{t\lvert \L \rvert \lvert E(H) \rvert(2g+3)} \L$ can be drawn in $\Sigma$ such that no crossing is contained in one edge with both ends in $\bigcap_{(X,Y) \in \L}Y$ and one edge with at least one end not in $\bigcap_{(X,Y) \in \L}Y$.
For each $(A,B) \in \L$ and each copy of $A$ in $H \wedge_{t\lvert \L \rvert \lvert E(H) \rvert(2g+3)} \L$, there are at most $t\lvert E(H) \rvert$ crossings contained in an edge with both ends in this copy of $A$.
So by greedily sacrificing some copies of $A$ for $(A,B) \in \L$, we know that $H \wedge_{2g+2} \L$ can be drawn in $\Sigma$ such that every crossing is contained in two edges with both ends in $\bigcap_{(X,Y) \in \L}Y$ or in two edges with both ends in a copy of $A$ for some $(A,B) \in \L$.
Since $\L$ is essential, in this drawing of $H \wedge_{2g+2} \L$, for each $(A,B) \in \L$ with $\lvert A \cap B \rvert=2$, there exists a curve contained in the drawing of a copy of $A$ connecting the two vertices in $A \cap B$.
So the central torso of $\L$ can be drawn in $\Sigma$ such that every edge contains at most $t$ crossings, and no peripheral edge contains a crossing.
Hence Statement 1 holds.

In addition, for each $(A,B) \in \L$, by using the aforementioned drawing of $H \wedge_{2g+2} \L$, we know $H[A] \wedge_{2g+2} (A \cap B)$ can be drawn in $\Sigma$ such that every crossing is contained in two edges in the same copy of $H[A]$ but not contained in any edge with both ends in $A \cap B$.
By adding a vertex on each crossing of this drawing of $H[A] \wedge_{2g+2} (A \cap B)$, we obtain a graph $W$ that can be draw in $\Sigma$ with no crossing, and every component of $W-(A \cap B)$ is contained in the subdrawing of one copy of $H[A]-(A \cap B)$.
So by \cite[Theorem 1.1]{a_add_genus} (or \cite[Theorem 1]{m_add_genus}), the subgraph of $W$ induced by $H[A] \wedge_2 (A \cap B)$ and the vertices corresponding to their crossings can be drawn in the sphere with no crossing.
Hence the peripheral torso corresponding to $(A,B)$ can be drawn in the sphere such that every edge contains at most $t$ crossings, and if $\lvert A \cap B \rvert =2$, the edge with both ends in $A \cap B$ does not contain any crossing.
So Statement 2 holds.
\end{pf}

\begin{corollary} \label{cor_per_cr}
Let $\Sigma$ be a surface.
Let $t$ be a nonnegative integer.
Let $\G$ be the class of $(\Sigma,t)$-planar graphs.
Let $H$ be a graph in $\G$.
Then $\ex(H,\G,n) = \Theta(n^k)$, where $k$ is the maximum size of an essential independent collection of separations of $H$ of order at most $2$ such that
	\begin{enumerate}
		\item the central torso of $\L$ can be drawn in $\Sigma$ such that every edge contains at most $t$ crossings, and no peripheral edge contains a crossing, and 
		\item for every $(A,B) \in \L$, the peripheral torso of $\L$ corresponding to $(A,B)$ can be drawn in the sphere such that every edge contains at most $t$ crossings, and if $\lvert A \cap B \rvert=2$, then the edge with both ends in $A \cap B$ does not contain any crossing.
	\end{enumerate}
\end{corollary}

\begin{pf}
Let $q$ be an integer such that some graph in $\G$ has a $(1,6\lvert V(H) \rvert+1)$-model for a $2\lvert V(H) \rvert$-shallow $K_{3,q}$-minor.
Then some $([6\lvert V(H) \rvert] \cup \{0\})$-subdivision $L$ of $K_{3,q}$ belongs to $\G$ since each of the branch sets corresponding to vertices in $K_{3,q}$ with degree greater than 3 has size 1.
So $L$ can be drawn in $\Sigma$ with at most $t\lvert E(L) \rvert$ crossings.
Since $L$ is a $([6\lvert V(H) \rvert] \cup \{0\})$-subdivision of $K_{3,q}$, $K_{3,q}$ can be drawn in $\Sigma$ with at most $t\lvert E(L) \rvert \leq t(6\lvert V(H) \rvert+1)\lvert E(K_{3,q}) \rvert$ crossings.
So $K_{3,q}$ is $t(6\lvert V(H) \rvert+1)$-close to $\Sigma$.
By \cite[Lemma 4.8]{oow}, $q$ is upper bounded by a constant only depended on $\Sigma$, $t$ and $\lvert V(H) \rvert$. 
That is, there exists a large integer $q^*$ such that no graph in $\G$ contains a $(1,6\lvert V(H) \rvert+1)$-model for a $2\lvert V(H) \rvert$-shallow $K_{3,q^*}$-minor.
By Proposition \ref{avg_cr_bdd_expan}, Corollary \ref{cor_bdd_expan} and Proposition \ref{make_essential}, $k$ equals the maximum size of an essential $\G$-duplicable independent collection of separations of $H$ of order at most $2$.
Then this corollary follows from Lemma \ref{per_cr_dup}.
\end{pf}

\bigskip

Proposition \ref{avg_cr_bdd_expan}, Corollary \ref{cor_bdd_expan} and Proposition \ref{make_essential} also immediately lead to the following corollary for the class of graphs $t$-close to $\Sigma$; but the explicit descriptions for $\G$-duplicable collection is more complicated, so we omit it.

\begin{corollary} \label{cor_avg_cr}
Let $\Sigma$ be a surface.
Let $t$ be a nonnegative real number.
Let $\G$ be the class of graphs that are $t$-close to $\Sigma$.
Let $H$ be a graph in $\G$.
Then $\ex(H,\G,n) = \Theta(n^k)$, where $k$ is the maximum size of an essential $\G$-duplicable independent collection of separations of $H$ of order at most $2$. 
\end{corollary}

The case $t=0$ in either Corollaries \ref{cor_per_cr} or \ref{cor_avg_cr} coincides the result for the class of graphs of bounded Euler genus proved in \cite{hjw}.
And the special case that $H$ is a tree in either Corollaries \ref{cor_per_cr} or \ref{cor_avg_cr} was proved in \cite{hw}.

\subsubsection{Higher dimensional objects} \label{subsubsec:stack_queue}

For a positive integer $k$, a graph $G$ has a {\it $k$-stack layout} if there exist an ordering $\sigma$ of $G$ and a partition $\{E_1,E_2,...,E_k\}$ of $E(G)$ into $k$ (not necessarily non-empty) sets such that for every $i \in [k]$, there do not exist edges $u_1v_1,u_2v_2 \in E_i$ such that $\sigma(u_1)<\sigma(u_2)<\sigma(v_1)<\sigma(v_2)$.
The {\it stack number} of a graph $G$ is the minimum $k$ such that $G$ has a $k$-stack layout.
Each stack layout is also called a {\it book embedding} which is an embedding of $G$ in a union of half-planes in ${\mathbb R^3}$ sharing a common line and all vertices are embedded in this line.
And stack number is also called {\it book thickness} and {\it page number}.
Stack number and book embeddings have been extensively studied in theoretical computer science and graph drawing. (See \cite{dw} for a survey.)
It is not hard to see that graphs with stack number 1 are exactly the outerplanar graphs (for example, see \cite{bk}). 
So graphs with bounded stack number can be viewed as generalizations of outerplanar graphs.
In addition, even though classes of graphs with bounded stack numbers are not topological minor-closed families, Blankenship \cite{b} proved that every graph in a fixed minor-closed family has bounded stack number.

As for extremal aspects, Eppstein \cite{e} asks the following question. 

\begin{question}[\cite{e}]
What are the graphs $H$ for which $\ex(H,\G,n)=O(n)$ when $\G$ is the class of graphs with stack number at most $k$?
\end{question}

Corollary \ref{cor_bdd_expan} gives an answer to this question.

\begin{corollary} \label{cor_stack}
Let $t$ be a positive integer.
Let $\G$ be the class of graphs with stack number at most $t$.
Then for every graph $H$ with stack number at most $t$, $\ex(H,\G,n)=\Theta(n^k)$, where $k$ is the maximum size of an essential $\G$-duplicable independent collection of separations of $H$.
\end{corollary}

\begin{pf}
Clearly, $\G$ is monotone.
By \cite[Theorem 8.4]{now}, $\G$ has bounded expansion.
So this corollary immediately follows from Corollary \ref{cor_bdd_expan} and Proposition \ref{make_essential}.
\end{pf}

\bigskip

Stacks are one of the most fundamental data structures.
Another fundamental structure is a queue.
For a positive integer $k$, a graph $G$ has a {\it $k$-queue layout} if there exist an ordering $\sigma$ of $G$ and a partition $\{E_1,E_2,...,E_k\}$ of $E(G)$ into $k$ (not necessarily non-empty) sets such that for every $i \in [k]$, there do not exist edges $u_1v_1,u_2v_2 \in E_i$ such that $\sigma(u_1)<\sigma(u_2)<\sigma(v_2)<\sigma(v_1)$.
The {\it queue number} of a graph $G$ is the minimum $k$ such that $G$ has a $k$-queue layout.
See \cite{dw} for a survey about queue numbers and queue layouts.

Queue number and stack number share similar features.
For example, classes of graphs with bounded queue number are not (topological) minor-closed, and Dujmovi\'{c}, Joret, Micek, Morin, Ueckerdt and Wood \cite{djmmuw} proved that every proper minor-closed class has bounded queue number.
However, even though they look similar, the connection between these two notions are unclear.
It was recently shown that stack number cannot be upper bounded by queue number \cite{dehmw}.
But it remains open whether queue number can be upper bounded by stack number or not.

We also have the analogy of Corollary \ref{cor_stack} for queue number.

\begin{corollary} \label{cor_queue}
Let $t$ be a positive integer.
Let $\G$ be the class of graphs with queue number at most $t$.
Then for every graph $H$ with queue number at most $t$, $\ex(H,\G,n)=\Theta(n^k)$, where $k$ is the maximum size of an essential $\G$-duplicable independent collection of separations of $H$.
\end{corollary}

\begin{pf}
Clearly, $\G$ is monotone.
By \cite[Theorem 7.4]{now}, $\G$ has bounded expansion.
So this corollary immediately follows from Corollary \ref{cor_bdd_expan} and Proposition \ref{make_essential}.
\end{pf}

\bigskip

Classes of graphs with strongly sublinear separators are typical examples of classes of graphs with bounded expansion besides topological minor-closed classes.

For a function $f: {\mathbb N} \cup \{0\} \rightarrow {\mathbb R}$, we say that a graph $G$ admits {\it $f$-balanced separators} if for every subgraph $H$ of $G$, there exists a separation $(A,B)$ of $H$ with $\lvert A \cap B \rvert \leq f(\lvert V(H) \rvert)$ such that $\lvert A-B \rvert \leq \frac{2}{3}\lvert V(H) \rvert$ and $\lvert B-A \rvert \leq \frac{2}{3}\lvert V(H) \rvert$.
Note that for every subgraph $H$, every separation $(A,B)$ of $G[V(H)]$ is a separation of $H$, so the above condition can be simplified by only requiring every induced subgraph of $G$ having a desired separation.
A class $\G$ of graphs admits {\it strongly sublinear balanced separators} if there exist constants $c \geq 1$ and $0 < \delta \leq 1$ and a function $f: {\mathbb N} \cup \{0\} \rightarrow {\mathbb R}$ defined by $f(x)=cx^{1-\delta}$ for every $x$, such that every graph in $\G$ admits $f$-balanced separators.

Strongly sublinear separators are useful for designing divide-and-conquer algorithms and inductive proofs.
A celebrated theorem of Alon, Seymour and Thomas \cite{ast} states that every minor-closed class admits strongly sublinear separators.
On the other hand, Dvo\v{r}\'{a}k and Norin \cite{dn} proved that every class of graphs admitting strongly sublinear separators have bounded expansion.
So Corollary \ref{cor_bdd_expan} can be applied to all hereditary classes of graphs admitting strongly sublinear separators.
Here we just include one example for such classes. 

The {\it intersection graph} of a set $\C$ is a graph $G$ such that there exists a bijection $\iota$ from $V(G)$ to $\C$ such that for any distinct $u,v \in V(G)$, $uv \in E(G)$ if and only if $\iota(u) \cap \iota(v) \neq \emptyset$.
For positive integers $k$ and $d$, {\it a $k$-ply neighborhood system in ${\mathbb R}^d$} is a set of closed balls in ${\mathbb R}^d$ such that no point in ${\mathbb R}^d$ is contained in the interior of $k$ members of this set. 
Miller, Teng, Thurston and Vavasis \cite{mttv} proved that the class of intersection graphs of $k$-ply neighborhood systems in ${\mathbb R}^d$ admits strongly sublinear balanced separators.
Hence such as a class is hereditary and has bounded expansion.
Therefore, Corollary \ref{cor_bdd_expan} immediately implies the following.

\begin{corollary}
Let $k$ and $d$ be positive integers.
Let $\G$ be the class of intersection graphs of $k$-ply neighborhood systems in ${\mathbb R}^d$. 
Let $H$ be a graph in $\G$.
Then $\ex(H,\G,n) = \Theta(n^k)$, where $k$ is the maximum size of a $\G$-duplicable independent collection of separations of $H$.
\end{corollary}

Note that the $k$-ply neighborhood system is a special case of a much more general result of Dvo\v{r}\'{a}k, McCarty and Norin \cite{dmn} about intersection graphs of certain sets ensuring strongly sublinear balanced separators. 
We refer readers to \cite{dmn} for details.

\section{Concluding remarks}

In this paper, we determine the maximum number of homomorphisms in $\HH$ to an $n$-vertex labelled graph in $\G$ up to a multiplicative constant, where $\G$ is any hereditary class with bounded expansion and $\HH$ is a consistent set of homomorphisms from a fixed labelled graph $H$ to $\G$ (Theorem \ref{bdd_expan_intro} and Corollary \ref{cor_bdd_expan}).
Classes with bounded expansion are general classes of sparse graphs, including many graph classes with some geometric and combinatorial properties.
Our results are generalizations of a number of known results and solve a number of open questions, as shown in Sections \ref{sec:intro} and \ref{sec:concrete}, about counting the number of $H$-subgraphs, $H$-induced subgraphs, and homomorphisms from $H$.
Moreover, our machinery allows us to determine the exact value of the asymptotic logarithmic density for $\HH$ in nowhere dense classes $\G$ which are the most general classes of sparse graphs in sparsity theory (Theorem \ref{nowhere_dense_intro} and Corollary \ref{nowhere_dense}) and determine the maximum number of $H$-subgraphs in the classes of $d$-degenerate graphs with any fixed $d$ up to a multiplicative constant (Theorem \ref{bdd_degen_intro}).

For most of the applications of the aforementioned results stated in this paper, we only require that $\G$ is monotone and only consider the number of unlabelled $H$-subgraphs.
The proof of our main results can be significantly simplified in those settings. 
But our focus is on proving more general results.

In addition, we put no efforts on optimizing the hidden constants in our theorems.
On the other hand, as our result determines the number of $H$-subgraphs and the number of homomorphisms from $H$ up to a multiplicative constant, now we know how to normalize them to define the $H$-subgraph density or the density for homomorphisms from $H$, for any fixed graph $H$.
This immediately leads to the question whether one can develop the theory about flag algebra or homomorphism inequalities for graphs in bounded expansion classes (or even nowhere dense classes) to compute the number of $H$-subgraphs up to lower order terms.
More detailed discussions can be found in Section 8 in \cite{hjw}.

\end{document}